\documentclass[numbook,envcountsame,envcountreset]{svjour3}
\usepackage{fix-cm}

\usepackage[pass]{geometry}

\usepackage{mathptmx}  
\usepackage[T1]{fontenc}
\usepackage{textcomp}
\usepackage{amsmath,amssymb,enumerate}
\usepackage{graphics}
\usepackage{rotating}
\usepackage{afterpage,pdflscape}
\usepackage{epsfig} 
\usepackage[colorlinks=true,pdfstartview=FitH,pdfpagemode=UseNone]{hyperref}
\hypersetup{urlcolor=blue, citecolor=red}

\smartqed
\newcommand{\bpf}{\begin{proof}}
\newcommand{\epf}{\hfill \qed \end{proof}}

\renewcommand{\epsilon}{\varepsilon}
\renewcommand{\phi}{\varphi}
\renewcommand{\leq}{\leqslant}
\renewcommand{\geq}{\geqslant}

\newcommand{\eq}[1]{(\ref{#1})}
\newcommand{\be}{\begin{equation}}
\newcommand{\ee}{\end{equation}}

\newcommand\udot{\dot{u}}

\newcommand\cO{\mathcal{O}}
\newcommand\N{\mathbb{N}}
\newcommand\R{\mathbb{R}}
\newcommand\Z{\mathbb{Z}}

\newcommand{\union}{\cup}
\newcommand{\leappr}{\kern0.5ex\vcenter{\hbox{$\scriptstyle\lessapprox$}}\kern0.5ex}

\journalname{J Dyn Diff Equat}

\begin{document}

\title{Periodic Solutions of a Singularly Perturbed Delay Differential Equation With Two State-Dependent Delays}

\titlerunning{Periodic Solutions of a Singularly Perturbed State-Dependent DDE}

\author{A.R. Humphries \and D.A. Bernucci \and R. Calleja \and \\ N. Homayounfar \and M. Snarski}

\institute{A.R. Humphries \and D.A. Bernucci  \and N. Homayounfar \and M. Snarski
       \at Department of Mathematics and Statistics, McGill University, Montreal, Quebec H3A 0B9, Canada.\\
       \email{tony.humphries@mcgill.ca}
   \and D.A. Bernucci
        \at \emph{Present address:}
        School of Mathematics, Georgia Institute of Technology, Atlanta, GA 30332-0160 USA.\\
        \email{dbernucci3@math.gatech.edu}
   \and R. Calleja
        \at
         Depto.\ Matem\'aticas y Mec\'anica, IIMAS, Universidad Nacional Aut\'onoma
         de M\'exico, 01000 M\'exico. \\
        \email{calleja@mym.iimas.unam.mx}
   \and N. Homayounfar
       \at \emph{Present address:}
         Department of Statistics, University of Toronto, Toronto, Ontario M5S 3G3 Canada.\\
        \email{namdar.homayounfar@mail.utoronto.ca}
   \and M. Snarski
         \at \emph{Present address:}
         Division of Applied Mathematics, Brown University, Providence, RI 02912 USA.\\
         \email{Michael\_Snarski@Brown.edu}}

\date{\today}    

\maketitle

\begin{abstract}
Periodic orbits and associated bifurcations of
singularly perturbed state-dependent delay differential equations (DDEs) are studied when the
profiles of the periodic orbits contain jump discontinuities in the singular limit. A definition
of singular solution is introduced which is based on a continuous parametrisation of the
possibly discontinuous limiting solution. This reduces
the construction of the limiting profiles to an algebraic problem. A model
two state-dependent delay differential equation is studied in detail and periodic singular solutions
are constructed with one and two local maxima per period. A complete characterisation of the conditions
on the parameters for these singular solutions to exist facilitates an investigation of bifurcation structures
in the singular case revealing folds and possible cusp bifurcations.
Sophisticated boundary value techniques are used to numerically compute the bifurcation
diagram of the state-dependent DDE when the perturbation parameter is close to zero.
This confirms that the solutions and bifurcations constructed in the singular case persist when the perturbation parameter is nonzero, and hence demonstrates that the solutions constructed using our singular
solution definition are useful and relevant to the singularly perturbed problem.
Fold and cusp bifurcations are found very close to the parameter values predicted by the singular solution theory,
and we also find period-doubling bifurcations as well as periodic orbits with more than two local maxima per period, and
explain the alignment between the folds on different bifurcation branches.

\end{abstract}

\section{Introduction}
\label{secintro}

We consider singularly perturbed periodic solutions of the scalar
state-dependent delay-differential equation (DDE)
\be \label{eps2del}
\epsilon\udot(t)=-u(t)-K_1 u(t-a_1-cu(t))-K_2u(t-a_2-cu(t)),
\ee
which has two linearly state-dependent delays, and no other nonlinearity apart from
the state-dependency of the delays. We consider
$\epsilon\geq0$, $c>0$, $a_i>0$, $K_i>0$, and without loss of generality
we order the terms so that $a_2>a_1>0$. Equation~\eq{eps2del} is an example of
a singularly perturbed scalar DDE
with $N$ state-dependent delays of the form
\be \label{epsf}
\epsilon\udot(t)=f(t,u(t),u(\alpha_1(t,u(t))),\ldots,u(\alpha_N(t,u(t)))), \qquad u(t)\in\R.
\ee
We will define a concept of singular solution for \eq{epsf} based on a continuous parametrisation.
This essentially entails defining a singular limit for the equation
\eq{epsf}, resulting in an equation whose solutions can in principle be found algebraically.
In the case of \eq{eps2del} we construct several such classes of singular periodic
solutions, and investigate the codimension-one and -two bifurcations that arise.

DDEs arise in many applications including engineering, economics, life sciences
and physics \cite{E09,M89,Smith10}. There is a well established theory for functional
differential equations as infinite-dimensional dynamical systems 
\cite{DGVLW95,H&L93},
which encompasses DDEs with constant or prescribed delay.
However, many problems that arise in applications have delays which depend on the state of the system
(see for example \cite{FM09,IST07,LR13,W03}). Such state-dependent DDEs fall outside of the scope of the previously
developed theory and have been the subject of much study in recent years. See \cite{HKWW06}
for a relatively recent review of the general theory of state-dependent DDEs.

The study of singularly perturbed DDEs already stretches over several decades. As early as 1985
Magalha\~es \cite{Mag85} recognised that singularly perturbed discrete DDEs
have different asymptotics to singularly perturbed distributed DDEs. For
equations with a single constant delay, in the singular limit the DDE reduces to a map
(see \eq{singmap} below) which describes the asymptotic behaviour when the limiting
profiles are functions \cite{JMPRN86,IS92,SMR93}.

\begin{figure}[t!]
\begin{center}
%
\scalebox{0.45}{\includegraphics{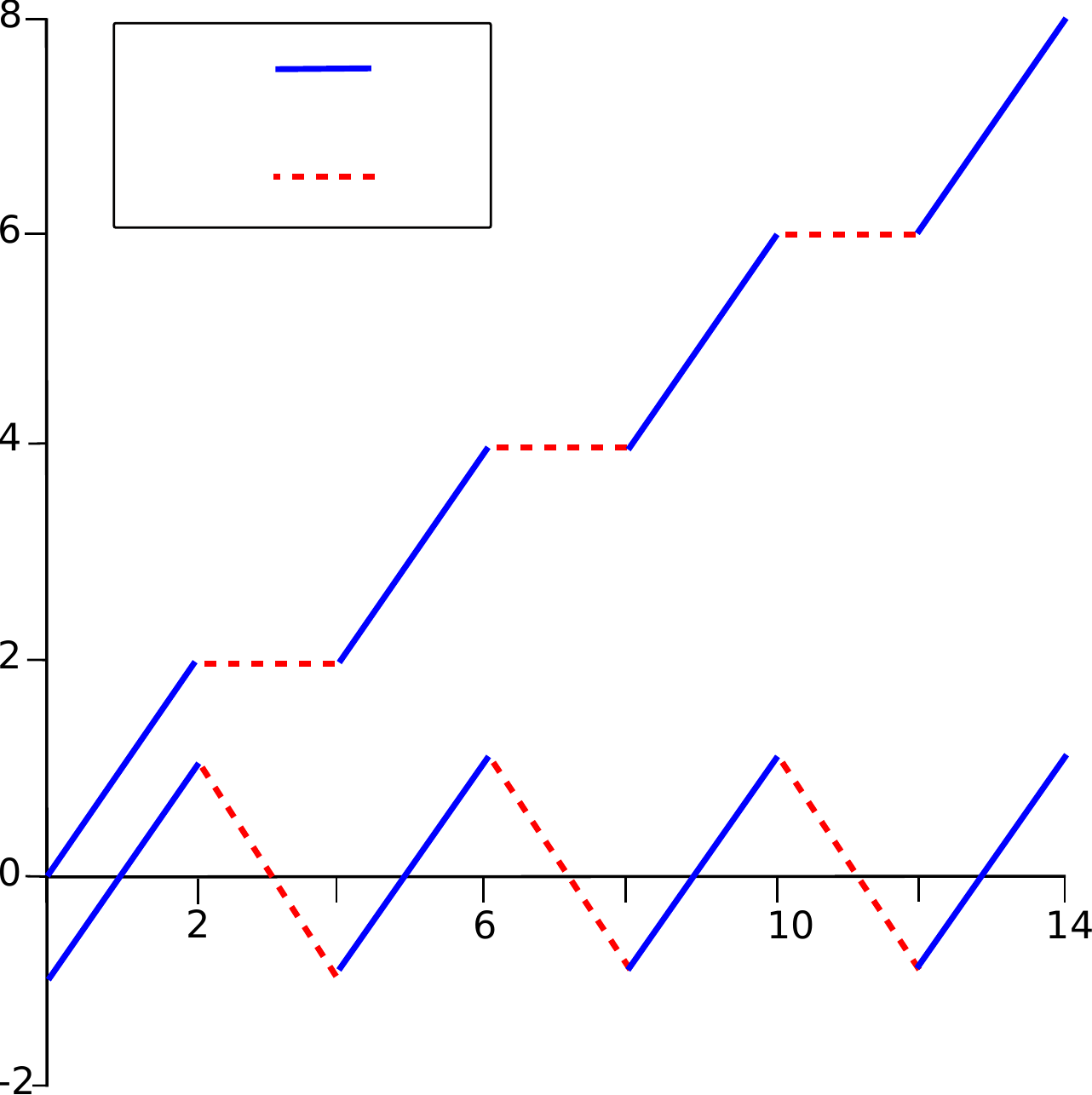}}\hspace{0mm} 
\scalebox{0.45}{\includegraphics{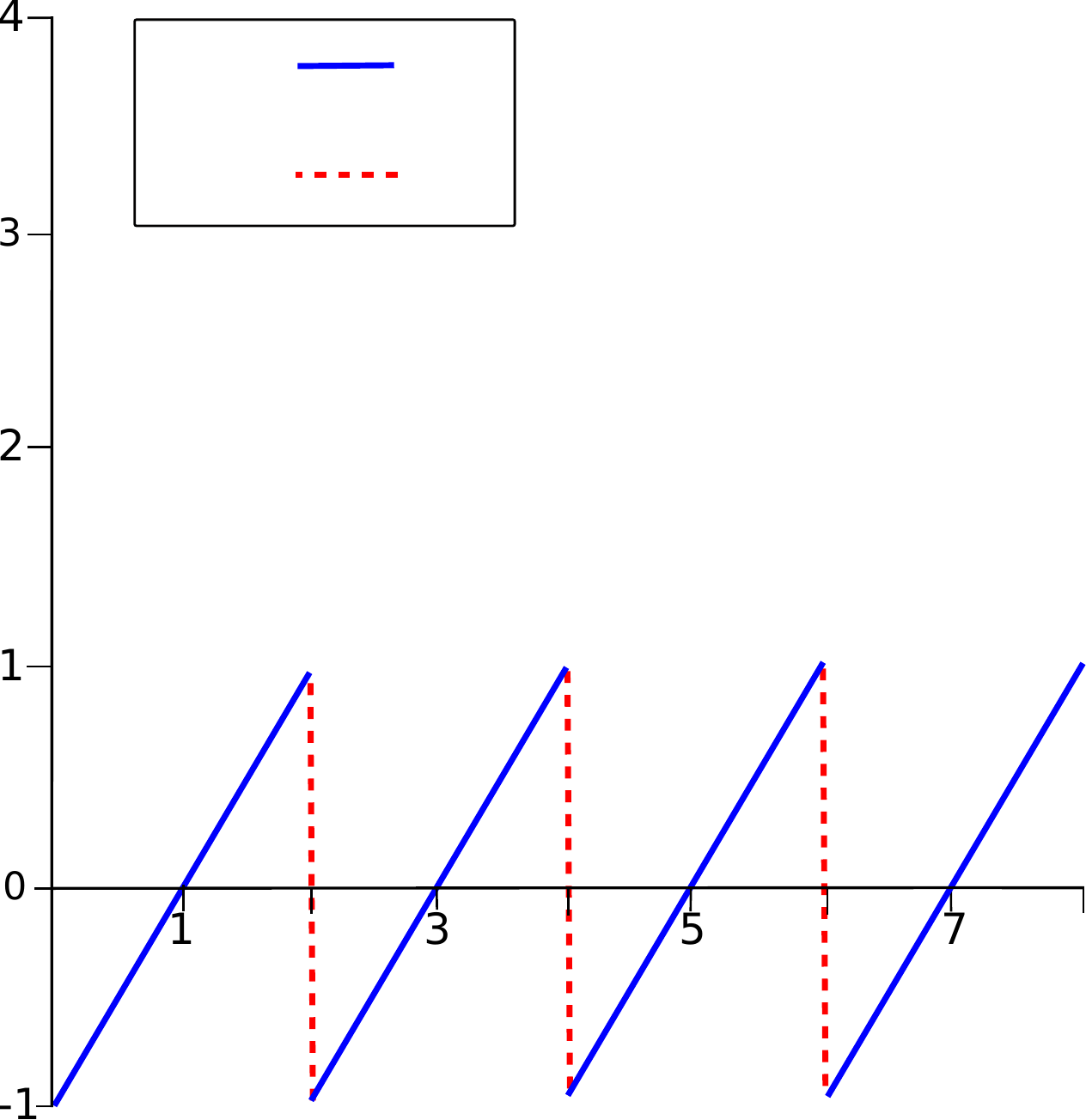}}   
\put(-178,15){$\mu$}
\put(-195,50){$u(\mu)$}
\put(-203,150){$t(\mu)$}
\put(-305,154){$I^*$}
\put(-305,138){$I^-$}
\put(-8,20){$t$}
\put(-167,120){$u$}
\put(-133,157){$\Gamma^*$}
\put(-133,141){$\Gamma^-$}
\put(-240,160){(i)}
\put(-60,160){(ii)}
\end{center}
\caption{An admissible singular solution. (i) The functions $t(\mu)$ and $u(\mu)$ for $\mu\in I=[0,14]$.
(ii) The corresponding admissible singular solution profile $\Gamma(I)=\Gamma^*\union\Gamma^-\in\R^2$.}
\label{figadmissprof}
\end{figure}
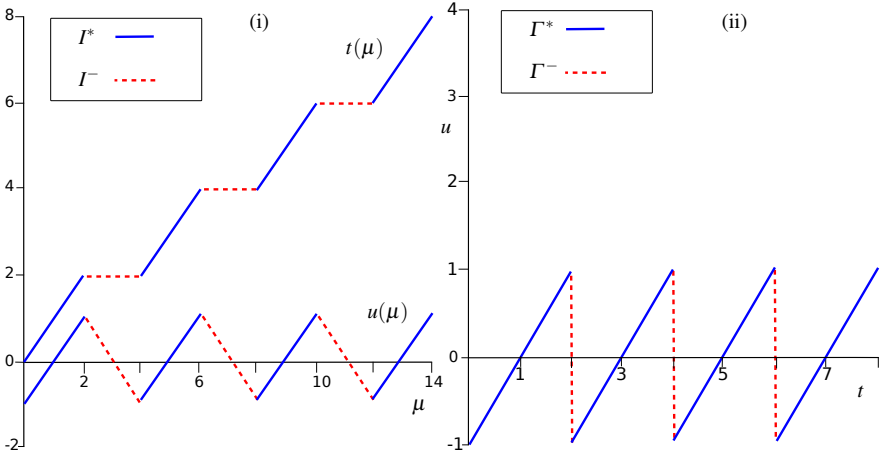

One of the main difficulties studying \eq{epsf} in the singular limit is that while
the solution $u(t)$ is a graph for any $\epsilon>0$, this need not be so in the limit as $\epsilon\to0$,
when derivatives can become unbounded, and the resulting limiting solution can have jump discontinuities.
Techniques for studying singularly perturbed DDEs with a single constant discrete delay
can be found in \cite{JMPRN86,CLMP89}. In \cite{JMPRN86} slowly oscillating periodic solutions (SOPS) are proved to
converge to a square wave in the singular limit, using layer equations to describe the solution
in the transition layer. In \cite{CLMP89} for monotone nonlinearities a homotopy method is used to show that
the layer equations have a unique homoclinic orbit. Mallet-Paret and Nussbaum,
in a series of papers \cite{JMPRNI,JMPRNII,JMPRNIII} extend the study of SOPS
to DDEs with a single state-dependent delay. In \cite{JMPRNI} SOPS
are shown to exist for all $\epsilon$ sufficiently small. These solutions are shown to have non-vanishing amplitude
in the singular limit in \cite{JMPRNII}, and under mild assumptions the discontinuity set of the limiting profile
is shown to consist of isolated points.
In \cite{JMPRNIII} Max-plus operators are introduced to study the shapes
of the limiting profiles. The DDE
\be \label{eps1del}
\epsilon\udot(t)=- u(t)-K u(t-a_1-u(t)),
\ee
is considered as an example in \cite{JMPRNIII}.
This corresponds to \eq{eps2del} with $K_2=0$ and $c=1$.
It is shown in \cite{JMPRNIII} that the limiting profile is the ``sawtooth''
shown in Fig.~\ref{figadmissprof}(ii) below.
In \cite{JMPRN11} the SOPS of \eq{eps1del} are studied in detail and the shape of the solution near the local maxima and
minima is determined for $0<\epsilon\ll1$ as well as the width of the transition layer, and the ``super-stability'' of the solution.
Other work on singularly perturbed state-dependent DDEs includes \cite{GH12}
where they
arise from the regularisation of
neutral state-dependent DDEs, and also \cite{PRMP14} where the
metastability of solutions of a singularly perturbed state-dependent DDE is studied
in the case where the state-dependency vanishes in the limit as $\epsilon\to0$.

The studies mentioned above all considered singularly perturbed DDEs with only one delay, and either considered
single solutions or a sequence of solutions as $\epsilon\to0$.
We will study the bifurcation diagram for the two-delay DDE
\eq{eps1del} when $0\leq\epsilon\ll1$, regarding $K_1$ as a bifurcation parameter.
Beyond those mentioned previously, the only other work we know of that tackles singularly perturbed
bifurcations in state-dependent DDEs is
\cite{KE14}, where the solutions of \eq{eps1del} with $a_1=c=1$ are studied close to the singular Hopf bifurcation.
On the other hand, singularly perturbed ODEs frequently arise through mixed mode oscillations on multiple time-scales and their
bifurcation analysis is well understood (see \cite{GKO12} for a review). Codimension-two bifurcations have also
been studied in singularly perturbed ODEs \cite{BKK13,Chiba11}.

The development of bifurcation theory for state-dependent DDEs has been difficult because the centre manifolds have not
been shown to have the necessary smoothness \cite{HKWW06}, and a rigorous
Hopf bifurcation theorem for state-dependent DDEs was first proved only in the
last decade \cite{E06} (see also \cite{HW10,Sieber12,GW13}). The numerical analysis of state-dependent DDEs
is more advanced with numerical techniques for solving both initial value problems \cite{BZ03,BMZG09} and
for computing bifurcation diagrams \cite{ETR02}.
DDEBiftool \cite{ETR02} is a very useful tool
for computing Hopf bifurcations and continuation of solution branches in state-dependent DDEs, and it has
been used to study the bifurcations that arise in \eq{eps2del} when $\epsilon=1$  \cite{DCDSA11,CHK15}.
John Mallet-Paret has presented numerical simulations of \eq{eps2del} in seminars, but the only other
published work of which we are aware that encompasses \eq{eps2del} is
\cite{JMPRNPP94}. There the existence of SOPs was proved for \eq{epsf} with suitable nonlinearities when
$\alpha_i(t,u(t))=t-\tau_i(u(t))$  with $\tau_i(0)=k\ne0$ for all $i$. Mallet-Paret and Nussbaum
have announced results for the existence of periodic orbits in state-dependent DDEs with two
delays including equations of the form \eq{eps2del}, but these results are as yet unpublished \cite{JMPRN15}.

In \cite{DCDSA11} a largely numerical investigation of \eq{eps2del} with $\epsilon=\cO(1)$ revealed fold
bifurcations on the branches of periodic orbits, resulting in parameter regions with bistability of
periodic orbits. While the stable periodic orbits usually had one local maxima per period,
the unstable periodic orbits in the these windows of bistability often had more than one local maxima
per period. In the current work we will investigate these fold bifurcations and the profiles
of the periodic orbits in the singular limit $\epsilon\to0$.

To study
\eq{epsf} in the singular limit $\epsilon=0$ when the limiting profile may have jump discontinuities,
we propose nested continuous parameterisations of the limiting singular solution.
We will not restrict our attention to slowly oscillating periodic orbits, but will consider both long and short period orbits.
We will study the case of the two delay state-dependent DDE
\eq{eps2del} in detail, and construct branches of singular periodic orbits
with fold and cusp bifurcations. We will then use the predictions of this theory to guide a
numerical study which will reveal branches of periodic orbits for $0<\epsilon\ll1$ with profiles
close to the singular limiting profiles and fold
and cusp bifurcations close to the predicted parameter values. We will also find period-doubling bifurcations in
the singularly perturbed problem.

Since our parametrisation technique is our main theoretical tool and crucial to all our results,
we will describe it here in detail.
For our outer parametrisation
we consider the solution profile as a parametric curve, $\Gamma(\mu)=(t(\mu),u(\mu))$.
This is a familiar concept from physics, where trajectories in space-time are parameterised, and
has been used in the study of the DDEs arising from Wheeler-Feynman Electrodynamics \cite{DLHR12}.
However, in the current work we use the parametric curve $\Gamma(\mu)$ to enable us to
consider continuous objects even in the singular limit. For any $\epsilon>0$ an injective
parametrisation of the solution must have $t(\mu)$ strictly monotonic,
but limiting profiles as $\epsilon\to0$ may have $t(\mu)$ merely monotonic.
This leads us to the parametric definition of an admissible singular solution profile in Definition~\ref{defadmiss}.
In Definition~\ref{defsingsol} we will introduce the inner parametrisation that allows us to define
singular solutions of \eq{epsf}.

\begin{definition} \label{defadmiss}
Let $\Gamma:I\to\R^2$ be a continuous injective parametric curve defined on a nonempty interval $I\subseteq\R$. For $\mu\in I$
let $\Gamma(\mu)=(t(\mu),u(\mu))$. Then if $t:I\to\R$ is monotonically increasing we say that $\Gamma(I)$ is an \emph{admissible singular
solution profile} for \eq{epsf}.
\end{definition}

Although $t(\mu)$ is not required to be a strictly monotonically increasing function
to be an admissible singular solution profile,
it is important to note that
on any subinterval $I_i$ on which $t(\mu)$ is constant,
the injectivity requirement ensures that $u(\mu)$ is strictly monotonic. Thus we partition the interval $I$ as
$I=I^+\union I^-\union I^*$ where
\begin{enumerate}
\item
$I^*$ a disjoint union of open intervals and $t(\mu)$ is strictly monotonically increasing on each interval,
\item
$I^\pm$ are each disjoint unions of closed intervals with $t(\mu)$ constant on each such interval, and $u(\mu)$
strictly monotonically decreasing (respectively increasing) on each interval of $I^-$ (resp. $I^+$).
\end{enumerate}
The partition of $I$ generates a corresponding partition of $\Gamma(I)$ as $\Gamma(I)=\Gamma^+\union \Gamma^-\union \Gamma^*$.
For~\eq{eps2del} we will find that $I^+=\emptyset$, and so $I^*$ and $I^-$ will both be unions of
disjoint intervals which we may write as
$$I^*=\bigcup_i I_{2i}=\bigcup_i(b_{2i},b_{2i+1}), \quad I^-=\bigcup_{i} I_{2i+1}=\bigcup_i[b_{2i+1},b_{2i+2}].$$
for a sequence of strictly increasing real numbers $b_i$. See Fig.~\ref{figadmissprof} for an example.

The partition of $\Gamma(I)$ as $\Gamma(I)=\Gamma^{+\!}\union \Gamma^{-\!}\union \Gamma^{*\!}$ is similar that
of $\Omega=\Omega^+\union \Omega^-\union \Omega^*$ introduced by Mallet-Paret and Nussbaum \cite{JMPRNII}
(see also Section 4 of \cite{JMPRNIII}). In their work $\Omega$ is defined as the limiting set for a sequence
of solutions as $\epsilon\to0$, while $\Omega^\pm$ are defined as the sets of points for which $\liminf\pm\epsilon\udot(t)>0$, which
results in $\Omega^{\pm\!}$ being relatively open subsets of $\Omega$. In contrast,
we define $\Gamma(I)$ and its partition directly from the parametrisation of the admissible singular solution profile,
with $\Gamma^{\pm\!}$ being closed subsets of $\Omega$. Now intuitively, since $\Gamma^{*\!}$ defines the parts of the
singular solution profile for which $\udot$ is finite, from \eq{epsf} it should correspond to the parts of
the solution for which $\lim_{\epsilon\to0}f(t,u(t),u(\alpha_1(t,u(t))),\ldots,u(\alpha_N(t,u(t))))=0$. \sloppy{Similarly
$\udot=\pm\infty$ on $\Gamma^{\pm\!}$ should imply that $\lim_{\epsilon\to0}f(t,u(t),u(\alpha_1(t,u(t))),\ldots,u(\alpha_N(t,u(t))))$
is respectively positive or negative.}
Rather than treating this process as $\epsilon\to0$ we introduce
an extra level of parametrisation, so that we can write the right-hand side of \eq{epsf} as a function
of a single parametrisation variable, which allows us to treat the $\epsilon=0$ case directly
in a continuous framework.

\begin{definition} \label{defsingsol}
Let $\Gamma$ be an admissible singular solution profile defined on $I\subseteq\R$ and let
$J\subseteq\R$ be a nonempty interval. Let $\mu_i:J\to I$ for $i=0,\ldots,N$ be continuous
functions with $\mu_0(\eta)$ monotonically increasing. Define $J^*=cl\{\eta:\mu_0(\eta)\in I^*\}$
and $J^\pm=int\{\eta:\mu_0(\eta)\in I^\pm\}$, and let
\be \label{singF}
F(\eta)=f(t(\mu_0(\eta)),u(\mu_0(\eta)),u(\mu_1(\eta)),\ldots,u(\mu_N(\eta))).
\ee
Then if
\be \label{singdel}
t(\mu_i(\eta))=\alpha_i(t(\mu_0(\eta)),u(\mu_0(\eta))),\quad\forall\eta\in J, \; \forall i=1,\ldots,N,
\ee
and
\begin{enumerate}
\item $F(\eta)=0$ for all $\eta\in J^*$,
\item $F(\eta)<0$ for all $\eta\in J^-$
\item $F(\eta)>0$ for all $\eta\in J^+$.
\end{enumerate}
we say that $\{\Gamma, \mu_0, \ldots, \mu_N\}$ define a \emph{singular solution} for \eq{epsf}
on the interval $t(\mu_0(J))$.
\end{definition}

In the definition, essentially one can think of $t(\mu_0(\eta))$ as the current time, and $t(\mu_i(\eta))$ as the delayed times.
Then \eq{singdel} simply says that the delayed times are given by the formula for $\alpha_i$ from the DDE~\eq{epsf}, while
\eq{singF} reduces the right-hand side of \eq{epsf} to a continuous function of the inner parametrisation variable.
Any solution of \eq{epsf} for $\epsilon>0$ can be similarly parameterised, resulting in
\be \label{ddeparam}
\epsilon\udot(t(\mu_0(\eta)))=F(\eta).
\ee
Now the conditions on $F(\eta)$ in the definition for a singular solution with $\epsilon=0$ follow from the remarks
on the sets $\Gamma^*$, $\Gamma^\pm$ before the definition.

This concept of singular solution generalises that of \cite{JMPRN86,IS92,SMR93}. To see this, consider the case
where equation
\eq{epsf} is autonomous with one fixed delay, so $N=1$ and $\alpha_1(t,u(t))=t-\tau$ for some constant $\tau>0$.
Suppose also that the limiting profile is a graph, so $\Gamma^-=\Gamma^+=\emptyset$.
Then we can define a singular solution following Definition~\ref{defsingsol} with
$\mu_0(\eta)=\eta$, $\mu_1(\eta)=\eta-\tau$, and $t=t(\mu)=\mu$. This parametrisation
respects \eq{singdel}, and since $\Gamma^-=\Gamma^+=\emptyset$ we have $J^*=J$ and require $F(\eta)=0$ for all $\eta\in J$.
But then
$$0=F(\eta)=f(u(\mu_0(\eta)),u(\mu_1(\eta)))=f(u(\eta),u(\eta-\tau))=f(u(t),u(t-\tau)),$$
and we are left to consider
\be \label{singmap}
f(u(t),u(t-\tau))=0,
\ee
which is the equation studied in
\cite{JMPRN86,IS92,SMR93}. Thus in the case that $\Gamma^-=\Gamma^+=\emptyset$ our definition
encompasses that of \cite{JMPRN86,IS92,SMR93}. However, in this work we will be interested in the case where
$\Gamma^-$ is not empty, and the delays are not constant.

If $J=\R$ and there exists $T>0$ and $\eta_T>0$ such that
$$t(\mu_i(\eta+\eta_T))=t(\mu_i(\eta))+T, \quad  u(\mu_i(\eta+\eta_T))=u(\mu_i(\eta)), \qquad
\forall i=1,\ldots,N, \quad \forall\eta\in\R,$$
we say that the singular solution is periodic. The period is the smallest such $T>0$. 

The main aim of this paper is to initiate a study of periodic solutions of the singularly perturbed two-delay DDE
\eq{eps2del}. We will construct singular periodic solutions
(as per Definition~\ref{defsingsol}), and will find both unimodal
sawtooth solutions that correspond to the profile seen in Fig.~\ref{figadmissprof}
and bimodal solutions which have two ``teeth'' per period. The labels {\em unimodal} and {\em bimodal} are
used throughout to indicate the number of local maxima of the solution per period.
Although
superficially the unimodal solutions look similar to those found in the one delay case,
the interaction between the two state-dependent delays adds both complications to
the derivations and richness to the dynamics observed. We will demonstrate numerically using
DDEBiftool \cite{ETR02}, a sophisticated numerical bifurcation package
for DDEs, that the singular solutions and associated
bifurcation structures that we find persist for $\epsilon>0$.

In Section~\ref{secsingsols} as an example we first consider \eq{eps1del} with one delay,
for which Mallet-Paret and Nussbaum \cite{JMPRNIII,JMPRN11} have already established the so-called
{\em sawtooth} limiting profile, as illustrated in Fig.~\ref{figadmissprof}(ii). We construct the
corresponding singular solution following Definition~\ref{defsingsol}. We then consider the two-delay problem
\eq{eps2del} and in Theorem~\ref{thm2delunimod} establish conditions on the parameters for this to have
a sawtooth solution. In \eq{eqtypeI} and \eq{eqtypeII} we introduce two admissible singular solution profiles
which have two local maxima per period. Theorems~\ref{thm2delI} and~\ref{thm2delII} present singular solutions
for these profiles and establish the constraints on the parameters for them to exist. Since these solutions
have two local maxima per period, we refer to them as type I and type II bimodal (periodic) solutions.

In Section~\ref{secbifs} we treat $K_1$ as a bifurcation parameter and in Theorems~\ref{thmunilegs},
\ref{thmIbif} and~\ref{thmIIbif} identify intervals of the parameter $K_1$ for which unimodal, type I
bimodal and type II bimodal solutions exist. We will also find singular fold bifurcations in
Theorem~\ref{thmIbif} where solutions transition between unimodal and type I bimodal solutions.
Theorem~\ref{thmIIbif} as well as identifying a singular fold bifurcation between the unimodal and type II bimodal solutions also
identifies a curve of parameter values at which a codimension-two singular cusp bifurcation occurs.
The fold bifurcation unfolds at this bifurcation and there is a transition between unimodal and type II bimodal solutions without a fold
in the bifurcation branch.

The definition of singular solution introduced above, and the resulting solutions found
are only useful if they tell us something about the dynamics of \eq{eps2del} when $0<\epsilon\ll1$.
In the case of one delay \eq{eps1del}, Mallet-Paret and Nussbaum \cite{JMPRNIII} proved the existence for $\epsilon>0$ of
a singular solution which is a perturbation of the sawtooth profile. It is not readily apparent
how to extend that proof to the two delay DDE \eq{eps2del}. So in Section~\ref{secnumerics}
we perform a numerical investigation of \eq{eps2del} with $1\gg\epsilon>0$ close
to the singular limit. We use DDEBiftool \cite{ETR02} to construct bifurcation diagrams
and show numerically
that there are periodic solutions of \eq{eps2del} for $0<\epsilon\ll1$ which are perturbations of the unimodal
and type I and II bimodal solutions that we constructed in Section~\ref{secsingsols}. Moreover, we find fold bifurcations
close to the values predicted by our singular solutions. We also investigate multimodal solutions,
which are more complex than the singular solutions that we constructed
algebraically. The existence of these seems to be generic on the unstable legs of the branches between folds.

In Section~\ref{sec:cusp} we investigate the first two
codimension-two cusp-like bifurcations identified in Theorem~\ref{thmIIbif}. For $0<\epsilon\ll1$ we find
cusp bifurcations very close to the values predicted by the singular theory. We also show that these cusp bifurcations are
one mechanism by which stable bimodal periodic solutions may arise, and identify differences between the first and
second cusp bifurcation.

In Section~\ref{sec:othersols}, guided by our results from Section~\ref{secbifs} we investigate other
periodic solutions of \eq{eps2del} for $0<\epsilon\ll1$. For $A<3$ when folds do not occur,
we find an unbounded leg of stable type II bimodal solutions, and also period-doubling bifurcations, leading to stable period-doubled orbits. We also show an example of multimodal solutions with fold bifurcations which are associated with transitions between such solutions. We also consider the
alignment of the fold bifurcations on different solution branches and explain this using our results from
Section~\ref{secbifs}.
We finish in Section~\ref{sec:conc} with brief conclusions.

\section{Singular Solutions}
\label{secsingsols}

Before constructing singular solutions for \eq{eps2del},
as an illustrative example we consider the singular solutions of the one delay DDE \eq{eps1del}
which we write as
\be \label{eps1dels2}
\epsilon\udot(t)=-u(t)-K u(\alpha(t,u(t))), \quad \alpha(t,u(t))=t-a_1-cu(t).
\ee
We will
construct periodic singular solutions following Definition~\ref{defsingsol}
for \eq{eps1dels2} when $K>1$ (required for instability of the trivial solution),
with the profile below. Here, and throughout we use $\N_0$ to denote the natural numbers
including zero.

\begin{definition}[Sawtooth Profile] \label{defsawtooth}
For any $n\in \N_0$ and period $T>0$ the \emph{sawtooth profile} is
an admissible periodic singular solution profile on $I=\R$ defined by
\begin{gather} \label{idelg1}
\left.
\begin{aligned}
t(\mu)&=(\mu-i)T \\
u(\mu)&=\frac{-a_1+(n+\mu-2i)T}{c}
\end{aligned}\quad\right\}\quad
\mu\in[2i,2i+1], \\
\label{idelg2}
\left.
\begin{aligned}
t(\mu)&=(i+1)T \\
u(\mu)&=\frac{-a_1+(n+1-(\mu-2i-1))T}{c}
\end{aligned}\quad\right\}\quad
\mu\in((2i+1),(2i+2)),
\end{gather}
for each $i\in\Z$.
\end{definition}

Fig.~\ref{figadmissprof} shows a part of this profile when $a_1=c=1$.
Notice that $I^*$ is the union of the intervals $(2i,2i+1)$ and on each such interval
$u$ increases from $(-a_1+nT)/c$ to $(-a_1+(n+1)T)/c$ while $t$ increases by $T$.
$I^-$ is the union of the intervals $[2i+1,2i+2]$ and on each such interval
$u$ decreases from $(-a_1+(n+1)T)/c$ to $(-a_1+nT)/c$ while $t$ is fixed.
Mallet-Paret and Nussbaum have considered this $\Gamma$ (but not our parametrisation of it)
extensively, and named it the ``sawtooth profile''
for the shape of $\Gamma$ in Fig.~\ref{figadmissprof}(ii) \cite{JMPRNIII,JMPRN11}.

The motivation for Definition~\ref{defsawtooth} comes from
numerical simulations, where we observe when $\udot(t)$ is finite that $\alpha(t,u(t))$ is (almost) constant.
The sawtooth profile can then be constructed for \eq{eps1del} by assuming that
$\alpha(t,u(t))$ is constant with $\alpha(t,u(t))\in t(I^-)$
when $\udot(t)$ is finite (that is $u(t)\in \Gamma^*$).
If the phase of the periodic solution is chosen so that $t(I^-)=\{jT:j\in\Z\}$
then for $t\in(0,T)$
we have $-nT=\alpha(t,u(t))=t-a_1-cu(t)$ for some $n\in\N_0$. Rearranging this
leads to the formula for $u$ in \eq{idelg1} with $i=0$.

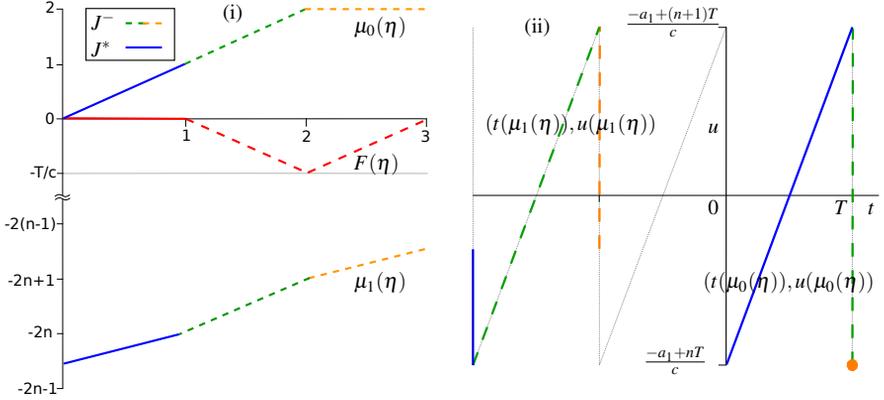
\begin{figure}[t!]
\begin{center}
\scalebox{0.4}{\includegraphics{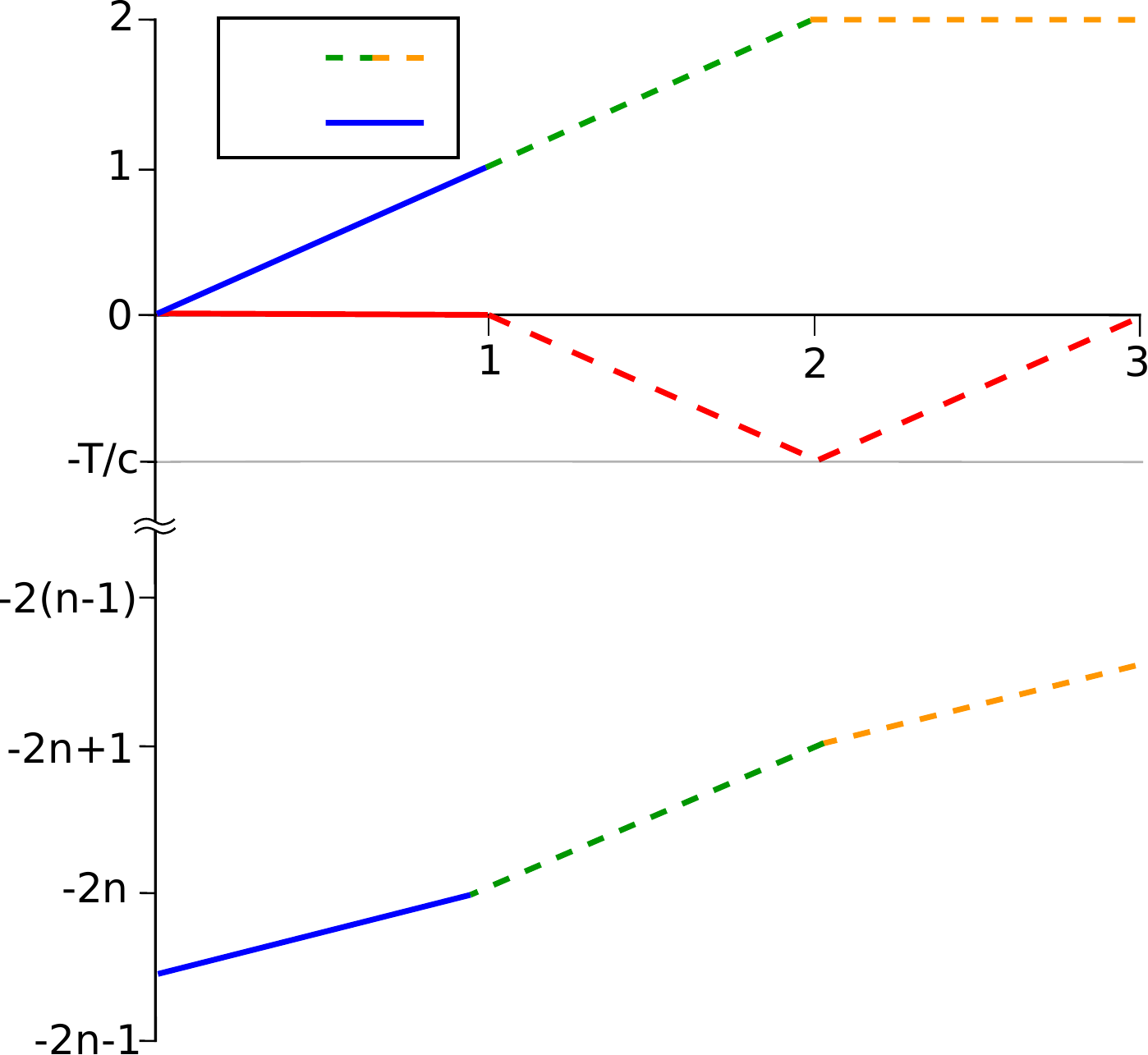}}\hspace{0.5cm}
\scalebox{0.4}{\raisebox{0.75cm}{\includegraphics{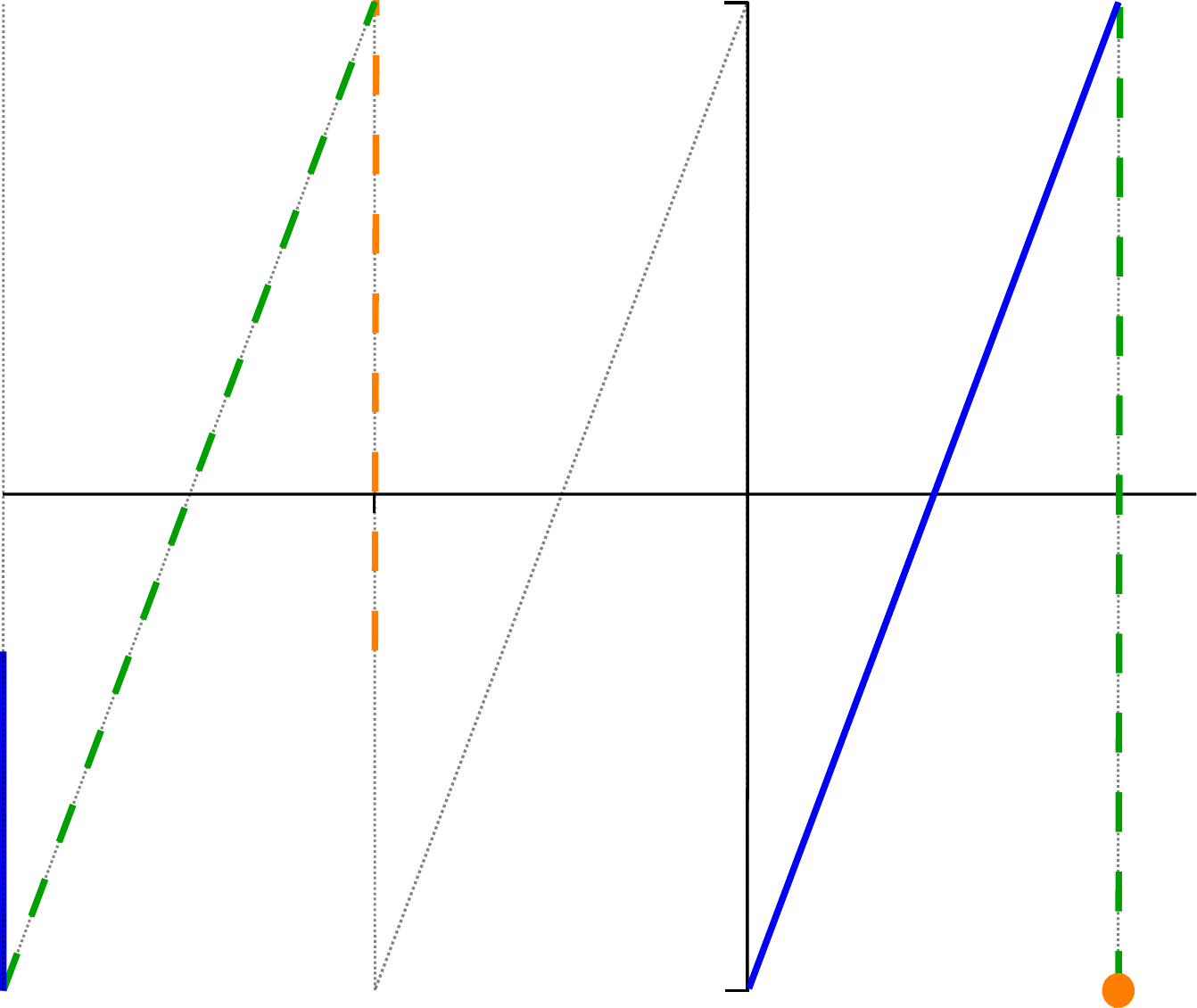}}}
\put(-250,142){(i)}
\put(-135,137){(ii)}
\put(-300,138){$J^-$}
\put(-300,128){$J^*$}
\put(-200,40){$\mu_1(\eta)$}
\put(-200,137){$\mu_0(\eta)$}
\put(-200,85){$F(\eta)$}
\put(-90,10){$\frac{-a_1+nT}{c}$}
\put(-65,68){$0$}
\put(-65,99){$u$}
\put(-17,68){$T$}
\put(-4,68){$t$}
\put(-67,40){$(t(\mu_0(\eta)),u(\mu_0(\eta))$}
\put(-150,100){$(t(\mu_1(\eta)),u(\mu_1(\eta))$}
\put(-97,138){$\frac{-a_1+(n+1)T}{c}$}
\end{center}
\caption{(i) $\mu_0(\eta)$, $\mu_1(\eta)$ and $F(\eta)$ for $\eta\in[0,3]$ for the singular solution of
\eq{eps1del}
defined by \eq{idelg1}--\eq{mu1del}. (ii) The corresponding periodic singular solution
$(t(\mu_0(\eta)),u(\mu_0(\eta)))$ and delayed solution
$(t(\mu_1(\eta)),u(\mu_1(\eta)))$
for $\eta\in[0,3]$.}
\label{fig1delpar}
\end{figure}

Each different $n$ will define a different singular solution,
with delay $t-\alpha(t,u(t))=a_1+cu(t)\in[nT,(n+1)T]$.
Here we will construct
singular solutions of \eq{eps1dels2} for all $n\in \N_0$ with period $T$ given by
\be \label{T1del}
T=\frac{a_1(1+K)}{1+n(1+K)}.
\ee
Later, we will construct periodic singular solutions of the two
delay equation \eq{eps2del} using the same sawtooth admissible solution profile.
To define a singular solution for \eq{eps1del} with this profile,
for $j\in\Z$ let
\be \label{mu1del}
\begin{tabular}{c|c|c}
$\mu_0(\eta)=$ & $\mu_1(\eta)=$  & $\eta\in$ \\ \hline
$2j+(\eta-3j)$ & \!$2(j-n)-1+[(\eta-3j)+(K-1)]/K$\! & $[3j,3j+1]$ \\
\!\!$2j+1+(\eta-3j-1)$\! & $2(j-n)+(\eta-3j-1)$ & \!$[3j+1,3j+2]$\!\! \\
$2j+2$ & \!$2(j-n)+1+(\eta-3j-2)(K-1)/K$\! &  \!$[3j+2,3j+3]$\!\!
\end{tabular}\\
\ee
Then $\mu_i(\eta)$ is continuous on the real line. It is a simple
but tedious algebraic exercise to check that \eq{singdel} holds for
all $\eta\in\R$. Notice in particular that for $\eta\in[3j,3j+1]$
we have $\mu_1(\eta)\in\bigl(2(j-n)-1,2(j-n)\bigr)$ provided $K>1$,
in which case
$$t(\mu_1(\eta))=(j-n)T=t(\mu_0(\eta))-a_1-cu(\mu_0(\eta))=\alpha\bigl(t(\mu_0(\eta)),u(\mu_0(\eta))\bigr),$$
as required to satisfy \eq{singdel}.
Before checking the conditions on $F(\eta)$, notice that $\mu_0(\eta)\in I^*$ for
$\eta\in(3j,3j+1)$, $\mu_0(\eta)\in int(I^-)$ for $\eta\in(3j+1,3j+2)$
and $\mu_0(\eta)\in \partial I^*=\partial I^-$ for $\eta\in[3j+2,3j+3]$ for each $j\in\Z$.
Hence $J^*$ is the union of the intervals $[3j,3j+1]$, while
$J^-$ is composed of intervals $(3j+1,3j+3)$. For $\eta\in J^*$ we have $F(\eta)=0$
provided \eq{T1del} holds (which is how $T$ was actually determined).
For $\eta\in(3j+1,3j+2]$ we have
\begin{align*}
F(\eta)&=-u(\mu_0(\eta))-Ku(\mu_1(\eta))\\
&=-\left[\frac{-a_1+(n+3j+2)T-\eta T}{c}\right]-K\left[\frac{-a_1+(n-3j+1)T+\eta T}{c}\right]\\
&=\frac{1}{c}(a_1-nT)(1+K)-\frac{T}{c}+\frac{1}{c}\bigl(\eta-(3j+1)\bigr)T(1-K)
=\frac{1}{c}\bigl(\eta-(3j+1)\bigr)T(1-K),
\end{align*}
and hence $F(\eta)<0$ for all $\eta\in(3j+1,3j+2]$, since $K>1$. Finally on the interval $[3j+2,3j+3]$,
we have $u(\mu_1(\eta))$ is a linear function of $\eta$, while $u(\mu_0(\eta))$ is constant, and
hence $F(\eta)$ is a linear function of $\eta$. By continuity and the previous
calculations $F(3j+2)=T(1-K)/c<0$ and $F(3j+3)=0$ hence $F(\eta)<0$ for
all $\eta\in J^-$ as required. Thus for each $n\in\N_0$ and each $K>1$ we have constructed a periodic singular
solution of \eq{eps1del} defined by
\eq{T1del}-\eq{mu1del}.
The parametrisation
leading to one of these solutions and the corresponding periodic singular solution is illustrated in
Fig.~\ref{fig1delpar}.

Using max-plus equations, in \cite{JMPRNIII} this $\Gamma$ is proved to
be the limiting profile of the slowly oscillating periodic solutions (corresponding to $n=0$)
of \eq{eps1del} as $\epsilon\to0$. In \cite{JMPRN11} higher order asymptotics reveal the shape of the
periodic solution for $0<\epsilon\ll1$. It is noted that that the asymptotic
forms of the periodic solution are very different near to the local maximum and minimum of the solution, with the
maximum corresponding to a regular point of the dynamics scaled by $\epsilon$, while the minimum can be interpreted
in the spirit of Fenichel as a turning point near a normally hyperbolic invariant manifold for an ordinary differential
equation with a time scaling of $\epsilon^2$ \cite{JMPRN11}. The singular solution
\eq{idelg1}--\eq{mu1del} also reveals a difference between the dynamics near to the maximum and minimum of the
periodic solution. The solution $u(t(\mu_0(\eta)))$ has its maximum when $\eta=3j+1$ (for any integer $j$),
which is at the boundary between two of the linear segments in the solution parametrisation \eq{mu1del},
corresponding to the boundary between $J^*$ and $J^-$. In contrast
$u(t(\mu_0(\eta)))$ takes its minimum value on the entire interval $\eta\in[3j+2,3j+3]$ (but $u(t(\mu_1(\eta)))$ is
not constant on this interval). Note that while at first sight it may have appeared more natural
in Definition~\ref{defsingsol} to define $J^*$
to be the set of $\eta\in J$ such that $\mu_0(\eta)\in I^*$ (or equivalently such that $u(\mu_0(\eta))\in\Gamma^*$),
such a definition would be problematical in the example above because $\mu_0(\eta)$ is constant on $\partial I^*$
on the interval $\eta\in[3j+2,3j+3]$. We will also find nontrivial intervals on which $\mu_0(\eta)$ is constant on $\partial I^*$
for singular solutions of \eq{eps2del}.

Now consider the two delay DDE \eq{eps2del}. We assume several conditions on the positive parameters. Without
loss of generality we assume that $a_2>a_1$ (if not we can either swap the
order of the terms, or reduce to an equation with one delay). Then letting
$\alpha_i(t,u(t))=t-a_i-cu(t)$ we see that $\alpha_2(t,u(t))<\alpha_1(t,u(t))$
with $\alpha_1(t,u(t))-\alpha_2(t,u(t))=a_2-a_1>0$, constant. So although
the arguments $\alpha_i(t,u(t))$ are both linearly state-dependent, the difference
between them is constant. The more general case where $\alpha_i(t,u(t))=t-a_i-c_iu(t)$
with $c_1\ne c_2$ so that the difference between the delays is nonconstant
would also be interesting, but in the current work we concentrate on understanding
the simpler case, which already leads to very complicated dynamics.

It is useful to
define the ratio $A=a_2/a_1>1$ which will play an important role later.
If $K_1+K_2<1$ the trivial solution
is asymptotically stable and there
are no stable periodic solutions, so we assume that $K_1+K_2>1$. Finally
we assume that
\be \label{K2<1}
K_2<1.
\ee
It is shown in \cite{DCDSA11} that \eq{K2<1} along with $A>1$
ensures that the DDE initial value problem is well-posed
for \eq{eps2del}, and in particular that the delay $\alpha_1(t,u(t)))<t$ and so
does not become advanced. It is also shown in \cite{DCDSA11} that when $\epsilon>0$ the function
$\alpha_i(t,u(t)))$ is a strictly monotonic increasing function of $t$ for
$t\geq a_2+a_1(K_1+K_2)$. Hence
$\alpha_i(t,u(t)))$ must be a strictly monotonic increasing function of $t$
on any periodic solution.
Thus we will construct singular periodic solutions for which all the $\mu_i(\eta)$ are monotonic
increasing functions of $\eta$ for all $i$, although Definition~\ref{defsingsol}
only requires that $\mu_0(\eta)$ be monotonic in general.


We first construct singular periodic solutions for \eq{eps2del} which
have the same sawtooth profile \eq{idelg1},\eq{idelg2} as the sawtooth solutions
of the one delay DDE \eq{eps1del}. Since these solutions have one local maxima per period
we refer to them as unimodal. We will then construct two types singular periodic solution
with two local maxima per step; type I and type II bimodal solutions. Each of the
solutions that we construct of each type will be characterised by a pair $(n,m)$ of
non-negative integers which will have the same meaning in each case. The first number $n$
is the integer number of periods in the past that the first delay falls, and the second
number $m$ is the integer number of periods between the two delay times $\alpha_1(t,u(t)))$ and $\alpha_2(t,u(t)))$.
So for a singular solution of period $T$ we always have
\be \label{nm}
t-\alpha_1\!(t,u(t))\in[nT\!,(n+\!1)T], \quad a_2-a_1=\alpha_1(t,u(t))-\alpha_2(t,u(t))\in(mT\!,(m+\!1)T).
\ee
Or using the parametrisation
\be \label{n}
t(\mu_0(\eta))-t(\mu_1\!(\eta))=t(\mu_0(\eta))-\alpha_1\!(t(\mu_0(\eta)),u(\mu_0(\eta)))
\in[nT\!,(n+1)T],\quad\forall\eta\in\R,\; n\in\N_0.
\ee
and
\begin{align} \notag
t(\mu_1(\eta))-t(\mu_2(\eta))
&=\alpha_1(t(\mu_0(\eta)),u(\mu_0(\eta)))-\alpha_2(t(\mu_0(\eta)),u(\mu_0(\eta)))\\ \notag
&=\bigl[t(\mu_0(\eta))-a_1-cu(\mu_0(\eta))\bigr]-\bigl[t(\mu_0(\eta))-a_2-cu(\mu_0(\eta))\bigr]\\
&=a_2-a_1 \in(mT,(m+1)T),\quad\forall\eta\in\R,\; m\in\N_0. \label{m}
\end{align}
With $n$ and $m$ defined by \eq{n} and \eq{m} to construct unimodal singular solutions
of \eq{eps2del} it is useful to define $\theta\in(0,1)$
by
$$t(\mu_1(\eta))-t(\mu_2(\eta))=a_2-a_1=(m+\theta)T, \quad \theta\in(0,1),$$
so $\theta$ is the fractional part of a period between the two delays, which is assumed to be non-zero.
(Although $n$ and $m$ will always have the same meaning, $\theta$ will be defined slightly
differently for each type of bimodal solution.)
As in the one delay case we will construct a solution
with
$t(\mu_1(\eta))=-nT$ while $t(\mu_0(\eta))\in(0,T)$. The following theorem establishes conditions for such a solution to exist.

\begin{theorem} \label{thm2delunimod}
Let $K_1>1>K_2>0$, $a_2>a_1>0$, $m,n\in\N_0$,
\be \label{T2delunimod}
T=\frac{a_1(1+K_1+K_2)+(a_2-a_1)K_2}{1+(m+1)K_2+n(1+K_1+K_2)}
\ee
and
\be \label{theta2delunimod}
\theta=\frac{a_2-a_1}{T}-m.
\ee
The DDE \eq{eps2del} has a periodic singular solution with profile \eq{idelg1},\eq{idelg2}
and period $T>0$ given by \eq{T2delunimod}
when the parameters are chosen so that
\be \label{thetabds}
\theta\in\left(\frac{K_2}{K_1+K_2-1},1\right).
\ee
\end{theorem}

\bpf
For $j\in\Z$ let $\mu_i(\eta)$ be defined by Table~\ref{tab:etaunimodal}.
%
By the conditions of the theorem, $\theta\in(0,1)$ and $K_1>1$ . From this it follows that each $\mu_i(\eta)$
is continuous and monotonically increasing.
For $\eta\in[5j+k,5j+k+1]$ for $k=0,1,2,3,4$, notice that each function $\mu_i(\eta)$ is linear in $\eta$,
and falls into a single subinterval of the sawtooth profile
defined by \eq{idelg1},\eq{idelg2}, and so $u(\mu_i(\eta))$ and $t(\mu_i(\eta))$
are linear functions for $\eta\in[5j+k,5j+k+1]$.
It follows that $F(\eta)$ is also linear in $\eta$ for $\eta\in[5j+k,5j+k+1]$ for each integer $k$.
It is straightforward to confirm that
\eq{singdel} holds, that is $t(\mu_i(\eta))=t(\mu_0(\eta))-a_i-cu(\mu_0(\eta))$ for $i=1,2$.

\afterpage{\begin{landscape}
\begin{table}
\centering    
\begin{tabular}{c|c|c|c}
$\mu_0(\eta)=$ & $\mu_1(\eta)=$ & $\mu_2(\eta)=$ & $\eta\in$ \\ \hline
$2j\!+\!\eta\!-\!5j$ & $2j\!-\!2n\!+\!(\!-\!1\!+\!\eta\!-\!5j)/K_1$
& $2(j\!-\!n\!-\!m)\!-\!1\!-\!\theta$ & $[5j,5j\!+\!1]$ \\
$2j\!+\!1\!+\!(\eta\!-\!5j\!-\!1)\theta$ & $2(j\!-\!n)\!+\!(\eta\!-\!5j\!-\!1)\theta$ &
$2(j\!-\!n\!-\!m)\!-\!1\!-\!\theta\!+\!(\eta\!-\!5j\!-\!1)\theta$ & $[5j\!+\!1,5j\!+\!2]$ \\
$2j\!+\!1\!+\!\theta$ & $2(j\!-\!n)\!+\!\theta$ & $2(j\!-\!n\!-\!m)\!-\!1\!+\!(\eta\!-\!5j\!-\!2)$ & $[5j\!+\!2,5j\!+\!3]$ \\
$2j\!+\!1\!+\!\theta\!+\!(\eta\!-\!5j\!-\!3)(1\!-\!\theta)$ & $2(j\!-\!n)\!+\!\theta\!+\!(\eta\!-\!5j\!-\!3)(1\!-\!\theta)$ &
$2(j\!-\!n\!-\!m)\!+\!(\eta\!-\!5j\!-\!3)(1\!-\!\theta)$ & $[5j\!+\!3,5j\!+\!4]$ \\
$2j\!+\!2$ & $2(j\!-\!n)\!+\!1\!+\!(1\!-\!1/K_1)(\eta\!-\!5j\!-\!4)$
& $2(j\!+\!1\!-\!n\!-\!m)\!-\!1\!-\!\theta$ & $[5j\!+\!4,5j\!+\!5]$ \\
\end{tabular}
\caption{Parameterization of the unimodal solution defined in Theorem~\ref{thm2delunimod} for the
sawtooth profile given in Definition~\ref{defsawtooth}.} \label{tab:etaunimodal}
\end{table}
\begin{table}
\vspace{-2ex}
\begin{center}
{\small\begin{tabular}{c|c|c|c}
$\mu_0(\eta)=$ & $\mu_1(\eta)=$ & $\mu_2(\eta)=$ & $\eta\in$ \\ \hline
$4j\!+\!(\eta\!-\!10j)$ & $4(j\!-\!n)\!-\!1\!+\!s_{13}\!+\!(1\!-\!s_{13})(\eta\!-\!10j)$
 & $4(j\!-\!n\!-\!m)\!-\!3\!-\!\theta$ & $[10j,10j\!+\!1]$ \\
$4j\!+\!1\!+\!(\eta\!-\!10j\!-\!1)s_{11}$ & $4(j\!-\!n)\!+\!(\eta\!-\!10j\!-\!1)\theta$
& $4(j\!-\!n\!-\!m)\!-\!3\!-\!\theta\!+\!(\eta\!-\!10j\!-\!1)\theta$ & $[10j\!+\!1,10j\!+\!2]$ \\
$4j\!+\!1\!+\!s_{11}$ & $4(j\!-\!n)\!+\!\theta$ & $4(j\!-\!n\!-\!m)\!-\!3\!+\!(\eta\!-\!10j\!-\!2)$ & $[10j\!+\!2,10j\!+\!3]$ \\
$4j\!+\!2\!+\!(\eta\!-\!10j\!-\!4)(1\!-\!s_{11})$ & $4(j\!-\!n)\!+\!\theta\!+\!(\eta\!-\!10j\!-\!3)T_2/T_1$
& $4(j\!-\!n\!-\!m)\!-\!2\!+\!(\eta\!-\!10j\!-\!3)$ & $[10j\!+\!3,10j\!+\!4]$ \\
$4j\!+\!2$ & $4(j\!-\!n)\!+\!\theta\!+\!T_2/T_1$ & $4(j\!-\!n\!-\!m)\!-\!1\!+\!(\eta\!-\!10j\!-\!4)s_{14}$ & $[10j\!+\!4,10j\!+\!5]$ \\
$4j\!+\!2\!+\!(\eta\!-\!10j\!-\!5)$ & $4(j\!-\!n)\!+\!\theta\!+\!T_2/T_1$ & $4(j\!-\!n\!-\!m)\!+\!(\!-\!1\!+\!\eta\!-\!10j\!-\!5)(1\!-\!s_{14})$
& $[10j\!+\!5,10j\!+\!6]$ \\
$4j\!+\!3\!+\!(\eta\!-\!10j\!-\!6)s_{12}$ & $4(j\!-\!n)\!+\!1\!+\!(\eta\!-\!10j\!-\!7)(1\!-\!\theta\!-\!T_2/T_1)$ &
$4(j\!-\!n\!-\!m)\!+\!(1\!-\!\theta\!-\!T_2/T_1)(\eta\!-\!10j\!-\!6)$ & $[10j\!+\!6,10j\!+\!7]$ \\
$4j\!+\!3\!+\!s_{12}$ & $4(j\!-\!n)\!+\!1\!+\!(\eta\!-\!10j\!-\!7)$ & $4(j\!-\!n\!-\!m)\!+\!1\!-\!\theta\!-\!T_2/T_1$ & $[10j\!+\!7,10j\!+\!8]$ \\
$4(j\!+\!1)\!+\!(\eta\!-\!10j\!-\!9)(1\!-\!s_{12})$ & $4(j\!-\!n)\!+\!2\!+\!(\eta\!-\!10j\!-\!8)$
& $4(j\!-\!n\!-\!m)\!+\!1\!-\!\theta\!+\!T_2/T_1(\eta\!-\!10j\!-\!9)$ & $[10j\!+\!8,10j\!+\!9]$ \\
$4(j\!+\!1)$ & $4(j\!-\!n)\!+\!3\!+\!(\eta\!-\!10j\!-\!9)s_{13}$ & $4(j\!+\!1\!-\!n\!-\!m)\!-\!3\!-\!\theta$
& \!$[10j\!+\!9,10(j\!+\!1)]$\!\!
\end{tabular}}
\caption{Parameterization of the type I bimodal solution defined in Theorem~\ref{thm2delI} for the
admissible profile defined by \eq{eqtypeI}, illustrated in Fig.~\ref{figtypeI}.} \label{tab:etatypeI}
\end{center}
\end{table} 
\nopagebreak
\begin{table}
\vspace{-2ex}
\begin{center}
{\small
\begin{tabular}{c|c|c|c}
$\mu_0(\eta)=$ & $\mu_1(\eta)=$ & $\mu_2(\eta)=$ & $\eta\in$ \\ \hline
$4j\!+\!(\eta\!-\!10j)$ & $4(j\!-\!n)\!+\!(\!-\!1\!+\!\eta\!-\!10j)(1\!-\!s_{23})$ & $4(j\!-\!n\!-\!m)\!-\!1\!-\!\theta$ & $[10j,10j\!+\!1]$ \\
$4j\!+\!1\!+\!(\eta\!-\!10j\!-\!1)$ & $4(j\!-\!n)\!+\!(\eta\!-\!10j\!-\!1)\theta T_2/T_1$ & $4(j\!-\!n\!-\!m)\!-\!1\!-\!\theta\!+\!\theta(\eta\!-\!10j\!-\!1)$
& $[10j\!+\!1,10j\!+\!2]$ \\
$4j\!+\!2$ & $4(j\!-\!n)\!+\!\theta T_2/T_1$ & $4(j\!-\!n\!-\!m)\!-\!1\!+\!(\eta\!-\!10j\!-\!2)s_{24}$ & $[10j\!+\!2,10j\!+\!3]$ \\
$4j\!+\!2\!+\!(\eta\!-\!10j\!-\!3)$ & $4(j\!-\!n)\!+\!\theta T_2/T_1$ &
$4(j\!-\!n\!-\!m)\!+\!(\eta\!-\!10j\!-\!4)(1\!-\!s_{24})$ & $[10j\!+\!3,10j\!+\!4]$ \\
$4j\!+\!3\!+\!(\eta\!-\!10j\!-\!4)s_{21}$ & $4(j\!-\!n)\!+\!1+\!(1\!-\!\theta T_2/T_1)(\eta\!-\!10j\!-\!5)$
& $4(j\!-\!n\!-\!m)+\!(1\!-\!\theta T_2/T_1)(\eta\!-\!10j\!-\!4)$ & $[10j\!+\!4,10j\!+\!5]$ \\
$4j\!+\!3\!+\!s_{21}$ & $4(j\!-\!n)\!+\!1\!+\!(\eta\!-\!10j\!-\!5)$ & $4(j\!-\!n\!-\!m)\!+\!1\!-\!\theta T_2/T_1$ & $[10j\!+\!5,10j\!+\!6]$ \\
$\!\!4j\!+\!3\!+\!s_{21}+\!(\eta\!-\!10j\!-\!6)(s_{22}\!-\!s_{21})$ & $4(j\!-\!n)\!+\!2\!+\!(\eta\!-\!10j\!-\!6)\theta$
& $4(j\!-\!n\!-\!m)\!+\!1+\!\theta T_2/T_1(\eta\!-\!10j\!-\!7)$ & $[10j\!+\!6,10j\!+\!7]$ \\
$4j\!+\!3\!+\!s_{22}$ & $4(j\!-\!n)\!+\!2\!+\!\theta$ & $4(j\!-\!n\!-\!m)\!+\!1\!+\!(\eta\!-\!10j\!-\!7)$ & $[10j\!+\!7,10j\!+\!8]$ \\
$4(j\!+\!1)\!+\!(\eta\!-\!10j\!-\!9)(1\!-\!s_{22})$ & $4(j\!-\!n)\!+\!2\!+\!\theta\!+\!(\eta\!-\!10j\!-\!8)(1\!-\!\theta)$ &
$4(j\!-\!n\!-\!m)\!+\!2\!+\!(1\!-\!\theta)(\eta\!-\!10j\!-\!8)$ & $[10j\!+\!8,10j\!+\!9]$ \\
$4(j\!+\!1)$ & $4(j\!-\!n)\!+\!3\!+\!(\eta\!-\!10j\!-\!9)s_{23}$  & $4(j\!+\!1\!-\!n\!-\!m)\!-\!1\!-\!\theta$ & $\![10j\!+\!9,10(j\!+\!1)]\!$
\end{tabular}}
\caption{Parameterization of the type II bimodal solution defined in Theorem~\ref{thm2delII} for the
admissible profile defined by \eq{eqtypeII}, illustrated in Fig.~\ref{figtypeII}.} \label{tab:etatypeII}
\end{center}\end{table}
\end{landscape}}

It remains to establish the conditions on $F$. First note that $J^*=\bigcup_{j\in\N}[5j,5j+1]$.
Now
\begin{align*}
F(5j)&=-u(\mu_0(5j))-K_1u(\mu_1(5j))-K_2u(\mu_2(5j))\\
&=-u(2j)-K_1u(2j-2n-1/K_1)-K_2u(2(j-n-m)-1-\theta)\\
&=-\Bigl(\frac{-a_1+nT}{c}\Bigr)-K_1\Bigl(\frac{-a_1+(n+1/K_1)T}{c}\Bigr)-K_2\Bigl(\frac{-a_1+(n+1-\theta)T}{c}\Bigr)
\end{align*}
hence
$$cF(5j)=(a_1-nT)(1+K_1+K_2)-T-(1-\theta)K_2T.$$
But multiplying \eq{T2delunimod} by its denominator, and noting that from \eq{theta2delunimod}
we have $a_2-a_1=(m+\theta)T$, we see that
$$(a_1-nT)(1+K_1+K_2)=-(a_2-a_1)K_2+T+(m+1)K_2T=T+(1-\theta)K_2T,$$
and hence $F(5j)=0$. It follows similarly that $F(5j+1)=0$, and
hence by linearity, $F(\eta)=0$ for all $\eta\in[5j,5j+1]$ and hence for all $\eta\in J^*$.

It remains to show that $F(\eta)<0$ for $\eta\in J^-=\bigcup_{j\in\N}(5j+1,5j+5)$.
Since $F(5j)=F(5j+1)=F(5j+5)=0$, by the linearity of $F(\eta)$ on each subinterval, it
is sufficient to show that $F(5j+2)<0$, $F(5j+3)<0$ and $F(5j+4)<0$. But similarly
to above we derive
$$cF(5j+2)=(1-K_1-K_2)\theta T, \quad cF(5j+4)=(1-K_1)T,$$
which are both negative since $K_1>1$, while
\be \label{cF5j3}
cF(5j+3)=cF(5j+2)+K_2T=\bigl[K_2-(K_1+K_2-1)\theta\bigr]T,
\ee
and $F(5j+3)<0$ provided $\theta>K_2/(K_1+K_2-1)$. Hence $F(\eta)<0$ for all $\eta\in J^-$,
which completes the proof.
\epf

Theorem~\ref{thm2delunimod} shows immediately that $\theta$ is bounded away from zero.
We will see in Section~\ref{secbifs} that only certain pairs of values of $m,n\in\N_0$
satisfy the bounds \eq{thetabds} in Theorem~\ref{thm2delunimod}.
In Theorem~\ref{thmunilegs} we will determine which pairs $(n,m)$ are possible and for
which parameter ranges the conditions of
Theorem~\ref{thm2delunimod} are satisfied to begin to construct a bifurcation
diagram of solution branches.
For now, we note that using \eq{T2delunimod} and \eq{theta2delunimod}
we can write
$$m+\theta=\frac{(A-1)\bigl(1+(m+1)K_2+n(1+K_1+K_2)\bigr)}{1+K_1+K_2+(A-1)K_2},$$
where $A=a_2/a_1$. Using this, the condition $\theta>K_2/(K_1+K_2-1)$ can be rewritten as
\be \label{Gnm<0}
G_{nm}(K_1)<0,
\ee
where
\be
G_{nm}(K_1)=\bigl[m-n(A-1)\bigr]\bigl((K_1+K_2)^2-1\bigr)
-K_1\bigl[(A-1)(1+K_2)-K_2\bigr]+K_2(1+K_2)+(A-1). \label{Gnm}
\ee

When the parameters are such that the bounds on $\theta$ in \eq{thetabds} are
violated other types of singular solution arise. We will construct
two such classes of solutions which we refer to as type I and type II bimodal solutions,
since each has two local maxima per period.

\begin{figure}[t!]
\vspace{-4.5ex}
\mbox{}\hspace{-0.5cm}\scalebox{0.48}{\includegraphics{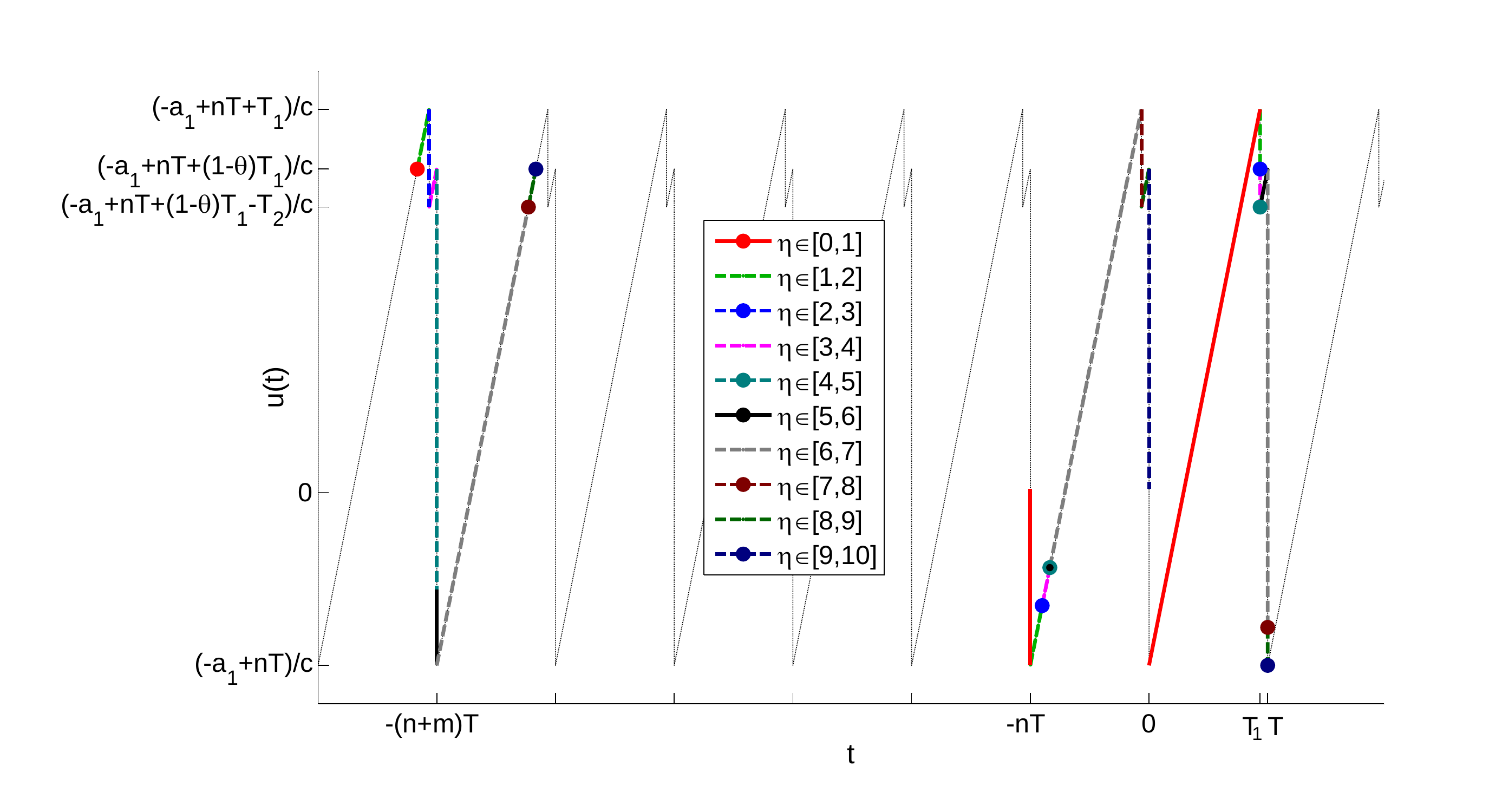}}
\put(-268,100){\rotatebox{79}{\footnotesize $(t(\mu_2(\eta),u(\mu_2(\eta))$}}
\put(-115,100){\rotatebox{79}{\footnotesize $(t(\mu_1(\eta),u(\mu_1(\eta))$}}
\put(-80,65){\rotatebox{79}{\footnotesize $(t(\mu_0(\eta),u(\mu_0(\eta))$}}

\vspace{-3ex}
\caption{A generic Type I admissible bimodal periodic solution profile as defined in \eq{eqtypeI}. Also
shown in colour are the ten stages of the parametrisation of $(t(\mu_i(\eta),u(\mu_i(\eta))$ from the
proof of Theorem~\ref{thm2delI} for $i=0,1,2$ and $j=0$. Where plotted the dots indicate that $(t(\mu_i(\eta),u(\mu_i(\eta))$
is constant for that stage of the parametrisation, with the multicoloured dot showing that $(t(\mu_1(\eta),u(\mu_1(\eta))$ is constant
for two successive stages for $\eta\in[4,6]$.}
\label{figtypeI}
\end{figure}

\begin{figure}
\vspace{-4.5ex}
\mbox{}\hspace{-0.5cm}\scalebox{0.48}{\includegraphics{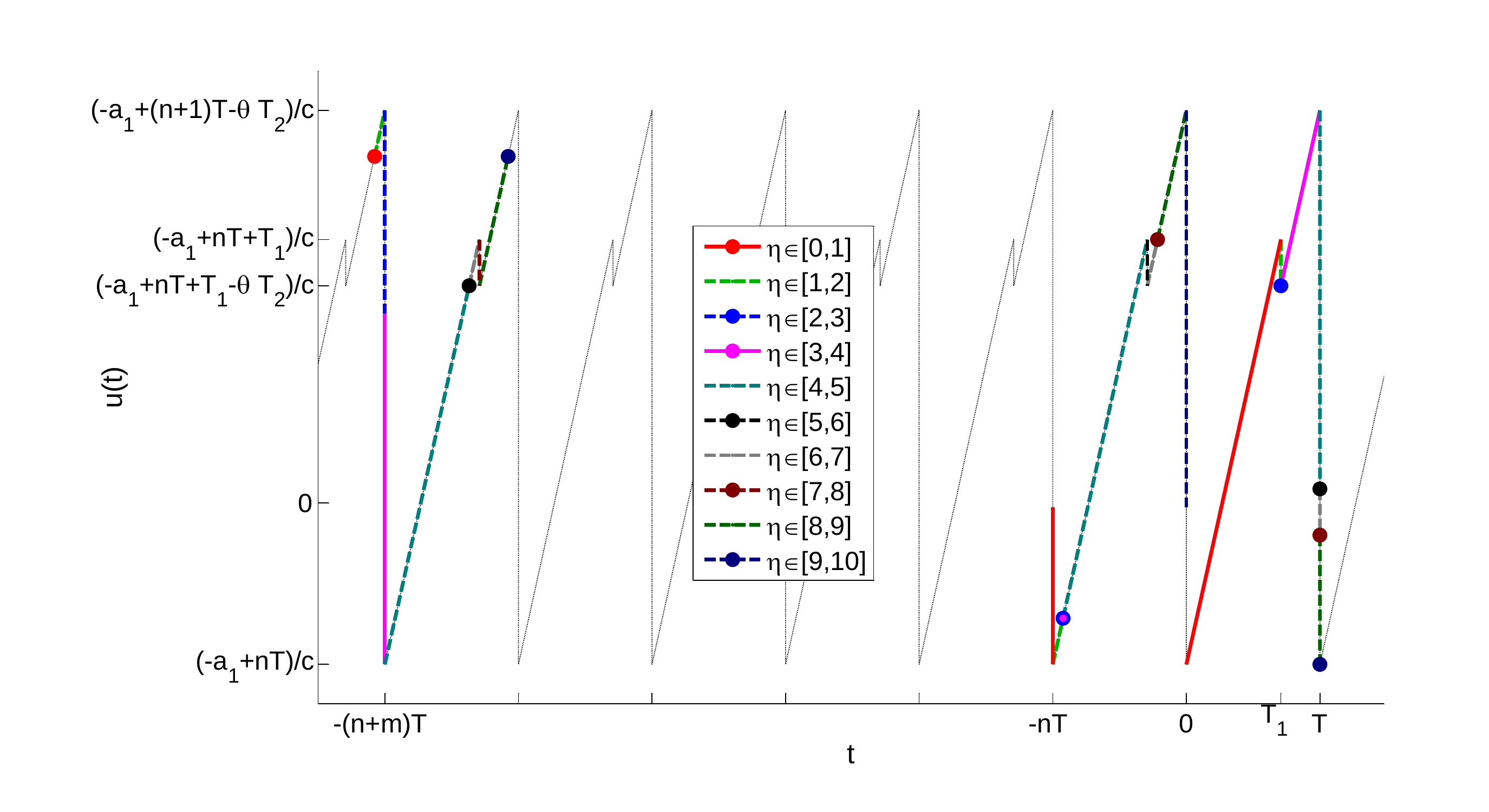}}
\put(-281,85){\rotatebox{79}{\footnotesize $(t(\mu_2(\eta),u(\mu_2(\eta))$}}
\put(-103,60){\rotatebox{79}{\footnotesize $(t(\mu_1(\eta),u(\mu_1(\eta))$}}
\put(-68,65){\rotatebox{79}{\footnotesize $(t(\mu_0(\eta),u(\mu_0(\eta))$}}
\vspace{-3ex}
\caption{A generic Type II admissible bimodal periodic solution profile as defined in \eq{eqtypeII}.
Also
shown in colour are the ten stages of the parametrisation of $(t(\mu_i(\eta),u(\mu_i(\eta))$ from the
proof of Theorem~\ref{thm2delII} for $i=0,1,2$ and $j=0$.}
\label{figtypeII}
\end{figure}

Let $n\in\N_0$ and $m\in\N_0$ be related to the delays and period $T$ as explained in \eq{nm}-\eq{m}.
For $\theta\in(0,1)$, $T=T_1+T_2$ where $T_i>0$,
the Type I and Type II bimodal periodic admissible singular solution profiles
are defined by
\be \label{eqtypeI}
\begin{tabular}{c|c|c}
$t(\mu)=$ & $u(\mu)=$ & $\mu\in$ \\ \hline
\rule[-4pt]{0pt}{15pt}$(\mu\!-\!4i)T_1\!+\!iT$ & $\frac{1}{c}(-a_1\!+\!nT\!+\!(\mu\!-\!4i)T_1)$ & $[4i,4i\!+\!1]$ \\
\rule[-4pt]{0pt}{15pt}$T_1\!+\!iT$ & $\frac{1}{c}(-a_1\!+\!nT\!+\!T_1\!-\!(\mu\!-\!4i\!-\!1)(T_2\!+\!\theta T_1)$ & $[4i\!+\!1,4i\!+\!2]$ \\
\rule[-4pt]{0pt}{15pt}$(\mu\!-\!4i\!-\!2)T_2\!+\!T_1\!+\!iT$ & $\frac{1}{c}(-a_1\!+\!nT\!+\!(1\!-\!\theta)T_1\!+\!(\mu\!-\!4i\!-\!3)T_2)$
& $[4i\!+\!2,4i\!+\!3]$ \\
\rule[-4pt]{0pt}{15pt}$(i\!+\!1)T$ & $\frac{1}{c}(-a_1\!+\!nT\!+\!(4i\!+\!4\!-\!\mu)(1\!-\!\theta)T_1)$  & $[4i\!+\!3,4i\!+\!4]$,
\end{tabular}
\ee
and
\be \label{eqtypeII}
\begin{tabular}{c|c|c}
$t(\mu)=$ & $u(\mu)=$ & $\mu\in$ \\ \hline
\rule[-4pt]{0pt}{15pt}$(\mu\!-\!4i)T_1\!+\!iT$ & $\frac{1}{c}(-a_1\!+\!nT\!+\!(\mu\!-\!4i)T_1)$ & $[4i,4i\!+\!1]$ \\
\rule[-4pt]{0pt}{15pt}$T_1\!+\!iT$ & $\frac{1}{c}(-a_1\!+\!nT\!+\!T_1\!-\!(\mu\!-\!4i\!-\!1)\theta T_2)$ & $[4i\!+\!1,4i\!+\!2]$ \\
\rule[-4pt]{0pt}{15pt}$(\mu\!-\!4i\!-\!2)T_2\!+\!T_1\!+\!iT$ & $\frac{1}{c}(-a_1\!+\!nT\!+\!T_1\!-\!\theta T_2+(\mu\!-\!4i\!-\!2)T_2)$
& $[4i\!+\!2,4i\!+\!3]$ \\
\rule[-4pt]{0pt}{15pt}$(i\!+\!1)T$ & $\frac{1}{c}(-a_1\!+\!nT\!+\!(T\!-\!\theta T_2)(4i\!+\!4\!-\!\mu)$ & $[4i\!+\!3,4i\!+\!4]$,
\end{tabular}
\ee
respectively. These profiles are illustrated in Figs.~\ref{figtypeI} and~\ref{figtypeII}.

We see from the \eq{eqtypeI} and \eq{eqtypeII} that both solutions have global minima with $u=(-a_1+nT)/c$.
If the phase of the periodic solution is chosen so that these minima occur when $t=jT$,
for integer $j$,
then for type I bimodal solutions the first local maximum which occurs when
$t=jT+T_1$ is also the global maximum, while for type II bimodal solutions the
second local maximum on the period is equal to the global maximum.

The following theorem identifies all the conditions on the parameters for
a type I bimodal solution to exist. In Theorem~\ref{thmIbif} we find parameter
ranges for which all these conditions are satisfied. The integers $n$ and $m$ in Theorem~\ref{thm2delI}
have similar geometrical meanings as for the sawtooth
solution, so $n$ and $m$ again satisfy \eq{nm}-\eq{m}. For the type I bimodal solution it is convenient to
define $\theta\in(0,1)$ by
$a_2-a_1=mT+T_2+\theta T_1$ where $m\in\N_0$, $\theta\in(0,1)$ and $T=T_1+T_2$
so
\begin{align*}
t(\mu_1(\eta))-t(\mu_2(\eta))
=\alpha_1(t(\mu_0(\eta)),\; & u(\mu_0(\eta)))-\alpha_2(t(\mu_0(\eta)),u(\mu_0(\eta)))\\
& =a_2-a_1=mT+T_2+\theta T_1\in(mT+T_2,(m+1)T).
\end{align*}
Thus
when $\alpha_1(t(\mu_0(\eta)),u(\mu_0(\eta)))=t(\mu_1(\eta))=-nT$ we have
\begin{align*}
\alpha_2(t(\mu_0(\eta)), \; & u(\mu_0(\eta)))=t(\mu_2(\eta))=-(n+m)T-T_2-\theta T_1\\
&=-(n+m+1)T+(1-\theta)T_1\in[-(n+m+1)T,-(n+m+1)T+T_1]
\end{align*}
and the second delay falls in the first leg of the periodic solution. 
The condition $T_2+\theta T_1<T_1$ which is implied by the conditions of
Theorem~\ref{thm2delI} ensures that when the second delay satisfies
$\alpha_2(t(\mu_0(\eta)),u(\mu_0(\eta)))=t(\mu_2(\eta))=-(n+m)T$
the first delay satisfies
$\alpha_1(t(\mu_0(\eta)),u(\mu_0(\eta)))=t(\mu_1(\eta))=-nT+T_2+\theta T_1\in(-nT,-nT+T_1)$
and hence also falls in the first leg of the periodic solution. 

\begin{theorem} \label{thm2delI}
Let $K_1>1>K_2>0$ and define
\be \label{T2delI}
T=\frac{a_1(1+K_1+K_2)+(a_2-a_1)(1-K_1)}{1-m(K_1-1)+n(1+K_1+K_2)},
\ee
and
\be \label{T1T2typeI}
T_2=\frac{a_1G_{nm}(K_1)}{1-m(K_1-1)+n(1+K_1+K_2)}
\ee
where $G_{nm}(K_1)$ is defined by \eq{Gnm},
and $T_1=T-T_2$. Let the parameters be chosen so that $T_2>0$,
\be \label{theta2delI}
\theta:=\frac{a_2-a_1-mT-T_2}{T_1},
\ee
satisfies $\theta\in(0,1-1/K_1)$
and
\be \label{tIthbds}
\frac{K_2}{K_1-1}T_2<\theta T_1<T_1-\frac{1}{K_2}T_2,
\ee
then \eq{eps2del} has a Type I bimodal singular solution of period $T>0$ with solution
profile given by \eq{eqtypeI}.
\end{theorem}

\bpf
Note that the upper bound on $\theta T_1$ implies that
$T_2< K_2(1-\theta)T_1$. Hence $0<T_2<T_1<T_1+T_2=T$, and
$T_2+\theta T_1< T_1$ (since $K_2<1$). It is also useful to notice that
\eq{T2delI} can be rearranged as
\be \label{T2delsimp}
T=(a_1-nT)(1+K_1+K_2)+(1-K_1)(a_2-a_1-mT)
\ee
and \eq{theta2delI} as
\be \label{theta2delIsimp}
a_2-a_1-mT=T_2+\theta T_1.
\ee
Now for $j\in\Z$ let the functions $\mu_i(\eta)$ for $i=0,1,2$ be defined by
Table~\ref{tab:etatypeI} where
\be \label{s1j}
s_{11}=\frac{\theta T_1}{T_2+\theta T_1},\; s_{12}=1-\frac{T_2}{(1\!-\!\theta)T_1}, \;
s_{13}=1-\frac{1}{K_1(1\!-\!\theta)}, \;  
s_{14}=1-\frac{T_2}{K_2(1\!-\!\theta)T_1}. 
\ee
Clearly $s_{11}\in(0,1)$,
while $1>s_{12}>s_{14}>0$, where the last inequality follows from the upper bound on $\theta T_1$
in \eq{tIthbds}.
The bound $\theta<1-1/K_1$ also implies that $s_{13}\in(0,1)$.
It follows that
each $\mu_i(\eta)$ is continuous and monotonically
increasing. Moreover for $\eta\in[10j+k,10j+k+1]$ with $k$ a non-negative single digit integer
each function $\mu_i(\eta)$ is linear in $\eta$ with range contained in an interval
on which $u(\mu)$ and $t(\mu)$ defined by \eq{eqtypeI} are linear. It follows that
$t(\mu_i(\eta))$ and $u(\mu_i(\eta))$ are linear functions of $\eta$ for $\eta\in[10j+k,10j+k+1]$,
for integers $j$ and non-negative single digit integers $k$, as illustrated in the colour version
of Fig.~\ref{figtypeI} with $j=0$. It then follows that $F(\eta)$ is
linear on each subinterval $\eta\in[10j+k,10j+k+1]$. It is straightforward to verify that
\eq{singdel} holds, that is $t(\mu_i(\eta))=t(\mu_0(\eta))-a_i-cu(\mu_0(\eta))$ for $i=1,2$
for all $\eta\in[10j,10(j+1)]$ and hence for all $\eta\in\R$.

It remains only to verify the conditions on $F$. First note that $J^*=\bigcup_{j}[10j,10j+1]\cup[10j+5,10j+6]$,
which defines the intervals on which $t(\mu_0(\eta))$ is non-constant.
Note also that $t(\mu_1(\eta))=(j-n)T$ and
$t(\mu_2(\eta))=(j-n-m)T-T_2-\theta T_1$ for all $\eta\in[10j,10j+1]$, while
$t(\mu_2(\eta))=(j-n-m)T$ and
$t(\mu_1(\eta))=(j-n)T+T_2+\theta T_1\in((j-n)T,(j-n)T+T_1)$
for all $\eta\in[10j+5,10j+6]$.

When $\eta=10j+6$ we have
\begin{gather*}
u(\mu_0(10j+6))=u(4j+3)=\frac1c(-a_1+nT+(1-\theta)T_1),\\
u(\mu_1(10j+6))=u(4j-4n+\theta+T_2/T_1)=\frac1c(-a_1+nT+T_2+\theta T_1),\\
u(\mu_2(10j+6))=u(4j-4n-4m)=\frac1c(-a_1+nT).
\end{gather*}
Hence
\begin{align*}
F(10j+6)&=-u(\mu_0(10j+6))-K_1u(\mu_1(10j+6))-K_2u(\mu_2(10j+6))\\
&=-\frac1c(-a_1\!+\!nT\!+\!(1\!-\!\theta)T_1)-\frac{K_1}{c}(-a_1\!+\!nT\!+\!T_2\!+\!\theta T_1)-\frac{K_2}{c}(-a_1\!+\!nT).
\end{align*}
Thus, since $T=T_1+T_2$,
$$cF(10j+6)=(a_1-nT)(1+K_1+K_2)-T+(1-K_1)(T_2+\theta T_1),$$
and $F(10j+6)=0$ using \eq{T2delsimp} and \eq{theta2delIsimp}. Similarly, using \eq{s1j},
$$cF(10j+5)=cF(10j+6)+T_2-K_2(1-\theta)(1-s_{13})T_1=0.$$
Hence by the linearity of $F$ we have $F(\eta)=0$ for all $\eta\in[10j+5,10j+6]$.
Also
\begin{align*}
F(10j+1)&=-u(\mu_0(10j+1))-K_1u(\mu_1(10j+1))-K_2u(\mu_2(10j+1))\\
&=-\frac1c(-a_1+nT+T_1)-\frac{K_1}{c}(-a_1+nT)-\frac{K_2}{c}(-a_1+nT+(1-\theta)T_1).
\end{align*}
Thus
\begin{align*}
cF(10j+1)&=(a_1-nT)(1+K_1+K_2)-T_1-K_2(1-\theta)T_1,\\
& = (a_1-nT)(1+K_1+K_2)-T_1-K_2T+K_2(a_2-a_1-mT),
\end{align*}
and by \eq{T2delsimp} we find that $F(10j+1)=0$.
Similarly, using \eq{s1j},
$$cF(10j)=cF(10j+1)+T_1\bigl[1-K_1(1-\theta)(1-s_{11})\bigr]=0.$$
The linearity of $F(\eta)$ for $\eta\in[10j,10j+1]$, now ensures that
$F(\eta)=0$ for all $\eta\in J^*$.

It remains to show that $F(\eta)<0$ for all $\eta\in J^-$. Again, calling on the linearity
of $F$ on each subinterval, it is sufficient to show that $F(10j+k)<0$
for $k=2,3,4,7,8,9$. But
\begin{gather*}
cF(10j+2)=cF(10j+1)-\theta T_1[K_1+K_2-1]<0,\\
cF(10j+3)=cF(10j+1)+(1-K_1)\theta T_1+K_2T_2=(1-K_1)\theta T_1+K_2T_2,\\
cF(10j+4)=cF(10j+3)-(K_1+K_2-1)T_2<cF(10j+3),\\
cF(10j+7)=cF(10j+6)-(K_1+K_2-1)((1-\theta)T_1-T_2)<0,\\
cF(10j+8)=cF(10j+3)-(K_1-1)(T_1-T_2)-\theta T_1(1+2K_2)<cF(10j+3),\\
cF(10j+9)=cF(10j+8)-(K_1+K_2-1)T_2)<cF(10j+8),
\end{gather*}
which establishes all the required conditions if $F(10j+3)<0$, but this
holds because
$\theta>\frac{K_2}{K_1-1}\frac{T_2}{T_1}$, which completes the proof.
\epf



Next we identify the conditions on the parameters for
a type II bimodal solution to exist. In Theorem~\ref{thmIIbif} we find parameter
ranges for which these conditions are satisfied. The integers $n$ and $m$
have the same geometric meaning as for the unimodal and type I bimodal solutions and
hence satisfy \eq{nm}--\eq{m}. For the type II bimodal solution we let
$a_2-a_1=mT+\theta T_2$ where $m\in\N_0$, $\theta\in(0,1)$ and $T=T_1+T_2$
so
\begin{align*}
t(\mu_1(\eta))-t(\mu_2(\eta))
&=\alpha_1(t(\mu_0(\eta)),u(\mu_0(\eta)))-\alpha_2(t(\mu_0(\eta)),u(\mu_0(\eta)))\\
&=a_2-a_1=mT+\theta T_2\in(mT,mT+T_2).
\end{align*}
Thus
when $\alpha_1(t(\mu_0(\eta)),u(\mu_0(\eta)))=t(\mu_1(\eta))=-nT$ we have
\begin{align*}
\alpha_2(t(\mu_0(\eta)),\; & u(\mu_0(\eta)))=t(\mu_2(\eta))=-(n+m)T-\theta T_2\\
&=-(n+m+1)T+T_1+(1-\theta)T_2\in[-(n+m+1)T+T_1,-(n+m)T]
\end{align*}
and the second delay falls in the second leg of the periodic solution, as indicated in Fig.~\ref{figtypeII}.
The condition $\theta T_2<T_1$ which is implied by the conditions of
Theorem~\ref{thm2delII} ensures that when the second delay satisfies
$\alpha_2(t(\mu_0(\eta)),u(\mu_0(\eta)))=t(\mu_2(\eta))=-(n+m)T$
the first delay satisfies
$\alpha_1(t(\mu_0(\eta)),u(\mu_0(\eta)))=t(\mu_1(\eta))=-nT+\theta T_2\in(-nT,-nT+T_1)$
and hence falls in the first leg of the periodic solution, as illustrated in Fig.~\ref{figtypeII}.

\begin{theorem} \label{thm2delII}
Let $K_1+K_2>1>K_2>0$,
let $T$ be defined by \eq{T2delI}, let
\be \label{T2typeII}
T_2=\frac{a_1H_{nm}(K_1)}{1-m(K_1-1)+n(1+K_1+K_2)},
\ee
where
\begin{align} \notag
H_{nm}(K_1)&=\bigl[m-n(A-1)\bigr](K_1+K_2+1)(K_1+2K_2-1) \\
&\mbox{}\quad-K_1\bigl[(A-1)(1+K_2)-K_2\bigr]+K_2(1+K_2)+(A-1)(1-K_2) \label{Hnm}
\end{align}
and $T_1=T-T_2$. Let the parameters be chosen so that $T_2>0$ and $\theta\in(0,1)$ where
\be \label{theta2delII}
\theta=\frac{a_2-a_1-mT}{T_2},
\ee
satisfies
$\theta<\frac{T_1}{T_2}+1-\frac{1}{K_2}$ and if $K_1\geq1$
\be \label{cF6}
\theta<\left(1-\frac{1}{K_1+K_2}\right)\frac{T_1}{T_2},
\ee
or if $K_1<1$ then
\be \label{s23cF8}
\theta<\min\left\{1-\frac{(1-K_1)}{K_1}\frac{T_1}{T_2},\Bigl(1-\frac{K_1}{2K_1+K_2-1}\Bigr)\frac{T_1}{T_2}\right\}.
\ee
Then \eq{eps2del} has a Type II bimodal singular solution of period $T>0$ with solution
profile given by \eq{eqtypeII}.
\end{theorem}

\bpf
This proof is similar to the proof of Theorem~\ref{thm2delI}, differing only in the details and conditions, due
to differences in the solution profiles and parameterisations. First note that \eq{T2delsimp} is also valid for
this solution, while \eq{theta2delII} can be rewritten as
\be \label{theta2delIIsimp}
a_2-a_1-mT=\theta T_2.
\ee
Also $\theta<\frac{T_1}{T_2}+1-\frac{1}{K_2}$ implies that $\theta T_2<T_1$ and hence $0<T_1<T$.
Now, for $j\in\Z$ define $\mu_i(\eta)$ by Table~\ref{tab:etatypeII}
where
\be \label{s2j}
\begin{split}
s_{21}=\frac{T_1-\theta T_2}{T-\theta T_2}, \;
s_{22}&=\frac{T-T_2}{T-\theta T_2}, \;
s_{23}=1-\frac{T_1}{K_1(T_1+(1-\theta)T_2)}, \\
s_{24}&=1-\frac{T_2}{K_2(T_1+(1-\theta)T_2)}.
\end{split}
\ee
Then $1>s_{22}>s_{21}>0$ and
clearly $s_{23}<1$ and $s_{24}<1$. If $K_1\geq1$ then $s_{23}>0$, while if $K_1<1$ we require
$\theta<1-\frac{(1-K_1)T_1}{K_1T_2}$ for $s_{23}>0$.
Finally $\theta<\frac{T_1}{T_2}+1-\frac{1}{K_2}$ implies that $s_{24}>0$. Under these conditions $s_{2j}\in(0,1)$
for all $j$ and it follows that
each $\mu_i(\eta)$ is continuous and monotonically
increasing. Moreover for $\eta\in[10j+k,10j+k+1]$ with $k$ a non-negative single digit integer
each function $\mu_i(\eta)$ is linear in $\eta$ with range contained in an interval
on which $u(\mu)$ and $t(\mu)$ defined by \eq{eqtypeII} are linear, as illustrated in Fig.~\ref{figtypeII}.
It follows that $F(\eta)$ is
linear on each subinterval $\eta\in[10j+k,10j+k+1]$. It is straightforward to verify that
\eq{singdel} holds, that is $t(\mu_i(\eta))=t(\mu_0(\eta))-a_i-cu(\mu_0(\eta))$ for $i=1,2$.

It remains only to verify the conditions on $F$. First note that $J^*=\bigcup_{j}[10j,10j+1]\cup[10j+3,10j+4]$.
Now,
\begin{align*}
F(10j+4)&=-u(\mu_0(10j+4))-K_1u(\mu_1(10j+4))-K_2u(\mu_2(10j+4))\\
&=-\frac1c(-a_1\!+\!nT\!+\!T_1\!+\!(1\!-\!\theta)T_2)-\frac{K_1}{c}(-a_1\!+\!nT\!+\!\theta T_2)-\frac{K_2}{c}(-a_1\!+\!nT)
\end{align*}
Hence,
$$cF(10j+4)=(a_1-nT)(1+K_1+K_2)-T_1-(1-\theta)T_2-K_1\theta T_2,$$
and $F(10j+4)=0$ using \eq{T2delsimp},\eq{theta2delIIsimp} and $T=T_1+T_2$. Similarly, using \eq{s2j},
$$cF(10j+3)=cF(10j+4)+T_2-K_2(T_1+(1-\theta)T_2)(1-s_{24})=0.$$
Hence by the linearity of $F$ we have $F(\eta)=0$ for all $\eta\in[10j+3,10j+4]$.
Similarly,
\begin{align*}
cF(10j+1)&=(a_1-nT)(1+K_1+K_2)-T_1-K_2[T_1+(1-2\theta)T_2]\\
&=(a_1-nT)(1+K_1+K_2)-T_1-K_2T+2K_2(a_2-a_1-mT),
\end{align*}
and using \eq{T2delsimp} and \eq{T2typeII} we find that $F(10j+1)=0$,
while, using \eq{s2j},
$$cF(10j)=cF(10j+1)+T_1-K_1(T_1+(1-\theta)T_2)(1-s_{23})=0.$$
The linearity of $F(\eta)$ for $\eta\in[10j,10j+1]$, now ensures that
$F(\eta)=0$ for all $\eta\in J^*$.

It remains only to show that $F(\eta)<0$ for all $\eta\in J^-$. Again, using linearity
on each subinterval, it is sufficient to show that $F(10j+k)<0$
for $k=2,5,6,7,8,9$. But
\begin{gather*}
cF(10j+2)=cF(10j+1)-(K_1+K_2-1)\theta T_2<0,\\
cF(10j+5)=cF(10j+4)-(K_1+K_2-1)(T_1-\theta T_2)<0,\\
cF(10j+6)=-(K_1+K_2-1)(T_1-\theta T_2)+K_1\theta T_2,\\
cF(10j+7)=cF(10j+8)-K_2\theta T_2<cF(10j+8),\\
cF(10j+8)=-(K_1+K_2-1)(T_1-\theta T_2)+\theta T_2,\\
cF(10j+9)=cF(10j)-S_{23}K_1(T_1+(1-\theta)T_2)<0.
\end{gather*}
Now if $K_1\geq1$ then $F(10j+8)\leq F(10j+6)<0$ by \eq{cF6}
and all required conditions are satisfied for $F(\eta)<0$ for all $\eta\in J^-$.
If $K_1<1$ then $F(10j+6)<F(10j+8)<0$ (using the right-hand inequality in
\eq{s23cF8}), and again $F(\eta)<0$ for all $\eta\in J^-$.
\epf

For type I bimodal solutions to exist we require $K_1>1$ in Theorem~\ref{thm2delI}.
This condition is used twice in an essential way in the proof of that theorem, to show
that $s_{13}>0$ and also that $F(10j+3)<0$, and so Type I bimodal solutions can only exist for $K_1>1$.
In contrast, Theorem~\ref{thm2delII} does not require $K_1>1$, and we will see examples later of
type II bimodal solutions which exist for $K_1<1$.

The type I and type II bimodal solutions were constructed so that when $\alpha_1=-nT$ the second delay
$\alpha_2$ falls in the first (type I) or second (type II) leg of the periodic solution, and for both solutions
when the second delay satisfies $\alpha_2=-(n+m)T$ the first delay satisfies $\alpha_1\in(-nT,-nT+T_1)$
so falls in the first leg of the solution. We also investigated solutions where the first delay satisfies
$\alpha_1\in(-nT+T_1,-(n-1)T)$ when $\alpha_2=-(n+m)T$ and so $\alpha_1$ falls in the second leg of the solution.
However, we did not find examples of such solutions on the branches, so will not present them here.

\section{Bifurcation Branches} \label{secbifs}

Theorems~\ref{thm2delunimod}, \ref{thm2delI} and~\ref{thm2delII} specify parameter conditions
for unimodal and Type I and II bimodal singular solutions to exist for \eq{eps2del}.
In this section we use those theorems to construct parts of the bifurcation
branches. We require $K_2<1$ to ensure that \eq{eps2del}
is well posed, while $K_1$ can be arbitrary large. Thus, it is natural to take $K_1$ as a bifurcation
parameter.

\begin{figure}[t!]
\vspace{-4ex}
\mbox{}\hspace{-1cm}\scalebox{0.52}{\includegraphics{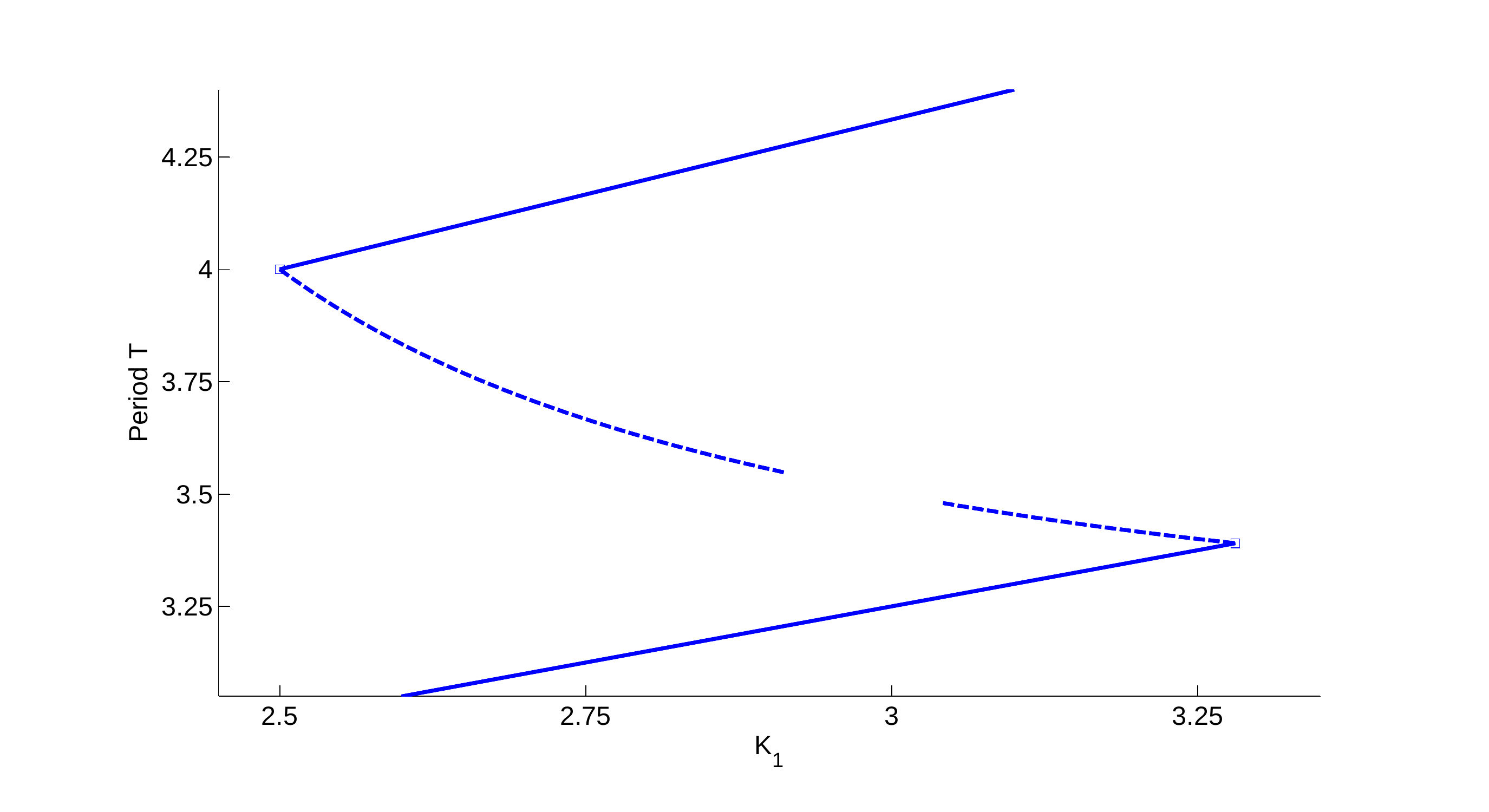}}  
\put(-230,180){\rotatebox{13}{Unimodal $n=m=0$}}
\put(-250,45){\rotatebox{9}{Unimodal $n=0$, $m=1$}}
\put(-150,88){\rotatebox{-8}{Type I, $n=0$, $m=1$}}
\put(-265,117){\rotatebox{-15}{Type II, $n=0$, $m=1$}}
\put(-347,162){$K_1=L_{00}$}
\put(-92,60){$K_1=M_{01}^+$}
\vspace{-3ex}
\caption{Example showing periods of unimodal and bimodal solutions satisfying
the conditions of Theorems~\ref{thm2delunimod}, \ref{thm2delI} and~\ref{thm2delII}
and forming a branch with two singular fold points at $K_1=L_{00}=2.5$ and $K_1=M_{01}^+\approx3.2808$
with $K_2=0.5$ and $a_2=A=5$ and $a_1=c=1$. The apparent gap in the branch near $K_1=3$
is studied in Section~\ref{secnumerics}.}
\label{figbranchex}
\end{figure}
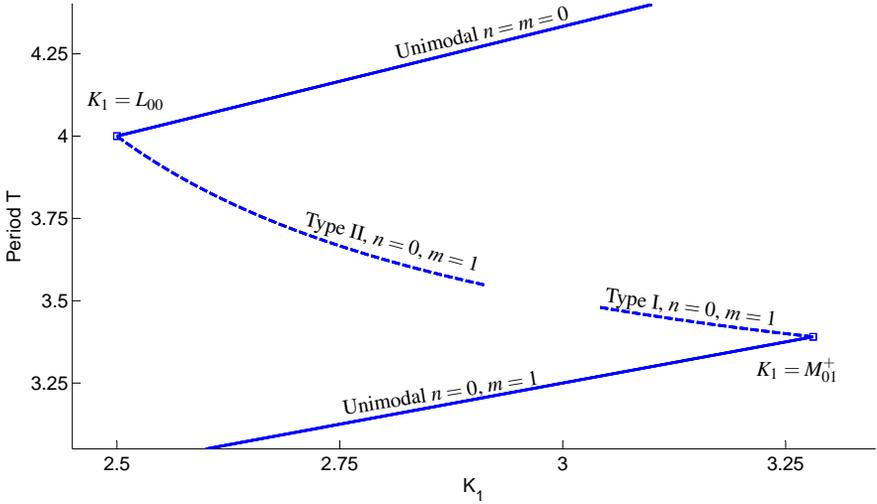

The unimodal and type I and type II bimodal solutions will be characterized by a pair of integers $(n,m)$ as
in the last section, where $n$ and $m$ are related to the delays via \eq{nm}--\eq{m}.
We will see that each value of $n$ defines a different branch of solutions,
with each branch mainly made up of segments of unimodal and type I and II bimodal singular solutions
for certain values of $m$. An example is shown in Fig.~\ref{figbranchex}. To explain
this example we need to study the parameter conditions from the three aforementioned theorems more closely.


First consider the bounds \eq{thetabds} on $\theta$ from Theorem~\ref{thm2delunimod} for the existence of unimodal solutions.
By \eq{theta2delunimod}, the bound $\theta<1$ is equivalent to $a_2-a_1<(m+1)T$. Using \eq{T2delunimod}
with $A=a_2/a_1$ this becomes
$$(A-1)(1+(m+1)K_2+n(1+K_1+K_2))<(1+m)(1+K_1+K_2+(A-1)K_2),$$
and hence $\theta<1$ is equivalent to
\be \label{unitheta1}
m-n(A-1)>-1+\frac{A-1}{1+K_1+K_2}.
\ee
We already showed that the bound $\theta>K_2/(K_1+K_2-1)$ can be written as $G_{nm}(K_1)<0$
where $G_{nm}(K_1)$ is defined by \eq{Gnm}. Notice that both bounds only depend on $n$ and $m$
through the common term $m-n(A-1)$. Let us consider the possible values of $m$ and $K_1>1$ that
satisfy these inequalities for fixed values of the other parameters. First define
\be \label{m*}
m^{*\!}(n)=n(A-1)+\frac{A-3-K_2}{2+K_2},
\ee
then when $m\geq m^{*\!}(n)$ we have
$$m-n(A-1)\geq\frac{A-3-K_2}{2+K_2}=-1+\frac{A-1}{2+K_2}>-1+\frac{A-1}{1+K_1+K_2},$$
and hence the bound $\theta<1$ is satisfied for all $K_1>1$. If
$m\in(n(A-1)-1,m^{*\!}(n))$ we find that \eq{unitheta1} is
satisfied for $K_1>L_{nm}$ where
\be \label{Lnm}
L_{nm}:=\frac{A-1}{m-n(A-1)+1}-(1+K_2)>1.
\ee
If $m=m^{*\!}(n)$ we have $L_{nm^*}=1$.
Finally there is no unimodal solution satisfying the conditions of Theorem~\ref{thm2delunimod}
if $m\leq n(A-1)-1$, since then it is impossible to satisfy \eq{unitheta1}.

Now to establish an interval of $K_1$ parameters on which a unimodal solution exists, we need to consider
both the bounds $\theta<1$ and $\theta>K_2/(K_1+K_2-1)$ together. Let $m^0(n)$ be the unique integer
for which $m^0(n)\in(n(A-1)-1,n(A-1)]$ and let
\begin{align} \notag
m^{**}(n)=n(A-1)+&\frac12\bigl[(A-1)\bigl((1+K_2)^2-K_2\bigr)+K_2\bigr]\\
&-\frac12\sqrt{\big(1+(1+K_2)^2\big)\big(((A-1)K_2+1)^2-1\big)} \label{m**}
\end{align}
In the following
theorem we establish that for $m=m^0(n)$ there is a unimodal solution satisfying
the conditions of Theorem~\ref{thm2delunimod} for all $K_1$
sufficiently large, while for each integer $m\in(m^0(n),m^{**}(n))$,
there is a non-empty bounded interval of values
of $K_1$ for which \eq{eps2del} has a unimodal solution.

\begin{theorem} \label{thmunilegs}
Let $A=a_2/a_1>1$, $K_2\in(0,1)$, and $n\in\N_0$. Let $m^{*\!}(n)$, $m^{**}(n)$, $L_{nm}$
be defined by \eq{m*}--\eq{m**}.
When $G_{nm}(K)$ defined by \eq{Gnm} has distinct roots denote them as $M_{nm}^-<M_{nm}^+$, and let
$M_{nm}$ be the root of $G_{nm}(K)$ when it has a unique root.
\begin{enumerate}[(i)]
\item
For $m=m^0(n)\in(n(A-1)-1,n(A-1)]$ there is a unimodal singular solution satisfying
the conditions of Theorem~\ref{thm2delunimod} for all $K_1>\max\{L_{nm},M_{nm}^+\}$
for any $A>1$ if $m^0(n)<n(A-1)$ and for all $K_1>\max\{L_{00},M_{00}\}$
when $A>1+\frac{K_2}{1+K_2}$ if $m^0(n)=n(A-1)$.
\item
For each integer $m\in(n(A-1),m^{*\!}(n))$ equation \eq{eps2del} has a unimodal singular solution satisfying
the conditions of Theorem~\ref{thm2delunimod} for all $K_1\in(L_{nm},M_{nm}^+)$ where
$1<L_{nm}<M_{nm}^+<\infty$.
\item
If $m=m^{*\!}(n)>m^0(n)$, then \eq{eps2del} has a unimodal singular solution satisfying
the conditions of Theorem~\ref{thm2delunimod} for all $K_1\in(1,M_{nm^*}^+)$, where
$1=L_{nm^*}=M_{nm^*}^-<M_{nm^*}^+$.
\item
For each integer $m>m^0(n)$ with $m\in(m^{*\!}(n),m^{**}(n))$ when $A>1+\frac{K_2}{1+K_2}$
we have $1\leq M_{nm}^-<M_{nm}^+<+\infty$
and \eq{eps2del} has a unimodal singular solution satisfying
the conditions of Theorem~\ref{thm2delunimod} for all $K_1\in(M_{nm}^-,M_{nm}^+)$
\item
There is no unimodal singular solution satisfying the conditions of Theorem~\ref{thm2delunimod}
if $m<m^0(n)$ or $m\geq m^{**}(n)$.
\end{enumerate}
\end{theorem}

\bpf
First consider the case when $m=m^0(n)\in(n(A-1)-1,n(A-1))$. If $m<m^{*\!}(n)$ (which is always the case if $A\geq3+K_2$)
then we have $\theta<1$ for $K_1>L_{nm}>1$. Now consider the polynomial $G_{nm}(K_1)$.
In this case the coefficient of the quadratic term is negative, and it is easy to verify that
$G_{nm}(0)>0>G_{nm}(L_{nm})$ and hence $M_{nm}^+<L_{nm}$ and $G_{nm}(K_1)<0$ for all $K_1\geq L_{nm}$. It follows that
\eq{thetabds} is satisfied for all $K_1>L_{nm}>1$. On the other hand if $m\geq m^{*\!}(n)$
then $\theta<1$ is satisfied for all $K_1>1\geq L_{nm}$ while the coefficient of the quadratic term of $G_{nm}(K_1)$
is still negative, but now $G_{nm}(1)>0$. In this case $G_{nm}(K_1)=0$ has a unique positive root $K_1=M_{nm}^+>1$
and \eq{thetabds} is satisfied for all $K_1>M_{nm}^+>1$.

Next consider the case when $m=m^0(n)=n(A-1)$, which can only arise when $A$ is rational or when $n=0$.
In this case the quadratic term in $G_{nm}(K_1)$ vanishes and the condition $G_{nm}(K_1)<0$ becomes
\be \label{Glin}
-K_1\bigl[(A-1)(1+K_2)-K_2\bigr]+K_2(1+K_2)+(A-1)<0,
\ee
which can only be satisfied for $K_1>1$ if $A>1+\frac{K_2}{1+K_2}$. In that case \eq{Glin} is
satisfied for
\be \label{Glinbd}
K_1>M_{00}:=\frac{K_2(1+K_2)+A-1}{(A-1)(1+K_2)-K_2}.
\ee
If we also set $m=n(A-1)$ in \eq{Lnm} we obtain
\be \label{L00}
K_1>L_{00}:=(A-2-K_2).
\ee
Now there are three cases. If $A\in(1+\frac{K_2}{1+K_2},3+K_2)$ then $M_{00}>1$
and by \eq{m*} we have $m>m^{*\!}(n)$ and hence \eq{thetabds} is satisfied
for all $K_1>M_{00}>1$. If $A=3+K_2$ then $L_{00}=M_{00}=1$ and \eq{thetabds} is satisfied
for all $K_1>1$. Finally if $A>3+K_2$ we have $L_{00}>1>M_{00}$ and \eq{thetabds} is satisfied
for all $K_1>L_{00}>1$. This completes the proof of (i).

To prove (ii), first note that if $A\leq 3+K_2$ then $m^{*\!}(n)\leq n(A-1)<m^0(n)+1$ and so there is no integer
$m\in(m^0(n),m^{*\!}(n))$ and nothing to prove. If $A>3+K_2$ then $m^{*\!}(n)>n(A-1)$ and for
$m\in(n(A-1),m^{*\!}(n))$ the bound $\theta<1$ is satisfied for all $K_1>L_{nm}>1$.
Moreover we find that $G_{nm}(L_{nm})<0$, while the coefficient of the quadratic term
is positive so \eq{thetabds} is satisfied for all $K_1\in(L_{nm},M_{nm}^+)$,
where $M_{nm}^+$ is the largest root of $G_{nm}(K_1)=0$.

In cases (iii) and (iv) we have $m>n(A-1)$ and $m\geq m^{*\!}(n)$.
The bound $\theta<1$ is satisfied for all $K_1>1$, while
$$G_{nm}(1)=K_2(2+K_2)(m-n(A-1))-K_2(A-3-K_2)\geq0$$
since $m\geq m^{*\!}(n)$. If and only if $K_1^*>1$ and $G_{nm}(K_1^*)<0$ where $G_{nm}'(K_1^*)=0$
there will exist a nonempty interval $(M_{nm}^-,M_{nm}^+)$ such that $1\leq M_{nm}^-<K_1^*<M_{nm}^+$,
$G_{nm}(M_{nm}^\pm)=0$ and \eq{thetabds} is satisfied for all $K_1\in(M_{nm}^-,M_{nm}^+)$.
But
$$G_{nm}'(K_1)=2(m-n(A-1))(K_1+K_2)-((A-1)(1+K_2)-K_2)$$
implies that
$$K_1^*=-K_2+\frac{(A-1)(1+K_2)-K_2}{2(m-n(A-1))},$$
and $K_1^*>1$ if and only if
\be \label{K1*>0}
m<n(A-1)+\frac12(A-1)-\frac{K_2}{2(1+K_2)}.
\ee
To establish (iii) note that $m=m^{*\!}(n)$ implies both \eq{K1*>0} and
$G_{nm}(1)=0$, thus $L_{nm^*}=M_{nm^*}^-=1<M_{nm^*}^+$.
Moreover $m>m^0(n)$ implies $m>n(A-1)$ so the quadratic term in $G_{nm}(K)$ has a positive coefficient,
and the conditions of Theorem~\ref{thm2delunimod} are satisfied for all $K_1\in(1,M_{nm^*}^+)$.

To establish (iv) let $\alpha=m-n(A-1)$ and $\beta=(A-1)(1+K_2)-K_2$.
The condition $A>1+\frac{K_2}{1+K_2}$ implies that $\beta>0$ while $m>m^0(n)$ implies
$m>n(A-1)$ and hence $\alpha>0$. The condition \eq{K1*>0} for
$K_1^*>1$ can be rewritten as $\alpha<\frac{\beta}{2(1+K_2)}$, and we also find that
\be \label{Gnm(K1*)}
G_{nm}(K_1^*)=-\frac{1}{\alpha}\left(\frac{\beta}{2}-\alpha(1+K_2)\right)^2
+K_2\bigl[(1+\alpha)(2+K_2)-(A-1)\bigr].
\ee
Then for $\alpha\in(0,\frac{\beta}{2(1+K_2)}]$ we see that $G_{nm}(K_1^*)$ is a strictly
monotonically increasing function of $\alpha$ with $G_{nm}(K_1^*)>0$ when $\alpha=\frac{\beta}{2(1+K_2)}$.
Moreover, $\lim_{\alpha\to0}G_{nm}(K_1^*)=-\infty$ and also
$G_{nm}(K_1^*)<0$ when $\alpha=\alpha^*=\frac{A-3-K_2}{2+K_2}$ (that is when $m=m^{*\!}(n)$), since
then $G_{nm}(1)=0>G_{nm}(K_1^*)$.
It follows that there
exists $\alpha^{**}\in\bigl(\max\{0,\alpha^*\},\frac{\beta}{2(1+K_2)}\bigr)$
such that $G_{nm}(K_1^*)<0$ and $K_1^*>1$
for all $\alpha\in(0,\alpha^{**})$ and $G_{nm}(K_1^*)\geq0$ and/or $K_1^*\leq1$ when
$\alpha\geq\alpha^{**}$.
Part (iv) follows on noting that $m=\alpha+n(A-1)$,
so $m^{**}(n)=\alpha^{**}+n(A-1)$. The formula \eq{m**} for $m^{**}(n)$ follows
from \eq{Gnm(K1*)} on noting that $\alpha G_{nm}(K_1^*)$ is quadratic in $\alpha$, and that $\alpha^{**}$
is given by the smaller root of $\alpha G_{nm}(K_1^*)=0$.

Finally to prove (v), note that $m<m^0(n)$ implies
$m\leq n(A-1)-1$, in which case it is not possible to satisfy \eq{unitheta1}, and
there is no unimodal solution satisfying the conditions of Theorem~\ref{thm2delunimod}.
The case of $m>m^{**}(n)$ was taken care of in the previous paragraph.
\epf

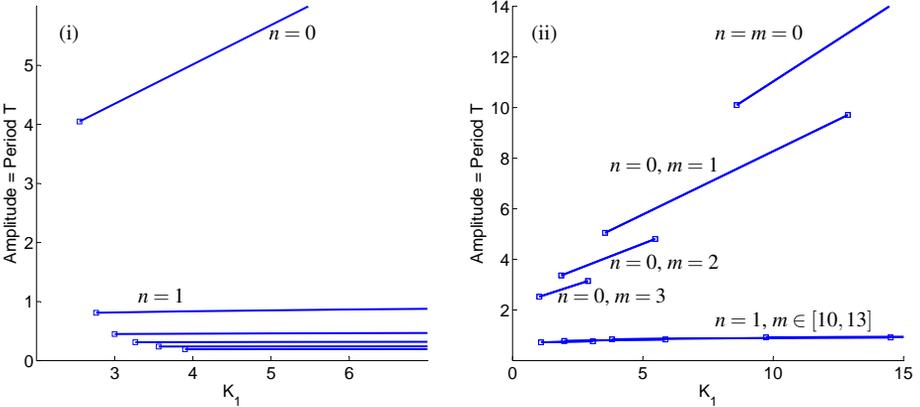
\begin{figure}
\vspace{-2ex}
\mbox{}\hspace{-0.4cm}\scalebox{0.38}{\includegraphics{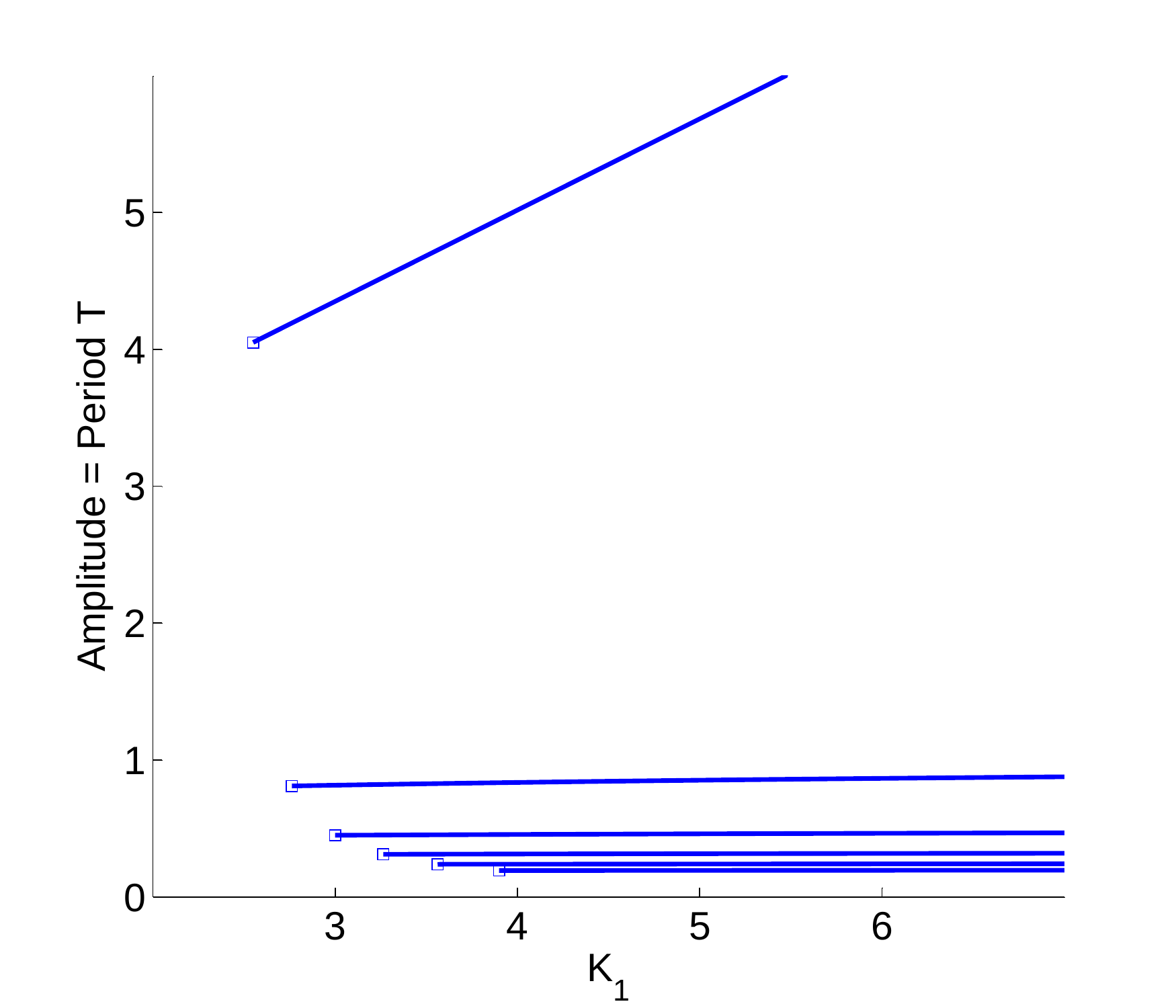}\hspace{-1cm}\includegraphics{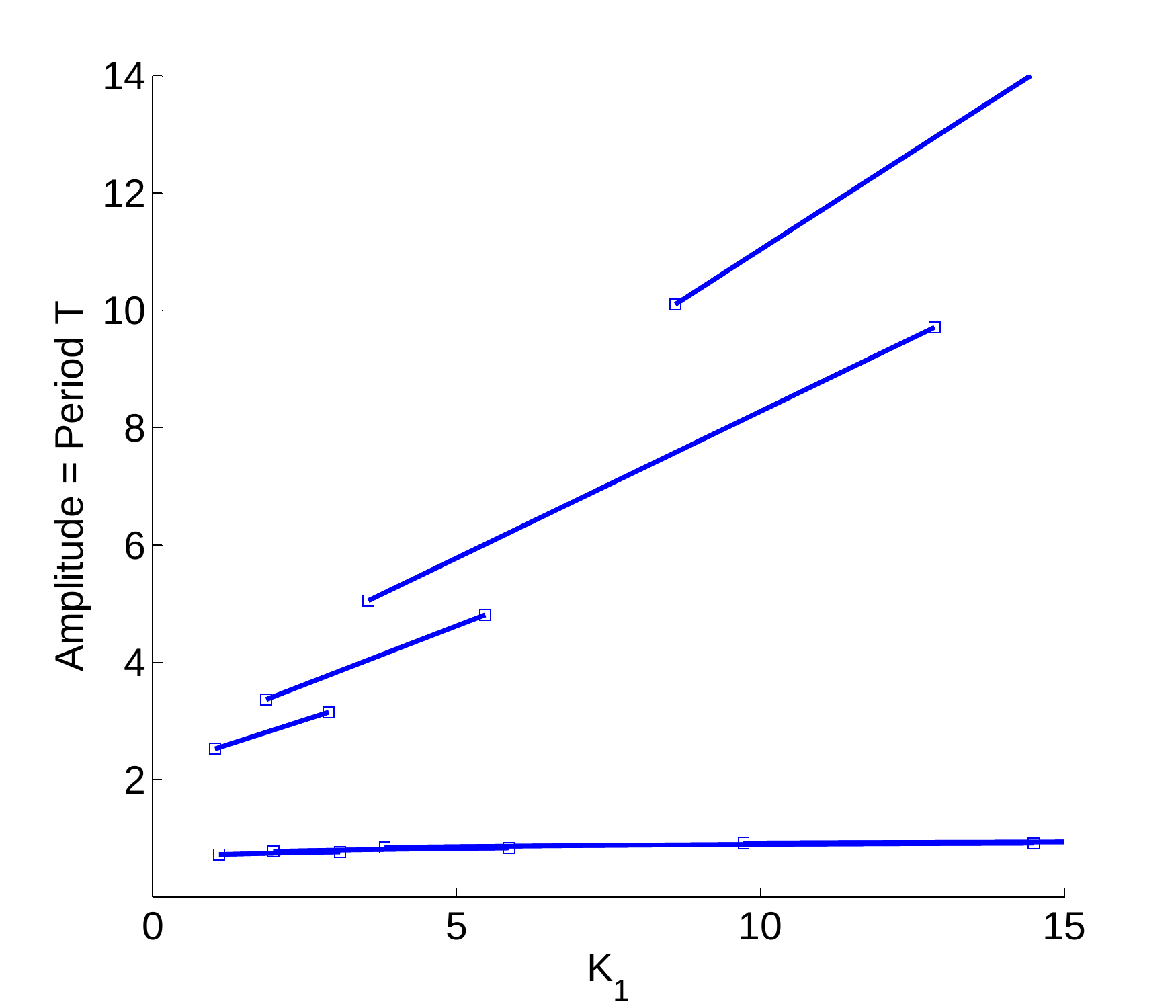}}
\put(-340,140){\footnotesize(i)}
\put(-160,140){\footnotesize(ii)}
\put(-260,140){\footnotesize$n=0$}
\put(-310,40){\footnotesize$n=1$}
\put(-90,140){\footnotesize$n=m=0$}
\put(-130,90){\footnotesize$n=0$, $m=1$}
\put(-130,53){\footnotesize$n=0$, $m=2$}
\put(-150,40){\footnotesize$n=0$, $m=3$}
\put(-90,31){\footnotesize$n=1$, $m\in[10,13]$}
\caption{(i) Periods of unimodal solutions for $m=m^0(n)$, $n=0,1,2,3,4,5$ (with decreasing period as $n$ increases)
satisfying Theorem~\ref{thmunilegs}(i), valid for all $K_1$ sufficiently large, with
$a_2=A=5.05$. (ii) Periods of legs of solutions satisfying the conditions of Theorems~\ref{thm2delunimod} and~\ref{thmunilegs}
for $a_2=A=11.1$, and
$n=0$ with $m=0,1,2,3$, and for $n=1$ with $m=10,11,12,13$. In both cases
$K_2=0.5$ and $a_1=c=1$.}
\label{figlegs}
\end{figure}

In Theorem~\ref{thmunilegs}(i)
we have shown that for $m=m^0(n)$, the smallest value of $m$ for which
a unimodal solution exists, the resulting solution exists for all $K_1$ sufficiently large.
This holds for each integer $n\geq0$ and hence, as illustrated in
Fig.~\ref{figlegs}(i), we have found the far end of infinitely
many solution branches. We note from \eq{T2delunimod} that the period $T$ increases linearly with $K_1$ on the
first ($n=0$) branch, but that for $n>0$ we have $\lim_{K_1\to\infty}T=a_1/n$.



The remainder of this work is devoted to the extension and study of these bifurcation branches
as well as their persistence for $\epsilon>0$.
Most of the rest of each solution branch will be composed of legs of other unimodal
solutions (with $m>n(A-1)$) and of bimodal solutions.
Theorem~\ref{thmunilegs}(ii)-(iv) identifies the parts of the solution branch
which are composed of unimodal solutions. This is illustrated in Fig.~\ref{figlegs}(ii).



\begin{figure}
\vspace{-4.5ex}
\mbox{}\hspace{-1cm}\scalebox{0.48}{\includegraphics{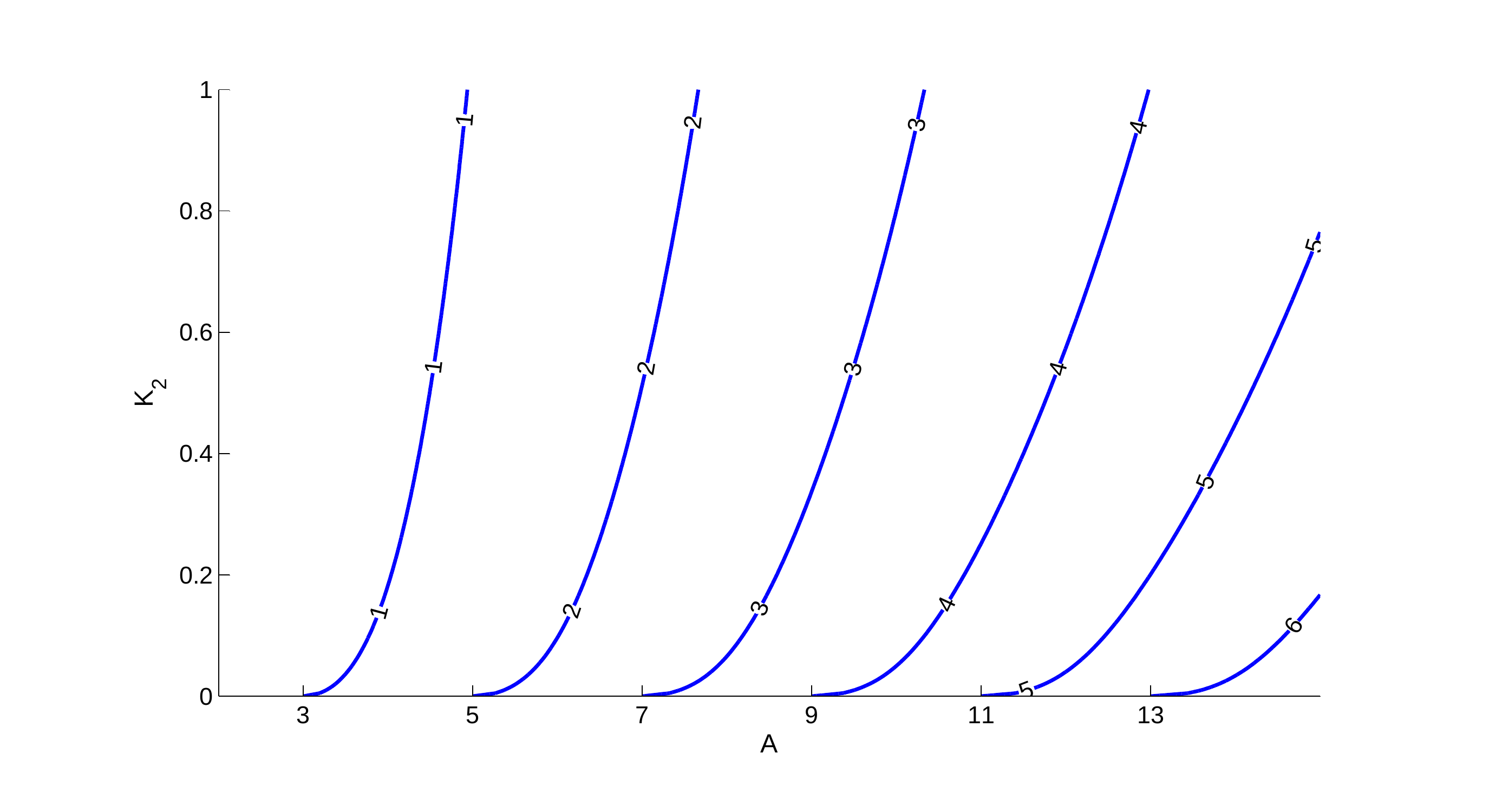}}
\vspace{-6ex}
\caption{Dependence of $m^{**}(n)$ given by \eq{m**} on $A=a_2/a_1$ and $K_2$ for the
principal $n=0$ branch. $m^{**}(0)>k$ to the right of each curve labelled $k$, and Theorem~\ref{thmunilegs}
then ensures there is a unimodal solution defined on an interval of $K_1$ values for $n=0$ and each
integer $m=0,1,\ldots,k$.}
\label{figMstst}
\end{figure}

From \eq{T2delunimod} we see that the unimodal solutions with the largest period occur on the
branch $n=0$. Let us consider this branch further. From Theorem~\ref{thmunilegs}(i) provided $A=a_2/a_1>1+\frac{K_2}{1+K_2}$
there is a leg of unimodal solutions for $n=m=m^{0}(0)=0$ for all $K_1>\max\{L_{00},M_{00}\}$.
In that case Theorem~\ref{thmunilegs} also ensures there will be legs of unimodal periodic solutions for
each integer $m$ between $0$ and $m^{**}(0)$. Hence we require $m^{**}(0)>1$ for there be a second leg
of unimodal solutions with $n=0$ and $m=1$. Fig.~\ref{figMstst} shows the dependence
of $m^{**}(0)$ on $A$ and $K_2$, from which we see that we require the ratio $A=a_2/a_1>3$
for there to be a second, $m=1$, leg of unimodal solutions for $K_2$ sufficiently small, while
for $A>5$ there is an $m=1$ leg of unimodal solutions for all $K_2\in(0,1)$. Arbitrary large
values of $m^{**}(0)$ are possible but require $A\gg1$. We will explore the case $A\gg1$ in Section~\ref{secnumerics}.
For other branches of solutions with $n>0$, note that $m^{0}(n)\in(m^{0}(0)+n(A-1)-1,m^{0}(0)+n(A-1)]$ and
from \eq{m**} we have $m^{**}(n)=m^{**}(0)+n(A-1)$, and so for fixed $A$ and $K_2$ essentially the same number of
legs of unimodal solutions appear for each value of $n$, but the corresponding values of $m$ are shifted by $n(A-1)$.


To show for a given value of $n$ that the legs for different values of $m$ form part of a connected branch
of solutions we need to join up the branches, which we will do using bimodal periodic solutions of type I and
II, and multi-modal solutions.
First we note that if there is a continuous branch for fixed $n$ then for $A$ sufficiently large it must
have fold bifurcations. To see this note that if $m>m^0(n)$ then the coefficient of the quadratic term in
\eq{Gnm} is positive, while by \eq{Lnm}
$$G_{nm}(L_{nm-1})=K_2(A-1)\left[2+K_2-\frac{A-2}{m-n(A-1)}\right].$$
Hence $G_{nm}(L_{nm-1})<0$ if $m<m^{*\!}(n)+\frac{1+K_2}{2+K_2}$. But since, $M_{nm}^+$ is the largest
zero of $G_{nm}$ this implies that $L_{nm-1}<M_{nm}^+$. By \eq{Lnm} we also have $L_{nm}<L_{nm-1}$.
Thus in the case of legs for $m=m^0(n),\ldots,M$
with $M<m^{*\!}(n)+\frac{1+K_2}{2+K_2}$ we have that $L_{nm-1}\in(L_{nm},M_{nm}^+)$ for $m=m^0(n)+1,\ldots,M$.
This implies that the $K_1$ values of legs for adjacent values of $m$ overlap, just as illustrated in
Fig.~\ref{figlegs}(ii). Hence if the legs form part of a continuous branch, the branch must have folds,
as seen in Fig.~\ref{figbranchex}. As Fig.~\ref{figbranchex} suggests, these will not be smooth fold bifurcations
in the classical sense, but we will see in Section~\ref{secnumerics} that they are the $\epsilon\to0$ singular limit of
smooth fold bifurcations of periodic orbits, and so we will refer to them as fold bifurcations anyway.
We will see that these singular fold bifurcations typically occur at $K_1=M_{nm}^+$ and $K_1=L_{nm}$.
To show this we need the following lemma which determines the sign of the denominator
of \eq{T2delI}, \eq{T1T2typeI} and \eq{T2typeII}.

\begin{lemma} \label{lemsk1}
Let
\be \label{sk1}
s(K_1):=1-m(K_1-1)+n(1+K_1+K_2).
\ee
If $m\in(m^0(n),m^{*\!}(n)+\frac{1+K_2}{2+K_2})$ then $s(L_{nm-1})<0$. Moreover,
if $n=0$ or if $n>0$ and $A\geq2-1/n$ then $s(K_1)<0$ for all $K_1\geq L_{nm-1}$.
Finally, if $m>m^{*\!}(n)+\frac{1+K_2}{2+K_2}$ then $s(L_{nm-1})>0$.
\end{lemma}

\bpf
First note that from \eq{Lnm}
$$s(L_{nm-1})=\frac{m\bigl[2-A+(2+K_2)(m-n(A-1)\bigr]}{m-n(A-1)}.$$
Now $m>m^0(n)$ implies that $m-n(A-1)>0$, while
$m<m^{*\!}(n)+\frac{1+K_2}{2+K_2}$ implies
$m-n(A-1)<(A-2)/(2+K_2)$ and hence $[2-A+(2+K_2)(m-n(A-1)]<0$ which shows that $s(L_{nm-1})<0$.
If $s(K_1)$ is a nonincreasing function of $K_1$, then it follows that
$s(K_1)<0$ for all $K_1\geq L_{nm-1}$. But this is trivially true
if $m\geq n$, which is always the case when $n=0$, since $m\in\N_0$.
For $n>0$,
provided $A\geq2-1/n$, we have $m\geq m^0(n)+1>n(A-1)\geq n-1$ and hence $m\geq n$.
Finally $m>m^{*\!}(n)+\frac{1+K_2}{2+K_2}$ implies
$[2-A+(2+K_2)(m-n(A-1)]>0$ which shows that $s(L_{nm-1})>0$.
\epf

The following theorem establishes the existence of a fold bifurcation of periodic
singular solutions at $K_1=M_{nm}^+$. As noted before Lemma~\ref{lemsk1}, this
will not be a smooth bifurcation, but rather a leg of unimodal solutions and a leg of
type I bimodal solutions will both exist for $K_1\in(M_{nm}^+-\delta,M_{nm}^+)$ and
these solutions will coincide in the limit as $K_1\to M_{nm}^+$. By coincide, we mean
that the limiting profiles and periods of both solutions will be identical.

\begin{theorem} \label{thmIbif}
Let $A=a_2/a_1>1$, $K_2\in(0,1)$, $n\in\N_0$ and $m^\dagger(n)=m^{*\!}(n)+\min\Bigl\{\frac{1+K_2}{2+K_2},1-\frac{A-1}{(2+K_2)(3+K_2)}\Bigr\}$.
If $m\in(n(A-1),m^\dagger(n))$
then there exists $\delta>0$ such that for $K_1\in(M_{nm}^+-\delta,M_{nm}^+)$
there is
\begin{enumerate}[a)]
\item a leg of unimodal solutions satisfying the conditions of Theorem~\ref{thm2delunimod}
\item a leg of type I bimodal solutions satisfying conditions of Theorem~\ref{thm2delI}
\end{enumerate}
and these solutions coincide at $K_1=M_{nm}^+$.
\end{theorem}

\bpf
Theorem~\ref{thmunilegs} gives the existence of a leg of unimodal solutions
for $K_1\in(L_{nm},M_{nm}^+)$ or $K_1\in(M_{nm}^-,M_{nm}^+)$ when $m\in(n(A-1),m^{*\!}(n)+\frac{1+K_2}{2+K_2})$.
Next we show that if there exists a leg
of type I bimodal solutions for $K_1\in(M_{nm}^+-\delta,M_{nm}^+)$ then
the unimodal and type I bimodal solutions must coincide in the limit as $K_1\to M_{nm}^+$.
To see this, compare the profile of the type I bimodal solution in \eq{eqtypeI}
with the profile of the unimodal solution in equations \eq{idelg1},\eq{idelg2}.
Since $G_{nm}(M_{nm}^+)=0$, by \eq{T1T2typeI} the bimodal solution must
satisfy $\lim_{K_1\to M_{nm}^+}T_2=0$. But when $T_2=0$ the bimodal profile
corresponds to the unimodal profile. Elementary algebra then shows that
the period $T$ of the unimodal solution given by \eq{T2delunimod} equals the
period $T$ given by \eq{T2delI} for the bimodal profile when $G_{nm}(K_1)=0$.

Finally we confirm the existence of the type I bimodal solution for $K_1\in(M_{nm}^+-\delta,M_{nm}^+)$
by verifying the conditions of Theorem~\ref{thm2delI}.
Since $G_{nm}(M_{nm}^+)=0$, when
$K_1=M_{nm}^+$ by \eq{T1T2typeI} we have $T_2=0$, and $T_1=T$, where the value of $T$
is given by \eq{T2delI} or \eq{T2delunimod}.
Now from \eq{theta2delI}
$$\theta=\frac{a_2-a_1-mT}{T}=\frac{K_2}{K_1+K_2-1}\in(0,1),$$
using \eq{theta2delunimod},\eq{thetabds} and the definition of $M_{nm}^+$.
Thus the bounds \eq{tIthbds} are trivially satisfied when $K_1=M_{nm}^+$.
The bound $\theta<1-1/K_1$ also holds provided $K_2<(K_1-1)^2$, in particular whenever $M_{nm}^+>2$,
but $M_{nm}^+>L_{nm-1}$ and $m<m^{*\!}(n)+1-\frac{A-1}{(2+K_2)(3+K_2)}$ implies that $L_{nm-1}>2$.

Thus all the conditions for the existence of a type I unimodal solution
from Theorem~\ref{thm2delI} are satisfied when $K_1=M_{nm}^+$, and by continuity on an interval
containing this point, except possibly for the condition $T_2>0$. But $T_2=a_1G_{nm}(K_1)/s(K_1)$ by \eq{T1T2typeI}.
Now, noting that $G_{nm}(K_1)<0$ for $K_1\in(M_{nm}^-,M_{nm}^+)$, and $G_{nm}(K_1)>0$
for $K_1>M_{nm}^+$, provided $s(M_{nm}^+)\ne0$,
the conditions of Theorem~\ref{thm2delI} must be satisfied on some interval
$(M_{nm}^+-\delta,M_{nm}^+)$ or $(M_{nm}^+,M_{nm}^++\delta)$ by continuity of $s(K_1)$.
But, by Lemma~\ref{lemsk1} we have $s(K_1)<0$ for all $K_1>L_{nm-1}$, and since
$L_{nm-1}<M_{nm}^+$ it follows that for $\delta>0$ sufficiently small that $s(K_1)<0$
for $K_1\in(M_{nm}^+-\delta,M_{nm}^++\delta)$.
Thus there is a unimodal solution for $K_1\in(L_{nm},M_{nm}^+)$
and a bimodal solution on an interval $(M_{nm}^+-\delta,M_{nm}^+)$ which coincide
at a fold bifurcation at $K_1=M_{nm}^+$.
\epf

For values of $m$ outside the range for which Theorem~\ref{thmIbif} is valid, it can still be possible to
obtain type I bimodal and unimodal solutions which coincide at $K_1=M_{nm}^+$ without a fold bifurcation.
An example of this will be seen later in Fig.~\ref{figamp76}. We will not determine here the size of $\delta>0$
such that Theorem~\ref{thmIbif} applies. However, we note that since the theorem guarantees
the existence of the unimodal and type I bimodal solutions on some interval, it is a straightforward
task to check the conditions of Theorems~\ref{thm2delunimod} and~\ref{thm2delI} to determine the
interval on which each solution exists, and this is what we will do in later examples.

Since the proof of Theorem~\ref{thmIbif} is purely algebraic, it is interesting to consider the bifurcation
from a dynamical viewpoint. For the leg of unimodal solutions $\theta$ approaches
its lower bound in \eq{thetabds} as $K_1$ approaches $M_{nm}^+$ .
Indeed, since $M_{nm}^\pm$ are the zeros of $G$ defined by
\eq{Gnm}, it follows that $\theta\to K_2/(K_1+K_2-1)$ as $K_1\to M_{nm}^\pm$ for all of the unimodal solutions
identified in Theorem~\ref{thmunilegs}. At $K_1=M_{nm}^\pm$ we have $F(5j+3)=0$ in the proof
of Theorem~\ref{thm2delunimod}. The condition $F(5j+3)<0$ in that proof ensures that $F$ remains negative
while $u(\mu_2(\eta))$, the value of
$u$ at the second delay, decreases from its maximum value $(-a_1+(n+1)T)/c$ to its minimum
value $(-a_1+nT)/c$. If $F(\eta)=0$ for some $\eta\in(5j+2,5j+3)$ then the solution would reenter $J^*$ and we
would expect another interval on which $F(\eta)=0$. This is exactly what happens in the bifurcation to the
type I bimodal solution in Theorem~\ref{thmIbif}. For the type I bimodal solution which exists for $K_1<M_{nm}^+$,
from the proof of Theorem~\ref{thm2delI} we see that for $\eta\in(10j+4,10j+5)$ the solution at the second delay, $u(\mu_2)$
decreases from its maximum value $(-a_1+(n+1)T)/c$ to $(-a_1+nT+(1-s_{14})(1-\theta)T_1)/c$  with $F(\eta)<0$, but
$F(\eta)\to0$ as $\eta\to10j+5$. For $\eta\in(10j+5,10j+6)$ we have $\eta\in J^*$, $F(\eta)=0$ and $u(\mu_2)$
further decreases to its minimum value $(-a_1+nT)/c$. However as $K_1$ approaches $M_{nm}^+$ we have $G_{nm}(K_1)\to0$, and
hence $T_2\to0$ and $s_{14}\to1$, so $(-a_1+nT+(1-s_{14})(1-\theta)T_1)/c$ is no longer larger than $(-a_1+nT)/c$
and $F$ does not become zero before $u(\mu_2)$ reaches its global minimum. Hence the second interval of $J^*$ for
$\eta\in(10j+5,10j+6)$ collapses, and as $s_{11}$ and $s_{12}$ both tend to $1$, we find that five of the intervals of the
parametrisation of the type II bimodal solution become trivial, and the remaining parts correspond to the unimodal solution.
Thus at the bifurcation between the unimodal solution and the type I bimodal solution $\theta$ hits its lower
bound for the unimodal solution, and $T_2\to0$ for the bimodal solution. It will be interesting to investigate below what
other bifurcations arise as other conditions in the theorems of Section~\ref{secsingsols} are violated.


Now consider the case of $K_1=L_{nm-1}$ at the left-hand end of the interval of unimodal
solutions for $K_1\in(L_{nm-1},M_{nm-1}^+)$. We show that at this point there is a fold bifurcation
and the solution transforms from a unimodal solution to a type II bimodal solution. By the definition of
$L_{nm}$ for $K_1>L_{nm-1}$ the unimodal solution satisfies $a_2-a_1=(m-1+\theta)T$ with $\theta\in(0,1)$ but as
$K_1\to L_{nm-1}$ we have $\theta\to1$. But if $\theta$ were equal to $1$, the difference $a_2-a_1$ between the two
delayed times would be exactly $m$ periods. Perhaps not surprisingly, as the following
theorem shows, this can result in a (type II) bimodal solution with the value of $m$ increased by $1$.

\begin{theorem} \label{thmIIbif}
Let $A=a_2/a_1>1$, $K_2\in(0,1)$ and $n\in\N_0$.
\begin{enumerate}[i)]
\item
If $p\in[m^0(n)+1,m^{*\!}(n)+\frac{1+K_2}{2+K_2})$ then
there exists $\delta>0$ such that for $K_1\in(L_{np-1},L_{np-1}+\delta)$
there is
\begin{enumerate}[a)]
\item a leg of unimodal solutions satisfying the conditions of Theorem~\ref{thm2delunimod} with $m=p-1$,
\item a leg of type II bimodal solutions satisfying conditions of Theorem~\ref{thm2delI} with $m=p$
\end{enumerate}
and these solutions coincide at $K_1=L_{np-1}$.
\item
If $p\in(m^{*\!}(n)+\frac{1+K_2}{2+K_2},m^{*\!}(n)+1]$ then
there exists $\delta>0$ such that
\begin{enumerate}[a)]
\item for $K_1\in(L_{np-1},L_{np-1}+\delta)$
there is a leg of unimodal solutions satisfying the conditions of Theorem~\ref{thm2delunimod} with $m=p-1$,
\item for $K_1\in(L_{np-1}-\delta,L_{np-1})$
there is a leg of type II bimodal solutions satisfying conditions of Theorem~\ref{thm2delI} with $m=p$
\end{enumerate}
and these solutions coincide at $K_1=L_{np-1}$.
\item
If $p\in\Bigl(m^{*\!}(n)+1,\min\bigl(m^{**}(n)+1,(n+\frac12)(A-1)\bigr)\Bigr)$ then
$1-K_2<L_{np-1}<1<M_{np-1}^-<M_{np-1}^+$ and there exists $\delta>0$ such that
\begin{enumerate}[a)]
\item for $K_1\in(M_{np-1}^-,M_{np-1}^+)$
there is a leg of unimodal solutions satisfying the conditions of Theorem~\ref{thm2delunimod} with $m=p-1$,
\item for $K_1\in(L_{np-1}-\delta,L_{np-1})$
there is a leg of type II bimodal solutions satisfying conditions of Theorem~\ref{thm2delI} with $m=p$
\end{enumerate}
and these solutions exist on disjoint parameter intervals.
\end{enumerate}
\end{theorem}

\bpf
Theorem~\ref{thmunilegs} gives the existence of a leg of unimodal solutions
with $m=p-1$
for $K_1\in(L_{np-1},+\infty)$ when $p=m^0(n)+1$,
for $K_1\in(L_{np-1},M_{np-1}^+)$ when $p\in(m^0(n)+1,m^{*\!}(n)+1]$,
and for $K_1\in(M_{np-1}^-,M_{np-1}^+)$ when $p\in[m^{*\!}(n)+1,m^{**}(n)+1)$.

To prove (i) and (ii), next we show that if there exists a leg
of type II bimodal solutions with $m=p$ for $K_1\in(L_{np-1}-\delta,L_{np-1})$
or $K_1\in(L_{np-1},L_{np-1}+\delta)$ then
the unimodal and type II bimodal solutions must coincide in the limit as $K_1\to L_{np-1}$.
To see this, compare the profile of the type II bimodal solution in \eq{eqtypeII}
with the profile of the unimodal solution in equations \eq{idelg1},\eq{idelg2}.
The two solutions will coincide in the limit as $K_1\to L_{np-1}$ if
both the unimodal and type II bimodal solution have the same limiting period $T$ and
for the type II bimodal solution $\theta\to0$ as $K_1\to L_{np-1}$. But it is simple
to check that the value of $T$ given by \eq{T2delI} for the type II bimodal solution
with $m=p$ agrees with that given by \eq{T2delunimod} for the unimodal solution with $m=p-1$.
The rest of this proof concerns the existence and properties of the type II
bimodal solution with $m=p$, so we can use $m$ and $p$ interchangeably.
To show that $\theta\to0$ as $K_1\to L_{np-1}$ for the type II bimodal solution with $m=p$, note that by
\eq{T2delI}
$$a_2-a_1-mT=a_1\bigl(A-1-\frac{mT}{a_1}\bigr)=\frac{a_1}{s(K_1)}\bigl[A-1-(m-n(A-1))(1+K_1+K_2)\bigr],$$
and from \eq{Lnm} we have
\be \label{num0}
(L_{nm-1}+K_2+1)(m-n(A-1))=A-1
\ee
Now, $s(L_{np-1})\ne0$ by Lemma~\ref{lemsk1}, since $p\ne m^{*\!}(n)+\frac{1+K_2}{2+K_2}$.
Hence
$$\lim_{K_1\to L_{np-1}}\theta=\lim_{K_1\to L_{np-1}}\frac{1}{T_2}(a_2-a_1-mT)=0,$$
as required, provided $T_2>0$.

To derive expressions for $T_1$ and $T_2$ when $K_1=L_{nm-1}$, using \eq{T2typeII} and \eq{num0}
\begin{align*}
T_2&=\frac{a_1}{s(L_{nm-1})}H_{nm}(L_{nm-1})\\
&=\frac{a_1}{s(L_{nm-1})}
\Bigl[(A-1)(K_1+2K_2-1)\\
& \qquad\qquad\qquad-K_1\bigl[(A-1)(1+K_2)-K_2\bigr]+K_2(1+K_2)+(A-1)(1-K_2)\Bigr]\\
&=\frac{a_1K_2}{s(L_{nm-1})}\Bigl[(A-1)(1-K_1)+(1+K_1+K_2)\Bigr]=K_2T.
\end{align*}
Thus $T_2=K_2T$ and $T_1=(1-K_2)T$ when $K_1=L_{nm-1}$.
Since $K_2\in(0,1)$ this implies that $T_1>0$ and $T_2>0$ in the limit as $K_1\to L_{np-1}$,
which establishes that the unimodal solution and type II bimodal solution have the
same limiting profiles as $K_1\to L_{np-1}$.

To prove (i) and (ii), it remains to verify the conditions of Theorem~\ref{thm2delII} to confirm the existence
of the type II bimodal solution. First note that since $T_1>0$ and $T_2>0$ when $K_1=L_{nm-1}$,
and $s(L_{nm-1})\ne0$ there exists $\delta>0$ such that $T_1$ and $T_2$ defined by
\eq{T2delI},\eq{T2typeII} vary continuously and are strictly positive for $K_1\in(L_{nm-1}-\delta,L_{nm-1}+\delta)$.
Now consider the condition $\theta>0$. From above
$$\theta=\frac{a_1}{T_2}\frac{1}{s(K_1)}\bigl[A-1-(m-n(A-1))(1+K_1+K_2)\bigr].$$
Under (i) we have $s(L_{nm-1})<0$ and hence $s(K_1)<0$ for $K_1\in(L_{nm-1},L_{nm-1}+\delta)$ for
$\delta$ sufficiently small. Also by \eq{num0} for $K_1>L_{nm-1}$ we have $A-1-(m-n(A-1))(1+K_1+K_2)<0$.
Hence $\theta>0$ for $K_1\in(L_{nm-1},L_{nm-1}+\delta)$.
Similarly under (ii) we have $s(L_{nm-1})>0$ and hence $s(K_1)>0$ for $K_1\in(L_{nm-1}-\delta,L_{nm-1})$ for
$\delta$ sufficiently small, and by \eq{num0} we have $A-1-(m-n(A-1))(1+K_1+K_2)>0$ for $K_1>L_{nm-1}$.
Thus under (ii) $\theta>0$ for $K_1\in(L_{nm-1}-\delta,L_{nm-1})$. Moreover since $\theta\to0$ as
$K_1\to L_{nm-1}$ in both cases, for $\delta>0$ sufficiently small we also have $\theta<1$.

Next we show that the condition $\theta<\frac{T_1}{T_2}+1-\frac{1}{K_2}$ holds. Since
$s(K_1)\bigl[A-1-(m-n(A-1))(1+K_1+K_2)\bigr]>0$ for $K_1\in(L_{nm-1}-\delta,L_{nm-1})$ under (i)
and for $K_1\in(L_{nm-1},L_{nm-1}+\delta)$ under (ii), in both cases we have
\begin{align*}
0&<\frac{a_1}{K_2s(K_1)}(K_1+K_2-1)\bigl[(A-1)-(m-n(A-1))(1+K_1+K_2)\bigr]\\
&=\frac{m+1}{s(K_1)}\bigl[a_1(1+K_1+K_2)+(a_2-a_1)(1-K_1)\bigr]-(a_2-a_1)-\frac{a_1H_{nm}(K_1)}{K_2s(K_1)}\\
&=(m+1)T-(a_2-a_1)-\frac{1}{K_2}T_2\\
&=T-\theta T_2-\frac{1}{K_2}T_2, \qquad\mbox{since }\theta T_2=a_2-a_1-mT.
\end{align*}
Hence $\theta T_2<T_1+T_2-\frac{1}{K_2}T_2$, and since $T_2>0$
we have $\theta<\frac{T_1}{T_2}+1-\frac{1}{K_2}$ as required.

It remains only to establish \eq{cF6} or \eq{s23cF8}. But for $p<m^{*\!}(n)+1$,
since $L_{np-1}>1$ by Theorem~\ref{thmunilegs}, for $\delta>0$ sufficiently small $K_1>1$ for
all $K_1\in(L_{np-1}-\delta,L_{np-1}+\delta)$, and so only \eq{cF6} is required. But the right-hand side
of \eq{cF6} is strictly positive since $K_1>1>K_2>0$, while from above $\theta\to0$ as $K_1\to L_{nm-1}$ so
this inequality also holds for $\delta>0$ sufficiently small. On the other hand if $p=m^{*\!}(n)+1$ then
$L_{np-1}=L_{nm^*}=1$ and we need to verify \eq{s23cF8} for $K_1\in(L_{np-1}-\delta,L_{np-1})=(1-\delta,1)$.
But both expressions on the right-hand side of \eq{s23cF8} are strictly positive for all $K_1$ sufficiently
close to $1$, while we already showed that $\lim_{K_1\to L_{np-1}}\theta=0$ and so
\eq{s23cF8} is satisfied for $K_1\in(1-\delta,1)$. This establishes (i) and (ii).

To prove (iii) it remains only to establish the existence of the type II bimodal solution in that case,
but this is similar to above, where we note that $p<(n+\frac12)(A-1)$ implies $L_{np-1}>1-K_2$
and choosing $\delta$ sufficiently small so that $L_{np-1}-\delta>1-K_2$ ensures that
$K_1+K_2>1$ for all $K_1\in(L_{np-1}-\delta,L_{np-1})$. This implies that the second term on the
right-hand side of \eq{s23cF8} is strictly positive, while the first expression tends to $(L_{np-1}+K_2-1)/(L_{np-1}K_2)>0$
in the limit as $K_1\to L_{np-1}$. Again, since $\lim_{K_1\to L_{np-1}}\theta=0$,
for $\delta>0$ sufficiently small
equation \eq{s23cF8} is satisfied for $K_1\in(L_{np-1}-\delta,L_{np-1})$.
\epf

Theorem~\ref{thmIIbif}(i) establishes the existence of a fold bifurcation when $K_1=L_{np-1}$ for
$p\in\bigl(m^0(n),m^{*\!}(n)+\frac{1+K_2}{2+K_2}\bigr)$. Interestingly, for $p\in(m^{*\!}(n)+\frac{1+K_2}{2+K_2},m^{*\!}(n)+1]$
the fold disappears, but
the two legs of periodic solutions continue to exist and coincide at $K_1=L_{np-1}$, but now
the type II bimodal solution exists for $K_1<L_{np-1}$ while the unimodal solution exists for $K_1>L_{np-1}$.
Essentially, the fold bifurcation unfolds suggesting a cusp-like bifurcation of periodic orbits,
which we will investigate in Section~\ref{sec:cusp}.

Theorem~\ref{thmIIbif}(iii) also reveals interesting behaviour. When $m=m^{*\!}(n)+1$, or
equivalently,
\be \label{K1<1}
A=1+\frac{m(2+K_2)}{1+n(2+K_2)}
\ee
we have $L_{nm-1}=L_{nm^*}=1$. Noting that
$m^{*\!}(n)+1<(n+1/2)(A-1)$, for $m\in(m^{*\!}(n)+1,(n+1/2)(A-1))$, or equivalently for
$$A\in\left(1+\frac{2m}{1+2n},1+\frac{(2+K_2)m}{1+(2+K_2)n}\right),$$
we have $1>L_{nm-1}>1-K_2$ and Theorem~\ref{thmIIbif}(iii) ensures the existence
of type II bimodal solutions for $K_1\in(L_{nm-1}-\delta,L_{nm-1})$ where $1-K_2<L_{nm-1}-\delta<K_1<L_{nm-1}<1$.
In contrast the construction of the unimodal and type I bimodal solutions in Theorems~\ref{thm2delunimod}
and~\ref{thm2delI} requires $K_1>1$ for those solutions to exist.

Dynamically, we see in the proof of Theorem~\ref{thmIIbif} that for both the unimodal and type II bimodal solution
we have $T\to(a_2-a_1)/m$ as $K_1\to L_{nm-1}$. For the unimodal solution
$(a_2-a_1)=(m-1+\theta)T<mT$ and $\theta\to1$ as $K_1\to L_{nm-1}$, while
for the type II bimodal solution $(a_2-a_1)=mT+\theta T_2>mT$ and $\theta\to0$ as $K_1\to L_{nm-1}$. Whether or
not there is a (non-smooth) fold bifurcation at $K_1=L_{nm-1}$ depends on whether $m$ is greater or
smaller than $m^{*\!}(n)+(1+K_2)/(2+K_2)$. That the value of $m$ increases close to $K_1=L_{nm-1}$ was already
observed numerically for $\epsilon>0$ in \cite{DCDSA11}.

Theorem~\ref{thm2delunimod} identified upper and lower bounds on $\theta$ for a unimodal solution to exist.
In Theorems~\ref{thmIbif} and~\ref{thm2delII} we have shown bifurcations to type I or type II bimodal
solutions when one of these bounds is violated. In Section~\ref{secnumerics}  we
will investigate the solutions that can
arise when the parameters bounds identified in Theorems~\ref{thm2delI} and~\ref{thm2delII} for type I and
type II bimodal solutions are violated.

\section{Singularly Perturbed Solution Branches}
\label{secnumerics}

We are interested in solutions of \eq{epsf} when $0<\epsilon\ll1$. However, so far we have only constructed
$\epsilon=0$ singular solutions, in the sense of Definition~\ref{defsingsol}. It would be desirable
to prove that \eq{epsf} has solutions close to the constructed singular solutions for all $\epsilon$ sufficiently
small. Mallet-Paret and Nussbaum \cite{JMPRN11} proved that the sawtooth is indeed the limiting profile as
$\epsilon\to0$ for the state-dependent DDE \eq{eps1del} which has only one delay.
However, for the two delay problem \eq{eps2del}, Theorems~\ref{thmIbif} and~\ref{thmIIbif} lead
us to expect fold bifurcations of periodic orbits. Indeed such bifurcations and resulting intervals of
co-existing stable periodic solutions were already observed for $\epsilon=\cO(1)$ in \cite{DCDSA11}.
`Superstability' is central to the results of \cite{JMPRN11}, and without further insight
it is difficult to see how to modify the techniques of \cite{JMPRN11} to
rigorously prove the persistence of the singular solutions for $\epsilon>0$
given that it is possible for \eq{eps2del} to have co-existing stable periodic orbits.
Nevertheless, Mallet-Paret and Nussbaum have announced as yet unpublished results \cite{JMPRN15}.

Given the analytical difficulties, in the current work we will instead use the algebraic results for $\epsilon=0$
from the previous sections to guide a numerical study of the
periodic solutions and bifurcation structures for \eq{eps2del} for $1\gg\epsilon>0$.
From the numerical solutions we will see that over wide parameter ranges
the singular solutions identified in the theorems
of Sections~\ref{secsingsols} and~\ref{secbifs} are indeed the limits of the solutions of the DDE \eq{eps2del}
as $\epsilon\to0$. Moreover we will find that that \eq{eps2del} has fold bifurcations of periodic
orbits at $K_1$ values which converge to
$K_1=L_{nm}$ and $K_1=M_{nm^+}$ in the limit as $\epsilon\to0$.

For $\epsilon>0$ we compute bifurcation branches numerically using DDEBiftool~\cite{ETR02}.
DDEBiftool is a suite of MATLAB \cite{matlab} routines for computing
solution branches and bifurcations of DDEs using path following and branch switching techniques.
Periodic orbits are found as the solution of a boundary value problem (BVP), using collocation techniques.
The numerical analysis details are well described in \cite{ETR02} and elsewhere, so we will not repeat
them here. We emphasise however, that periodic orbits are found by solving BVPs, and {\em not} by solving
initial value problems. This allows unstable orbits to be found just as easily as stable ones. DDEBiftool
can determine the stability of periodic orbits by computing their Floquet multipliers which also allows
us to detect bifurcations.
In \cite{DCDSA11} we already used DDEBiftool to investigate the dynamics of \eq{eps2del} in
the non-singular case $\epsilon=1$.


We will
mainly concentrate our attention on the principal branch of periodic solutions, which is the only
one on which we found large amplitude stable solutions for $\epsilon>0$. By the
principal branch, we mean the branch of periodic orbits which has the largest period among all
the Hopf bifurcations, both at the bifurcation and in the limit as $K_1\to\infty$. This usually
also corresponds to the Hopf bifurcation with the smallest value of $K_1$, but due to the vagaries of
the behaviour of the characteristic values in DDEs for $\epsilon$ very close to zero,
it is sometimes possible for a shorter period orbit to bifurcate first. If that happens
the orbit on the principal branch is initially unstable but we found numerically that it
becomes stable in a torus bifurcation while its amplitude is still very small. In the current
work we will not study small amplitude solutions or invariant tori (see \cite{CHK15} for a study of
the invariant tori of \eq{eps2del}, and \cite{KE14} for a study of \eq{eps1del} close to the singular Hopf bifurcation).
The principal branch will always correspond to the choice $n=0$ for the singular solutions and hence $m^0=0$.

Throughout this work, the amplitude of a periodic orbit of period $T$ is defined simply as the difference
between the maximum and minimum values of $u(t)$ over the period. We will take $c=1$ and $K_2=0.5$
in all our examples, and $a_1=1$, so $A=a_2/a_1=a_2$.


\begin{figure}[t]
\vspace{-5ex}
\mbox{}\hspace{-1.5cm}\scalebox{0.55}{\includegraphics{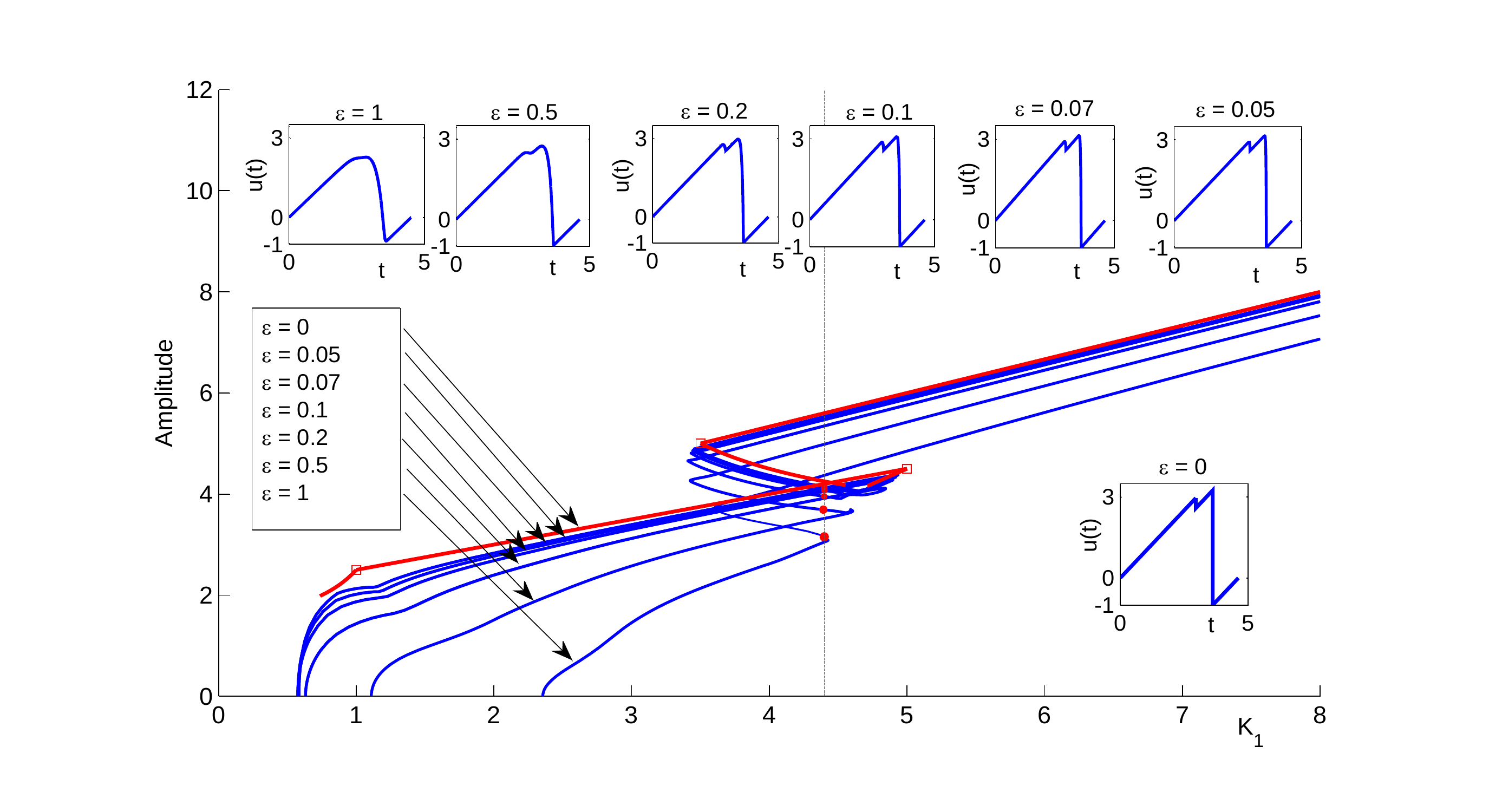}
\put(-445,200){\Large $L_{00}$}
\put(-307,180){\Large $M_{01}^+$}
\put(-670,132){\Large $L_{01}=M_{01}^-$}}
\vspace{-7ex}
\caption{Superimposed plots of the amplitude of the periodic solutions on the principal branch
created from a Hopf bifurcation as $K_1$ is varied for different values of $\epsilon$ with
$K_2 = 0.5$, $A=a_2=6$, $a_1=c=1$.
The legs of $n=0$ unimodal and type I and II bimodal singular $\epsilon=0$ solutions, are shown in red,
and the singular fold bifurcations are labelled by $L_{00}$ and $M_{01}^+$, using the notation of Section~\ref{secbifs}.
The branches for $\epsilon>0$
are computed using DDEBiftool. 
For $K_1 = 4.4$ the profile of the $\epsilon=0$ type II  singular solution on the middle leg is also shown, along
with corresponding numerically found periodic orbit profiles for $\epsilon>0$.}
\label{figA6epsto0}
\end{figure}

Fig.~\ref{figA6epsto0} illustrates the convergence of the principal solution branch as $\epsilon\to0$.
For $A=6$, the amplitude of the periodic solutions on the branch are plotted against the bifurcation parameter $K_1$
for different values of $\epsilon>0$. Also shown are the amplitudes of the $\epsilon=0$ singular
solutions following from the results of Sections~\ref{secsingsols} and~\ref{secbifs}.

For $A=6$ and $K_2=0.5$ with $n=m^0=0$, we have $m^*=1$ and $m^{**}\approx1.6$. Hence by Theorem~\ref{thmunilegs} there
are legs of unimodal singular solutions with $m=0$ for $K_1>L_{00}=3.5$ and for
with $m=1$ for $K_1\in(M_{01}^-,M_{01}^+)=(1,5)$. Theorems~\ref{thmIbif} and~\ref{thmIIbif} ensure that
there are legs of type I and type II bimodal solutions with $m=1$ for $K_1$ between $L_{00}$ and $M_{01}^+$.
By verifying the conditions of Theorems~\ref{thm2delI} and~\ref{thm2delII} we find that
the type I bimodal solutions exist for $K_1\in(4.7122,5)$ and the type II solutions for $K_1\in(3.5,4.5549)$.

Fig.~\ref{figA6epsto0} shows that
the branch converges over its entire length as $\epsilon\to0$,
and on the intervals where singular solutions exist
the amplitudes converge to those of the singular solutions. For $\epsilon>0$ the orbits have slightly smaller amplitude
than the singular solutions, which is to be expected since the singular solutions
are of sawtooth shape, while for $\epsilon>0$ the orbits are smooth, and some amplitude
is lost in the ``smoothing''. 
Insets in  Fig.~\ref{figA6epsto0}
show solution profiles on the unstable middle leg of solutions for $K_1=4.4$, $K_2=0.5$, converging to
a type II bimodal solution (also shown) as $\epsilon\to0$. Even for $\epsilon$ as large as $0.2$ the bimodal sawtooth structure
of the solution profile is very clearly seen. For larger $\epsilon$ the solution profiles are smoother, particularly
near the local maxima, but the fold structure on the solution branch persists even when $\epsilon=1$.

Fig.~\ref{figA6epsto0} is representative of the behaviour for other values of $A$. Not only do the singular solutions constructed as
in Section~\ref{secbifs} give the limiting amplitudes for the bifurcation branches as $\epsilon\to0$, but the points
$L_{nm}$ and $M_{nm}^+$ give the limiting locations of the fold bifurcations. Moreover, we will see in
Section~\ref{sec:cusp} that the singular solution theory is robust enough to show the location of codimension-two cusp-like bifurcations.

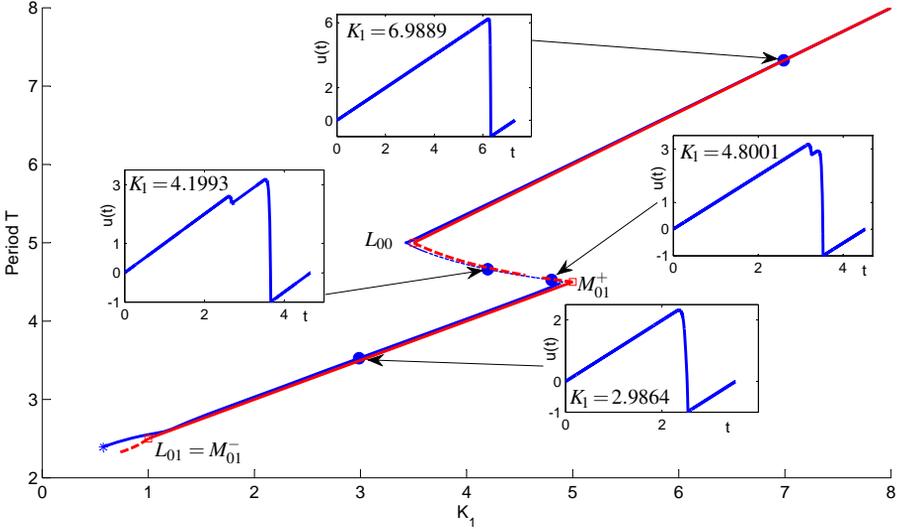
\begin{figure}[t]
\vspace{-5ex}
\mbox{}\hspace{-1.5cm}\scalebox{0.55}{\includegraphics{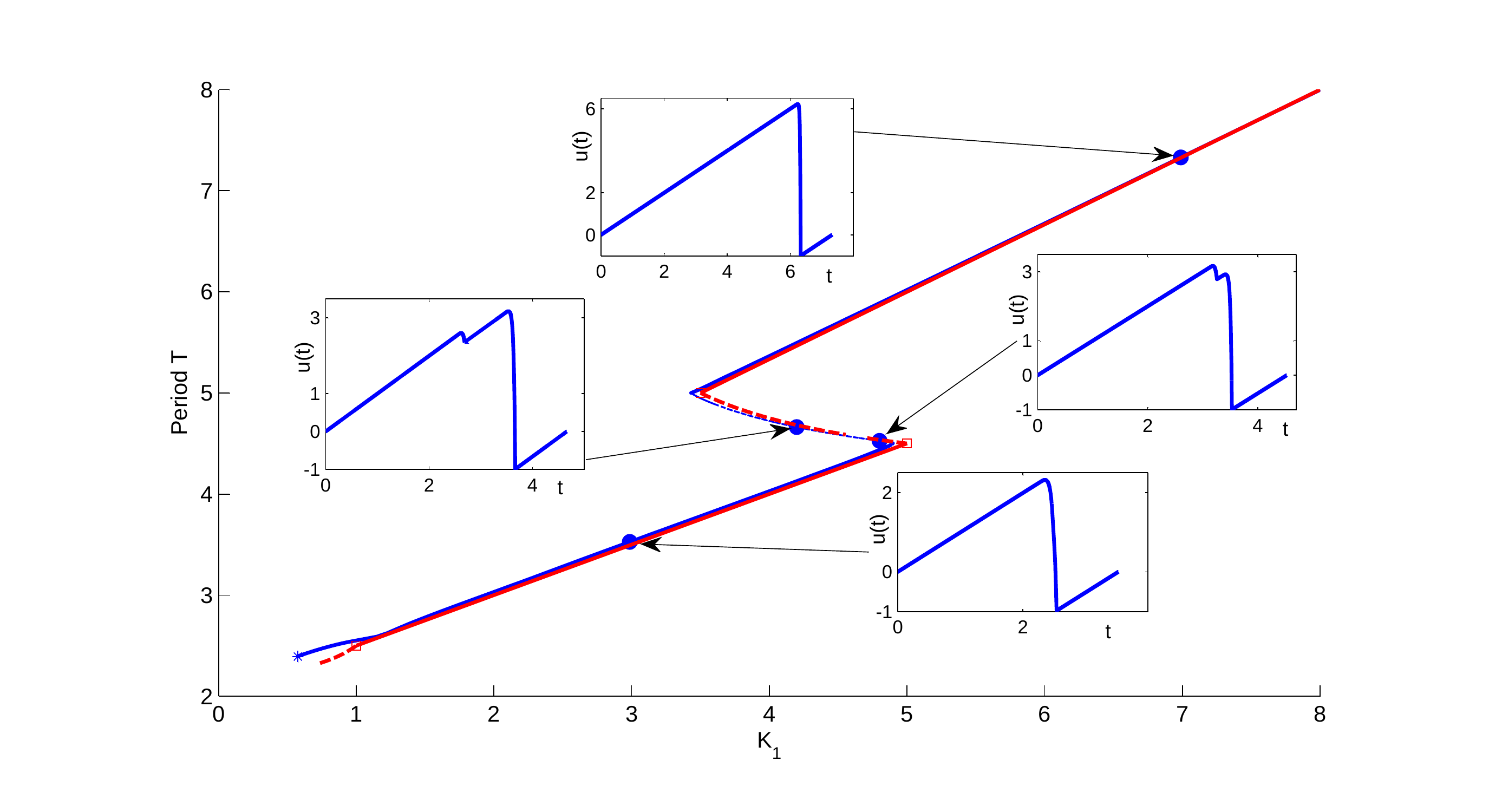}}
\put(-172,62){$K_{\!1}\!=\!2.9864$}
\put(-130,155){$K_{\!1}\!=\!4.8001$}
\put(-340,143){$K_{\!1}\!=\!4.1993$}
\put(-257,201){$K_{\!1}\!=\!6.9889$}
\put(-250,122){$L_{00}$}
\put(-169,105){$M_{01}^+$}
\put(-330,42){$L_{01}=M_{01}^-$}
\vspace{-5ex}
\caption{Periods of orbits on the principal branch
with $\epsilon=0.1$ (in blue), and the legs of $\epsilon=0$ singular unimodal and bimodal
solutions (in red). Other parameters have the same values as in Figure~\ref{figA6epsto0}.
Unimodal singular solutions
exist with $m=0$ for $K_1>L_{00}=3.5$, and with $m=1$ for $K_1\in(M_{01}^-,M_{01}^+)$ where
$M_{01}^-=L_{01}=1$ and $M_{01}^+=5$. Legs of bimodal solutions of type I exist for
$K_1\in(4.7122,5)$ and of type II for $K_1\in(3.5,4.5549)$.
Insets show solution profiles for $\epsilon=0.1$. For $\epsilon=0.1$ solutions are unstable
between the folds, and stable on the upper and lower legs of the branch.}
\label{figPereps01}
\end{figure}

Fig.~\ref{figPereps01} shows the bifurcation branch for $\epsilon=0.1$ plotted against the period, along
with the $\epsilon=0$ singular solutions for the same values of the other parameters as Fig.~\ref{figA6epsto0}.
Profiles of periodic orbits for $\epsilon=0.1$ are also shown. Where the singular solutions exist their period
is very close to that of the numerically computed $\epsilon=0.1$ solutions. The agreement in period is even
better than the agreement in amplitude between the singular and $\epsilon=0.1$ solutions seen in Fig.~\ref{figA6epsto0}.
For $\epsilon=0.1$ the periodic solutions before the first fold and after the second fold are stable with an
interval of bistability of periodic solutions between the two folds.
Periodic solutions on the leg of the branch between the two folds are always unstable, and are bimodal
for at least part of the leg.
On the leg of unstable solutions in Fig.~\ref{figPereps01}
we see two types of bimodal periodic solutions for $\epsilon=0.1$, where
either the first or second local maximum after the solution minimum is higher (see insets with $K_1=4.8001$ and $K_1=4.1993$).
Note the resemblance between the profiles of the $\epsilon=0.1$ solutions and
the singular solutions shown in Figs.~\ref{figtypeI} and~\ref{figtypeII}, which is not coincidental; the construction
of the singular solutions in Section~\ref{secsingsols} was guided by preliminary numerical computations.

\begin{figure}[t]
\vspace{-5ex}
\mbox{}\hspace{-1.5cm}\scalebox{0.55}{\includegraphics{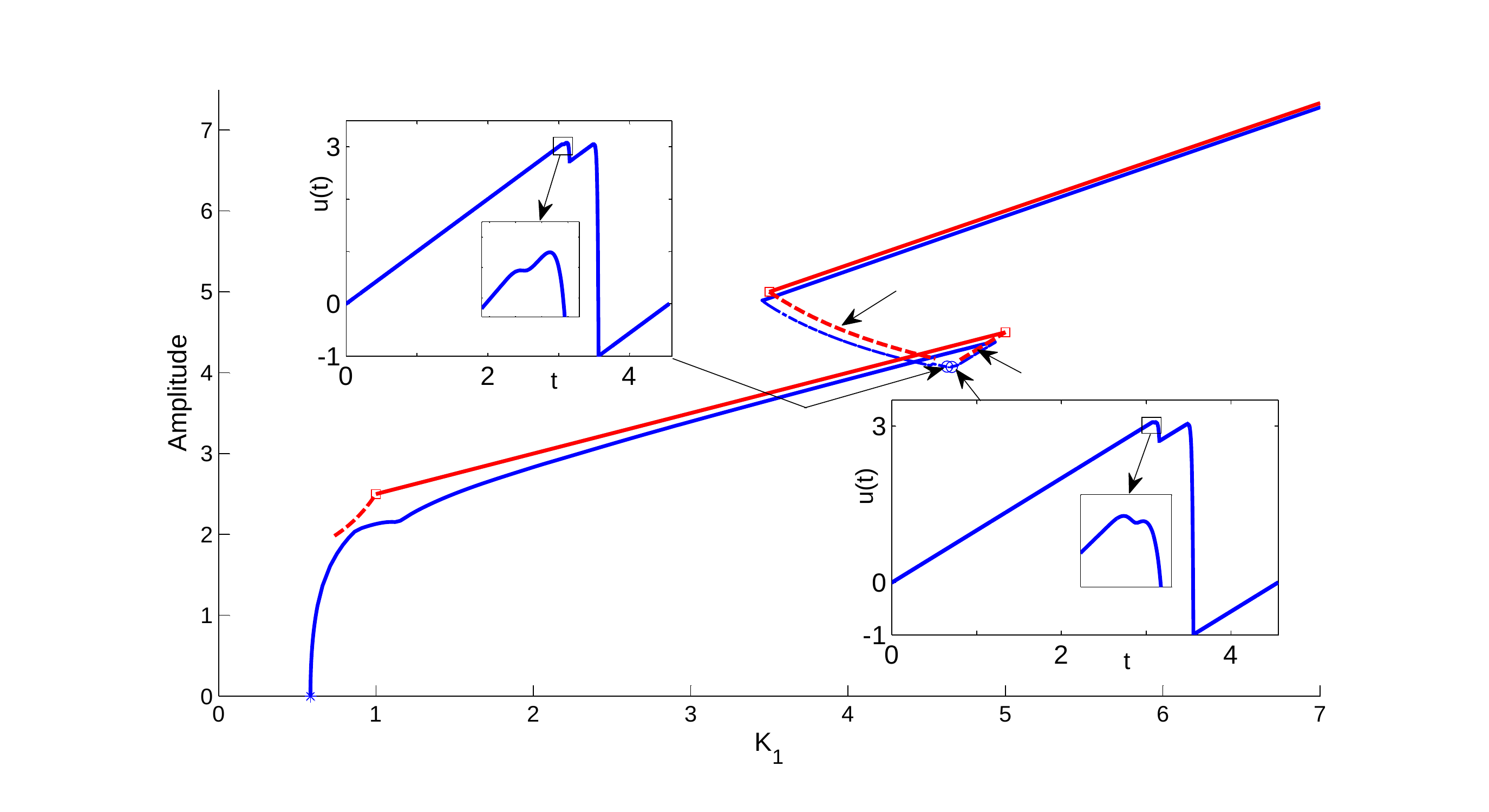}
\put(-600,355){\Large $K_{\!1}\!=\!4.6294$}
\put(-310,205){\Large $K_{\!1}\!=\!4.6590$}
\put(-313,275){\Large $\epsilon=0$, type II bimodal, $m=1$}
\put(-247,229){\Large $\epsilon\!=\!0$, type I bimodal, $m\!=\!1$}
\put(-320,310){\rotatebox{17}{\Large $\epsilon=0$, unimodal, $m=0$}}
\put(-280,298){\rotatebox{17}{\Large $\epsilon=0.05$}}
\put(-575,180){\rotatebox{14}{\Large $\epsilon=0$, unimodal, $m=1$}}

\put(-625,100){\Large $\epsilon=0.05$}

\put(-410,280){\Large $L_{00}$}
\put(-253,252){\Large $M_{01}^+$}
\put(-660,170){\Large $L_{01}=M_{01}^-$}
}
\vspace{-5ex}
\caption{Amplitude of $\epsilon=0.05$ and singular solutions, with the other parameters
taking the same values as in Fig.~\ref{figA6epsto0}. For $\epsilon=0.05$ periodic
solutions with three local maxima per period (trimodal solutions) are found on a small interval
for $K_1\in(4.6294,4.6668)$. Profiles of two of these trimodal periodic solutions
with $K_1=4.6294$ and $K_1=4.6590$ are shown as insets.} \label{figAmpeps005}
\end{figure}
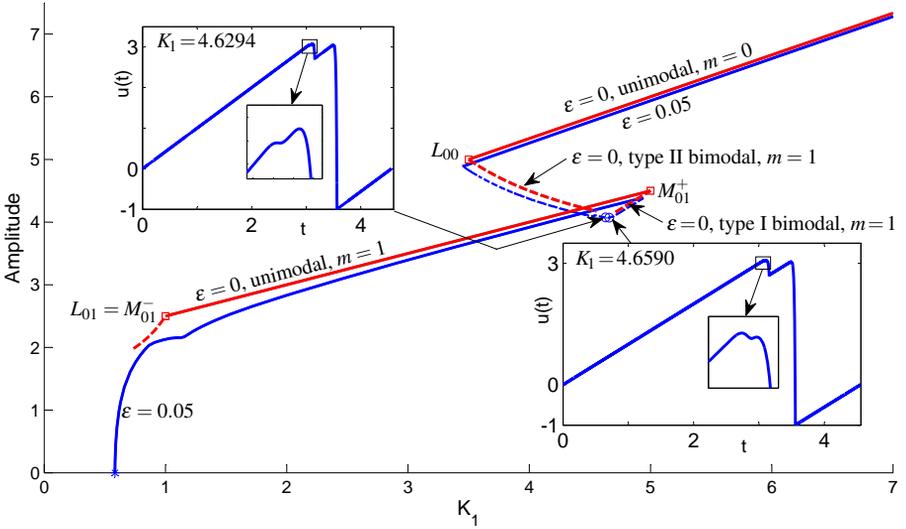


We see in Figs.~\ref{figA6epsto0} and~\ref{figPereps01}
that the $\epsilon>0$ branches pass continuously through
the gap in the singular solution branch.
With $\epsilon=0.1$ there is a smooth transition between the two types of bimodal solutions
along the unstable leg, whereas Fig.~\ref{figbranchex} suggests that there should be a gap
between the intervals where these two types of solutions exist in the limit as $\epsilon\to0$.
So we next investigate periodic solution profiles as
$\epsilon\to0$ paying particular attention to the legs between the fold bifurcations where those gaps occur.

For $\epsilon=0.05$ we find periodic solutions with three local maxima per period (which we call trimodal solutions)
on the unstable leg of the branch for an interval of $K_1$ values
which falls within the gap between the
type I and type II bimodal singular solutions.
Fig.~\ref{figAmpeps005} shows the profiles of two of these trimodal periodic solutions.
We note that both profiles are similar to bimodal solutions,
but that in both cases the first local maxima of the bimodal solution has split into
two local maxima. For parameter values close to where the type I bimodal solutions
exist (including $K_1=4.659$) the first two local maxima of the solution resemble those of
a type I bimodal solution (with the first local maxima after the global minima being the global maxima),
while for parameter values close to the type II bimodal solutions
(including $K_1=4.6294$) the first two local maxima of the solution resemble those of
a type II bimodal solution (with the second local maxima after the global minima being the global maxima).


\begin{figure}
\vspace{-2ex}
\mbox{}\hspace{-0.5cm}\scalebox{0.36}{\includegraphics{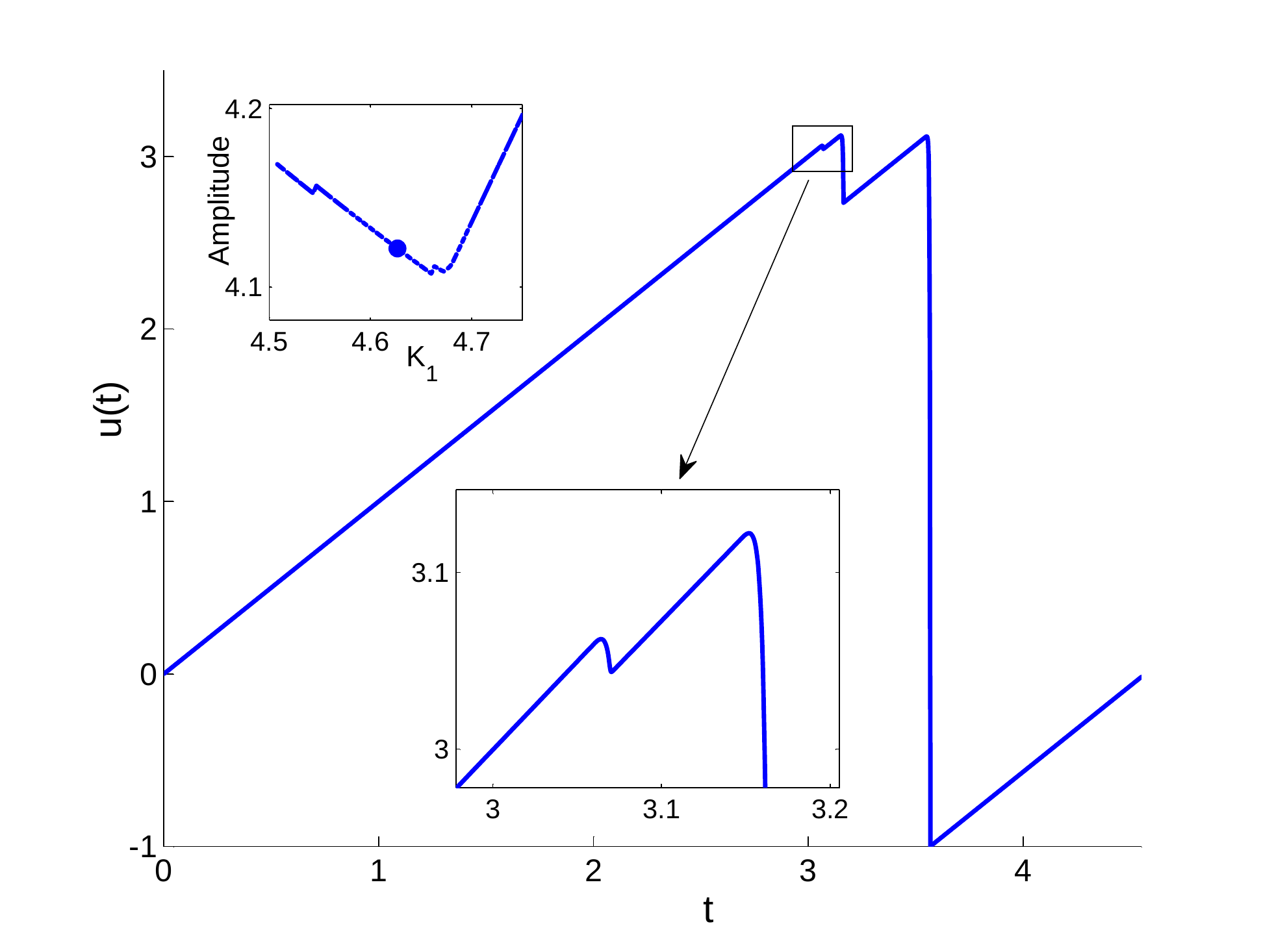}\hspace{-2.5cm}\includegraphics{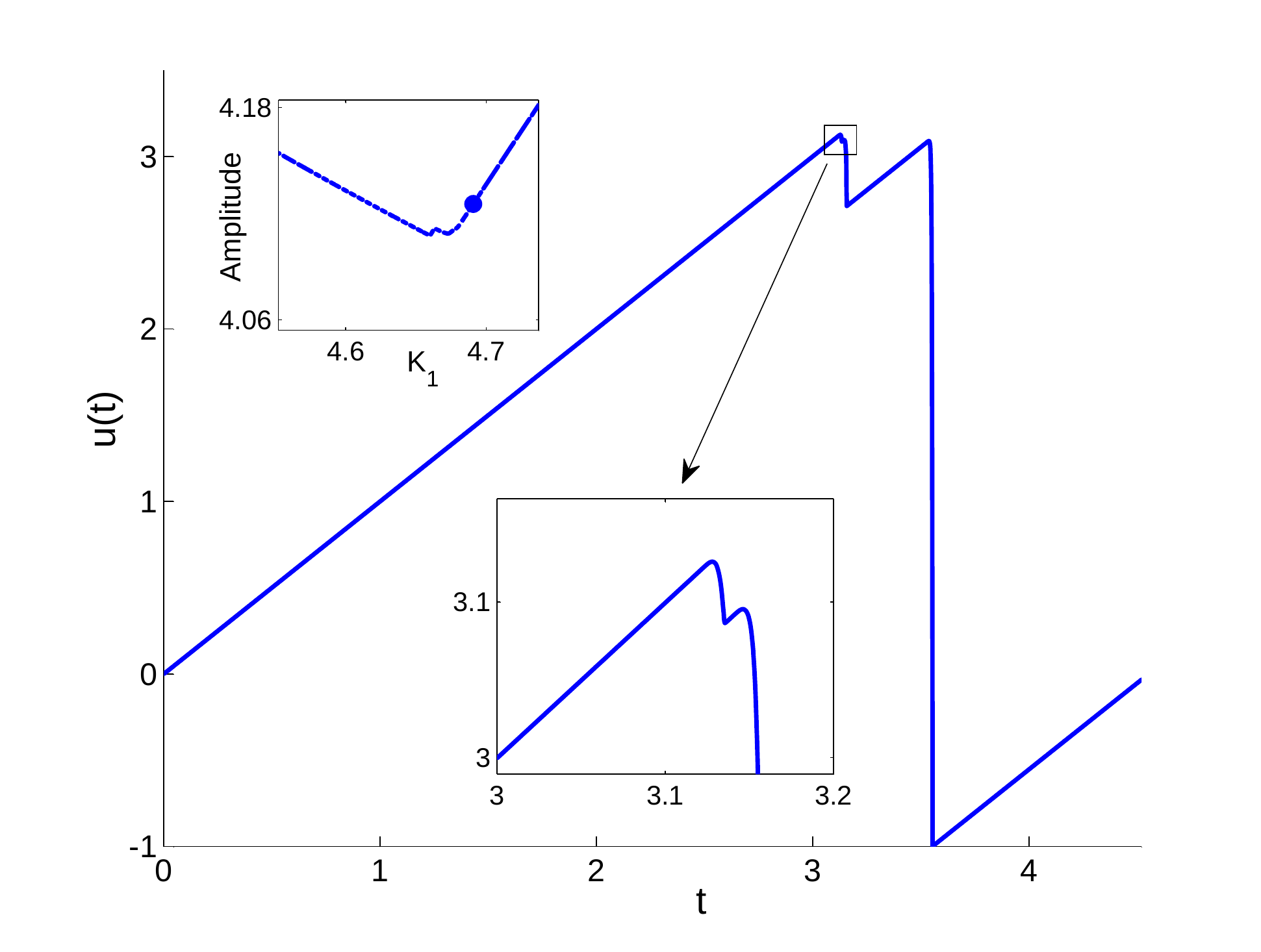}}
\put(-280,135){(i)}
\put(-105,135){(ii)}

\vspace{-1.8ex}

\begin{center}
\includegraphics[scale=0.45]{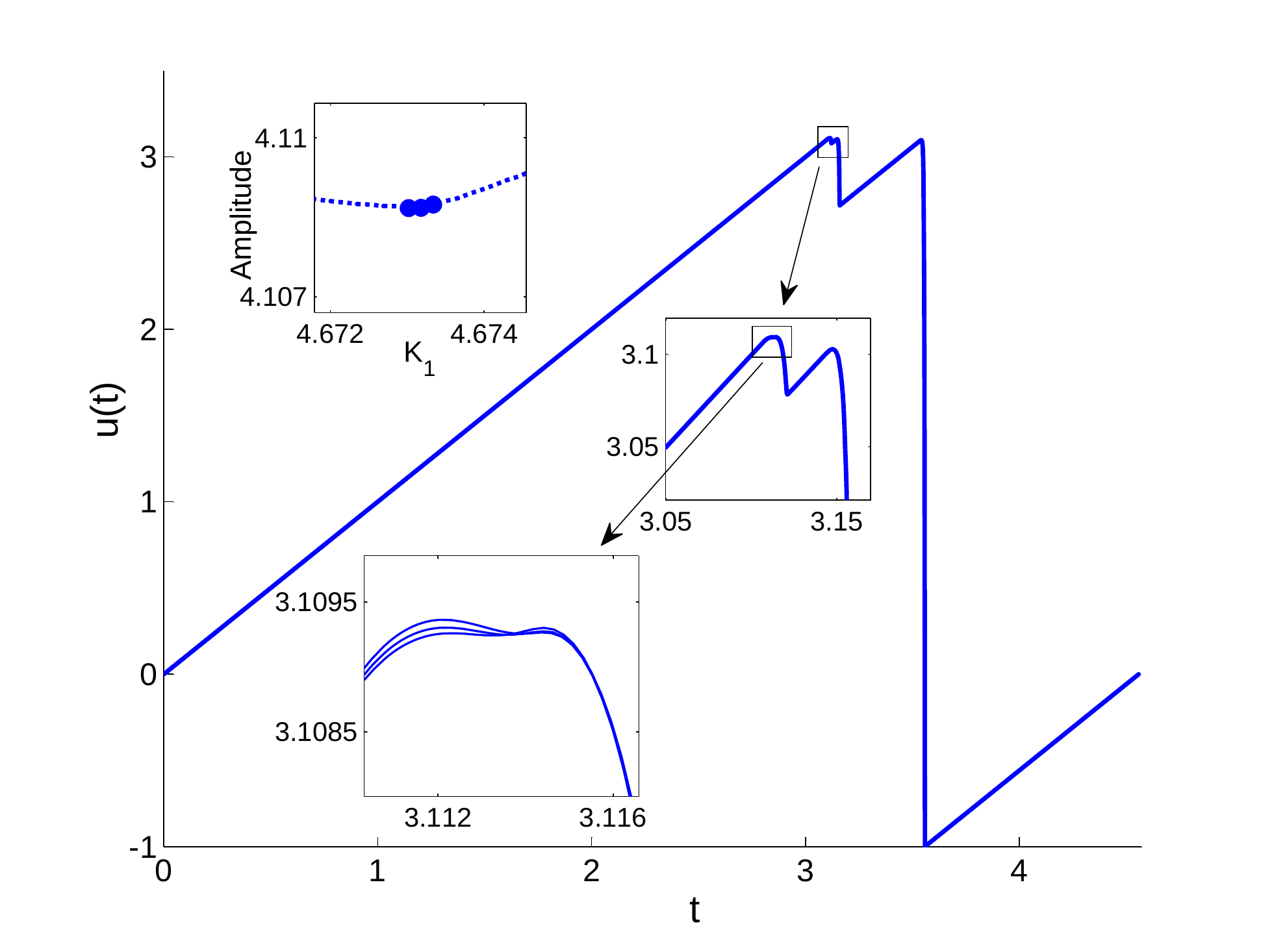}
\put(-130,167){(iii)}

\vspace{-5ex}

\mbox{}
\end{center}
\caption{For $\epsilon=0.01$ and the other parameters as in the preceding figures, trimodal periodic solutions are found for
$K_1\in(4.5747,4.6999)$ and within this interval solutions are quadrimodal for $K_1\in(4.6730,4.6733)$.
Shown are two example profiles of trimodal solutions (i) for $K_1=4.6908$, and (ii) for $K_1=4.6266$,
and (iii) three quadrimodal solutions (plotted on the same axis).}
\label{figprofA6eps001}
\end{figure}



With $\epsilon=0.01$ on the unstable leg of the branch,
Type I-like bimodal solutions occur in the
approximate range $K_1\in(4.6998,4.9802)$. At $K_1\approx4.6998$ there is a transition to a trimodal solution,
and trimodal solutions exist for $K_1\in(4.6735,4.6998)$.
The numerically found trimodal solution for $K_1=4.6908$ is shown in Fig.~\ref{figprofA6eps001}(i).
Again we see (in the inset) that it is the first maximum of the solution which splits into two to
form the trimodal solution.
Around $K_1\approx4.673$ there is a brief interval of quadrimodal solutions, where the first maximum of
the trimodal solution splits into two as shown in Fig.~\ref{figprofA6eps001}(iii).
There is then another interval of trimodal solutions for $K_1\in(4.5746,4.673)$, with the solution for
$K_1=4.6266$ shown in Fig.~\ref{figprofA6eps001}(ii). Finally for $K_1<4.5746$ the
solutions are bimodal (and type II-like). Comparing the trimodal solutions in Figs.~\ref{figAmpeps005}
and~\ref{figprofA6eps001} we see that the trimodality is much more clearly defined
for the smaller value of $\epsilon$ with the profiles in Fig.~\ref{figprofA6eps001} much more
`sawtooth-like' than the smoother profiles seen in Fig.~\ref{figAmpeps005}. Moreover, 
for both $\epsilon=0.05$ and $\epsilon=0.01$ the trimodal solution in the interval adjacent to the type I bimodal solutions
has a larger first peak than second peak, just as the type I bimodal solutions do, and similarly for type II bimodal solutions
and the trimodal solutions in the adjacent parameter interval the second peak is larger. This could lead us to
define type I and type II trimodal solutions which could 
be found algebraically in the $\epsilon=0$ limit
using our singular solution techniques. 
However each would involve
about 15 intervals of parametrisation, which would be tedious beyond belief.
Moreover, Fig.~\ref{figprofA6eps001} suggests that the quadrimodal solutions
also come in both types, and we suspect that as $\epsilon\to0$ there is a cascade
of solutions with arbitrarily many maxima, 
and some form of self-similarity to
the evolution of the periodic solution profile in the limit as $\epsilon\to0$.
However all these multimodal solutions lie on the unstable leg of the branch,
and we will not pursue their construction here.

Although we omit the algebraic construction of trimodal solutions, our construction of the bimodal solutions
is sufficient to compare the bimodal to trimodal solution transition with the unimodal to bimodal transitions/bifurcations
already studied algebraically in Section~\ref{secbifs}. The interval of $K_1$ values for
which type I bimodal solutions exist was found by checking the conditions of Theorem~\ref{thm2delI}. In all the
examples shown above, and indeed in the other examples of type I bimodal solutions that occur later in this
paper we only find two different behaviours which arise at the ends of these intervals. One case is when
$K_1=M_{nm}^+$ indicating a transition or bifurcation between
a unimodal and type I bimodal solution as studied in Theorem~\ref{thmIbif} and
the comments after that theorem. At the other end of the leg of type I bimodal solutions where the $\epsilon>0$ numerics
indicate a transition to a trimodal solution, algebraically in all the examples shown here we find that
the lower bound on $\theta T_1$ in \eq{tIthbds} fails. From the proof of Theorem~\ref{thmIbif} we see that
equality in this bound corresponds to $F(10j+3)=0$. This is similar to the transition from
a unimodal to type I bimodal solution at $K_1=M_{nm}^+$ as described after Theorem~\ref{thmIbif}.
Then we saw that the failure of the condition $F(5j+3)<0$ in the unimodal singular solution led to the creation
of a second subinterval of $J^*$ in the periodic orbit. In an analogous manner the failure of the condition $F(10j+3)<0$
in the type I bimodal singular solution can lead to a solution where $J^*$ consists of three disjoint intervals
per period and the resulting solution is trimodal.

The transition from type II bimodal to trimodal solutions also appears to be similar to the transition from unimodal to
type II bimodal solutions. After Theorem~\ref{thmIIbif} we noted that at the transition between unimodal
and type II bimodal solutions at $K_1=L_{nm-1}$ we have $\theta\to1$ for the unimodal solution and $\theta\to0$ for the
type II bimodal solution as $K_1\to L_{nm-1}$.
Checking the conditions of Theorem~\ref{thm2delII} we find that in all
the examples above $\theta\to0$ as $K_1\to L_{nm-1}$ and $\theta\to1$ at the other end of each leg
of type II bimodal solutions. We would expect the solution to transition to a type II trimodal solution with
$\theta=0$ at this point.



\section{Cusp-like Bifurcations}
\label{sec:cusp}

Here we investigate the cusp-like bifurcations, identified by Theorem~\ref{thmIIbif}, where
fold bifurcations of singular solutions disappear
when $m=m^{*\!}(n)+\frac{1+K_2}{2+K_2}$ and $K_1=L_{nm-1}$,
or equivalently when
\be \label{cusp}
A=1+\frac{1+m(2+K_2)}{1+n(2+K_2)},\qquad K_1=L_{nm-1}=1+\frac{1+n(2+K_2)}{m-n}.
\ee
On the principal branch $n=0$, so the cusps occur when $A=2+m(2+K_2)$ and $L_{0m-1}=1+1/m$.
Taking $K_2=0.5$, the cusp-like bifurcations for the singular solutions occur
when $A=4.5,7,9.5,...$ and $L_{0m-1}=2,3/2,4/3,\ldots$.
Here we will investigate the first two such bifurcations both in the
singular case with $\epsilon=0$, and numerically for $0<\epsilon\ll1$.



\begin{figure}
\vspace{-5ex}
\mbox{}\hspace{-1.5cm}\scalebox{0.55}{\includegraphics{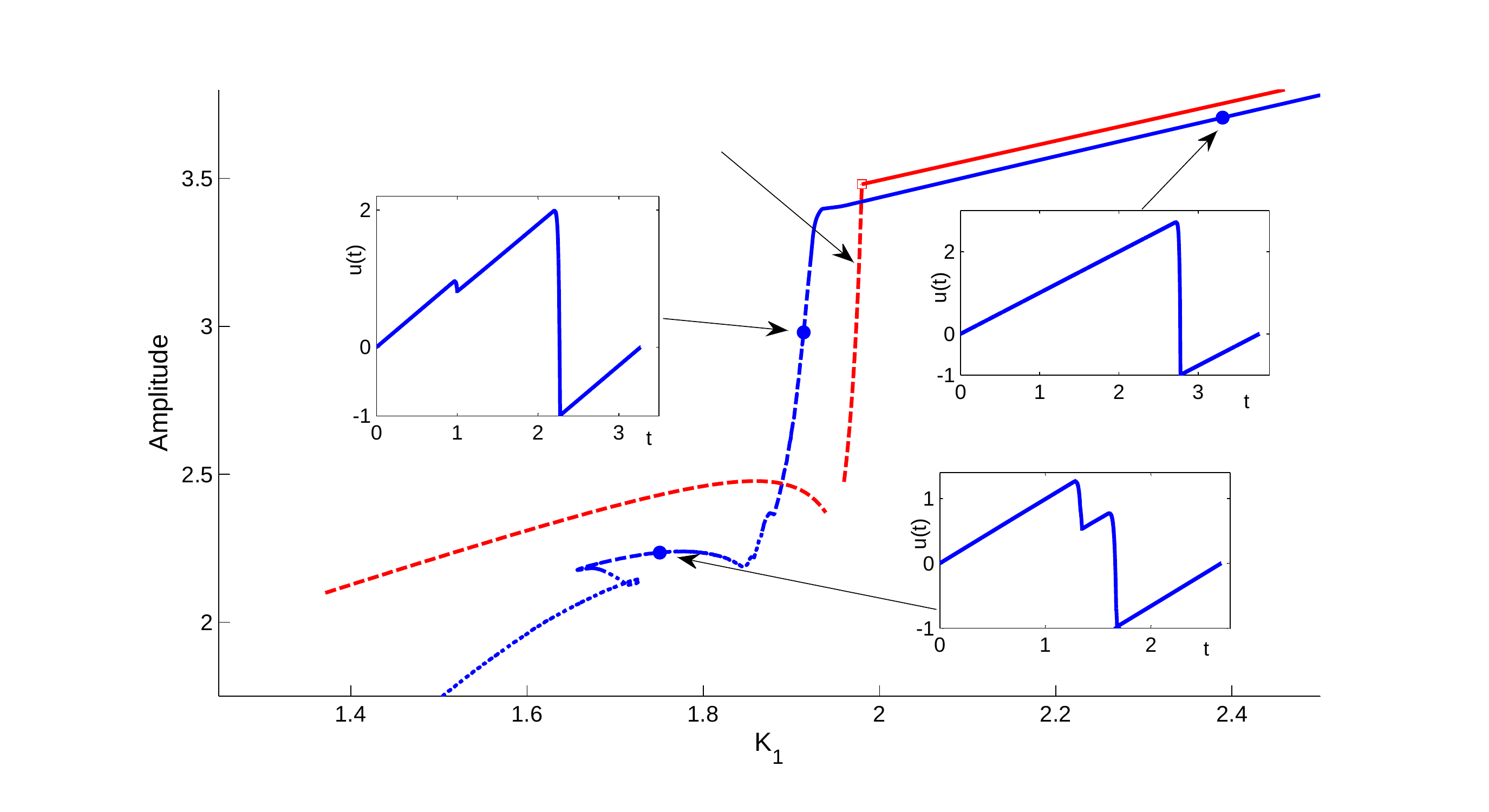}
\put(-260,240){\Large $K_1\!=\!2.3894$}
\put(-590,218){\Large $K_1\!=\!1.9140$}
\put(-485,355){\Large $\epsilon\!=\!0$, type II bimodal, $m\!=\!1$}
\put(-305,350){\rotatebox{12}{\Large $\epsilon\!=\!0$, unimodal, $m\!=\!0$}}
\put(-285,105){\Large $K_1\!=\!1.7508$}
\put(-625,125){\rotatebox{16}{\Large $\epsilon\!=\!0$, type I bimodal, $m\!=\!1$}}
}
\vspace{-5ex}
\caption{Amplitude of Type I and II bimodal singular solutions with $m=1$ and unimodal solutions
with $m=0$ for $A=4.48$, $K_2=0.5$ and $n=0$. Also shown is the numerically computed
principal branch of periodic orbits for the same parameters except $\epsilon=0.02$. Unimodal,
bimodal and multimodal periodic orbits for $\epsilon=0.02$ occur on the solid, dashed or dotted parts of
the curve, respectively. Insets show profiles of stable numerically computed orbits with $\epsilon=0.02$
corresponding to unimodal, and type I and II bimodal solutions.}
\label{figA448amp}
\end{figure}

\begin{figure}
\vspace{-5ex}
\mbox{}\hspace{-1.5cm}\scalebox{0.55}{\includegraphics{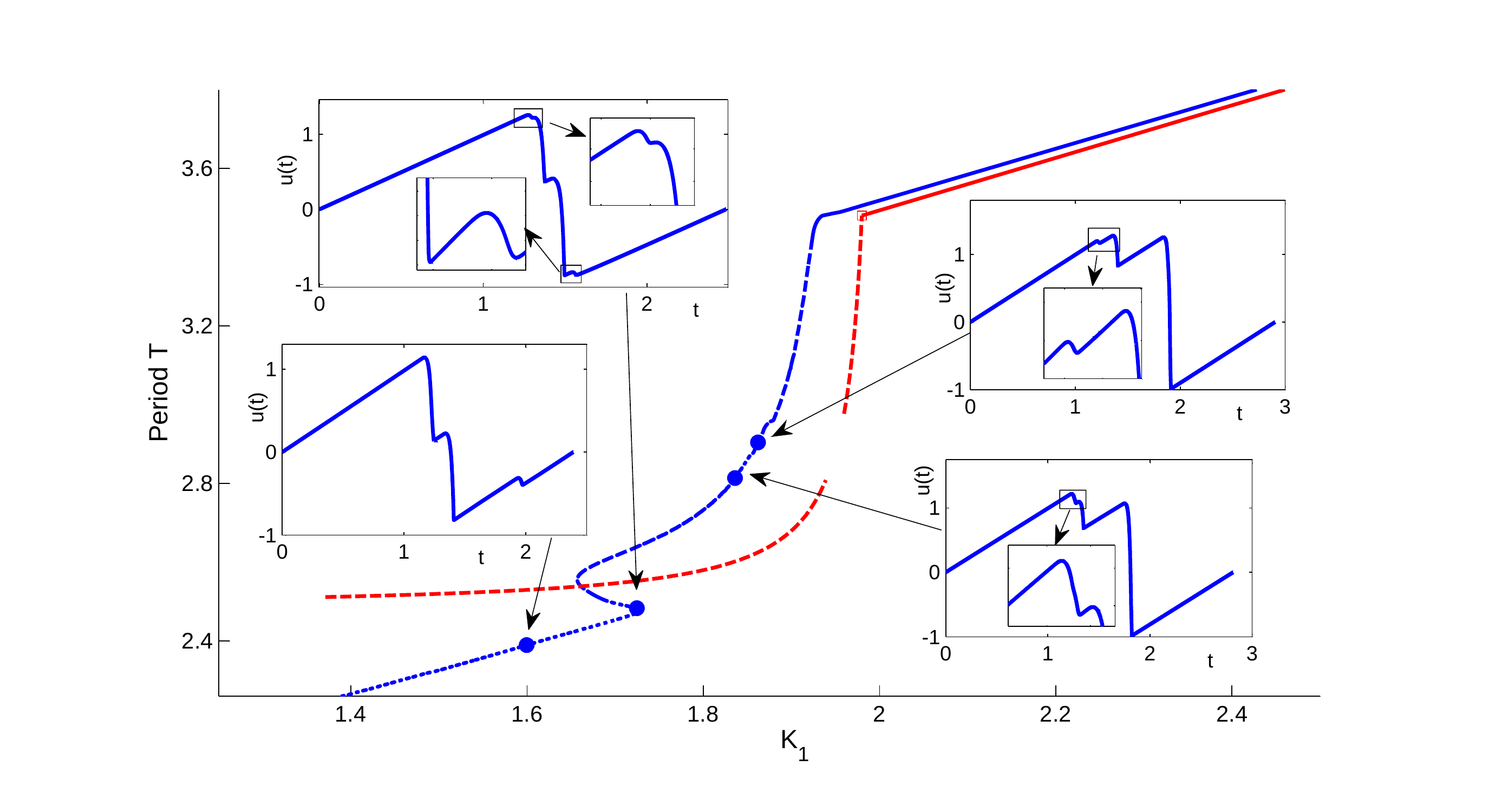}
\put(-636,155){\Large $K_1\!=\!1.5997$}
\put(-618,363){\Large $K_1\!=\!1.7248$}
\put(-205,174){\Large $K_{\!1}\!=\!1.8361$}
\put(-190,313){\Large $K_{\!1}\!=\!1.8622$}
\put(-365,335){\Large $L_{00}\!=\!1.98$}
}
\vspace{-5ex}
\caption{Periods of $\epsilon=0$ singular solutions and the numerically computed $\epsilon=0.02$
solution branch for the same parameters as Fig.~\protect\ref{figA448amp}. Insets show profiles of
numerically computed type I and II trimodal solutions for $K_1=1.8361$ and $1.8622$, as well as
examples of trimodal and quadrimodal solutions found earlier on the branch. All the inset solutions are
unstable with a complex conjugate pair of unstable Floquet multipliers.}
\label{figA448per}
\end{figure}

\begin{figure}
\vspace{-5ex}
\mbox{}\hspace{-1.5cm}\scalebox{0.55}{\includegraphics{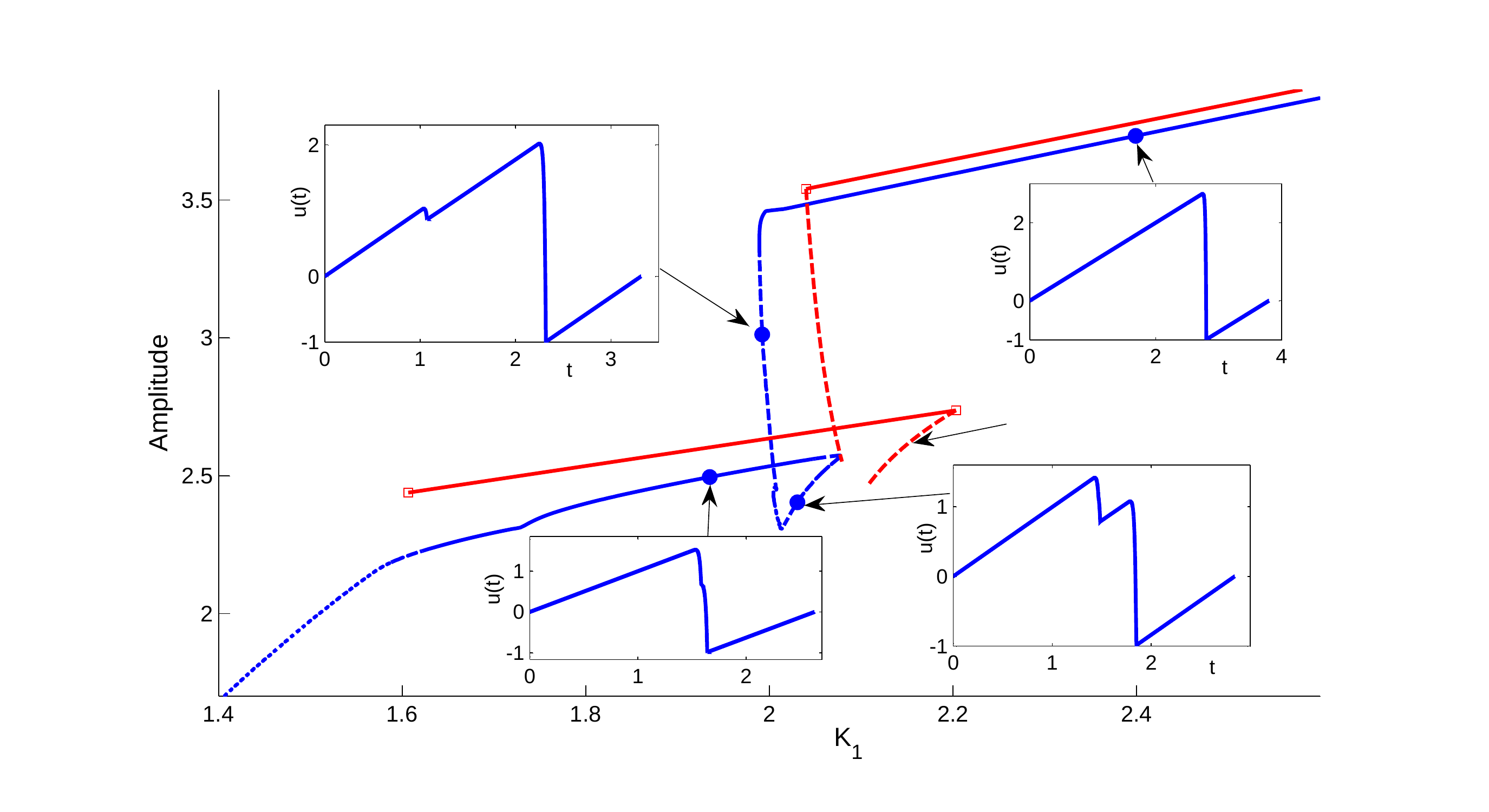}
\put(-235,260){\Large $K_1\!=\!2.3991$}
\put(-507,88){\Large $K_1\!=\!1.9351$}
\put(-605,260){\Large $K_1\!=\!1.9923$}
\put(-275,97){\Large $K_1\!=\!2.0306$}
\put(-350,348){\rotatebox{10}{\Large $\epsilon\!=\!0$, unimodal, $m\!=\!0$}}
\put(-550,185){\rotatebox{7}{\Large $\epsilon\!=\!0$, unimodal, $m\!=\!1$}}
\put(-250,205){\Large $\epsilon\!=\!0$, type I bimodal, $m\!=\!1$}
\put(-352,276){\Large $\epsilon\!=\!0$, type II}
\put(-352,262){\Large bimodal, $m\!=\!1$}
}
\vspace{-5ex}
\caption{Amplitudes of Type I and II bimodal singular solutions with $m=1$ and unimodal solutions
with $m=0$ and $m=1$ for $A=4.54$, $K_2=0.5$ and $n=0$. The numerically computed
principal branch of periodic orbits with $\epsilon=0.02$ is also shown.
Insets show profiles of stable numerically computed periodic orbits.
The unimodal orbits with $K_1=1.9351$ and $K_1=2.3991$ are stable, while the other two
orbits which correspond to type I and II bimodal solutions are unstable
with one real unstable Floquet multiplier.}
\label{figA454amp}
\end{figure}

\begin{figure}
\vspace{-5ex}
\mbox{}\hspace{-1.5cm}\scalebox{0.55}{\includegraphics{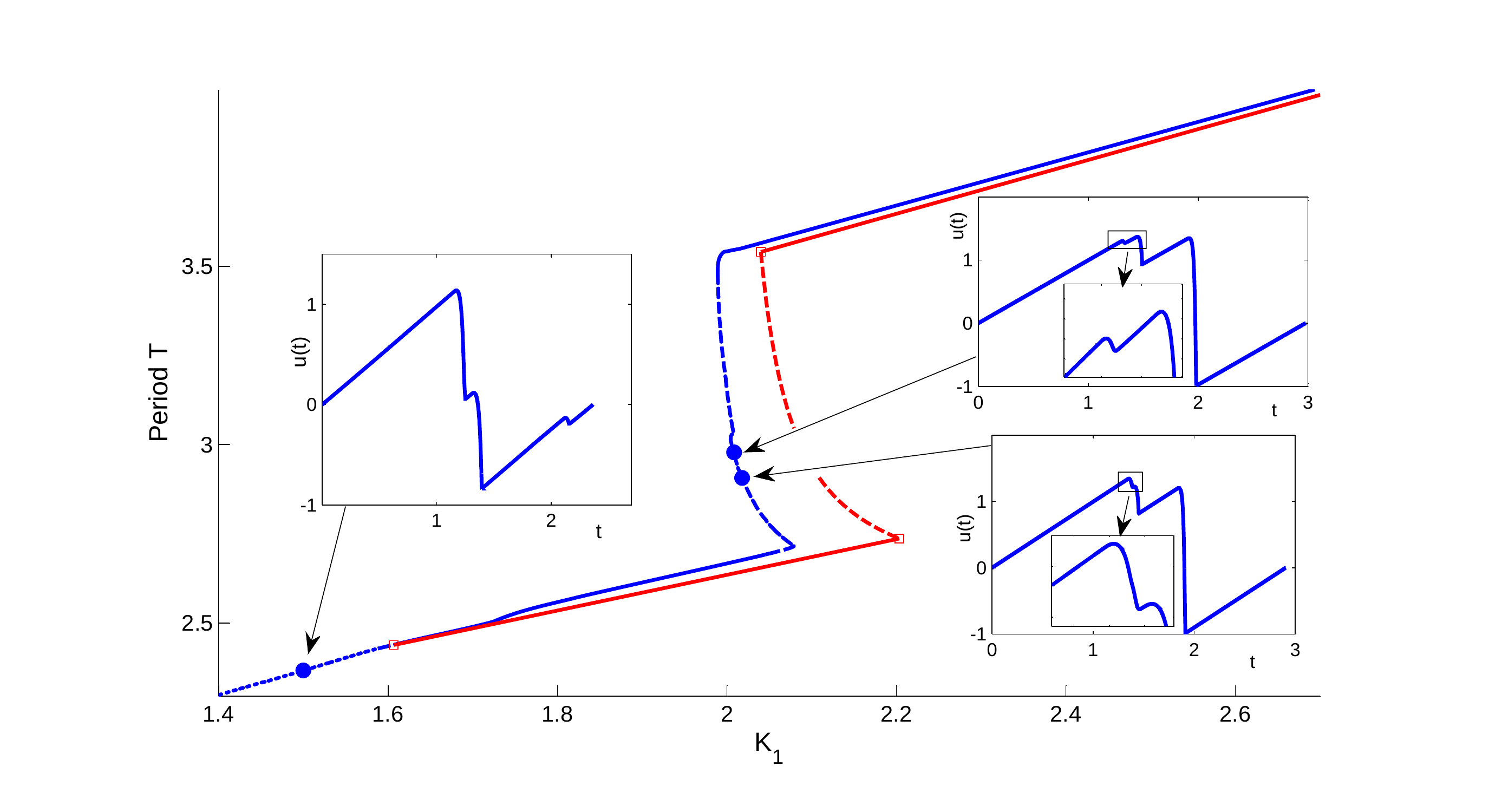}
\put(-260,188){\Large $K_1=2.0179$}
\put(-270,315){\Large $K_1=2.0084$}
\put(-535,282){\Large $K_1=1.5001$}

\put(-380,285){\Large $L_{00}\!=\!2.04$}
\put(-330,125){\Large $M_{01}^+$}
\put(-355,105){\Large $=2.2034$}
\put(-575,75){\Large $M_{01}^-\!=\!1.6066$}
}
\vspace{-5ex}
\caption{Periods of $\epsilon=0$ singular solutions and the numerically computed
$\epsilon=0.02$  solution branch for the same parameters as Fig.~\protect\ref{figA454amp}. Insets show profiles of
trimodal solutions on the $\epsilon=0.02$ branch, which are all found to be unstable.}
\label{figA454per}
\end{figure}

Figs.~\ref{figA448amp}-\ref{figA454per} illustrate the change in the dynamics near to $K_1=L_{00}$
as $A$ passes through $4.5$.  Amplitude and period plots of the unimodal and bimodal singular solutions
are shown for $A=4.48$ and $A=4.54$
along with the numerically computed principal branch of periodic solutions
with $\epsilon=0.02$. 
By Theorem~\ref{thmIIbif} we have a leg of unimodal singular solutions
for $K_1>L_{00}=A-2-K_2$ with period $T\to A-1$ as $K_1\to L_{00}$, and since $c=1$, amplitude equal to the period.
For $A=4.48$ and $K_2=0.5$ this gives $L_{00}=1.98$ with the period and amplitude of the unimodal solutions tending
to $3.48$ as $K_1\to L_{00}$. Moreover \eqref{m*} implies $m^*(0)+\frac{1+K_1}{1+K_2}<1$ and
Theorem~\ref{thmIIbif}(ii) gives the existence of a leg of type II bimodal solutions for $K_1<L_{00}$.
We see from Figs.~\ref{figA448amp} and~\ref{figA448per} that the solution branch for $\epsilon=0.02$
behaves similarly though the transition point is perturbed to $K_1\approx1.9221$ which is  slightly less than $L_{00}=1.98$.

For $A=4.54$, $m^0(0)=0$ and $m^*(0)+\frac{1+K_1}{1+K_2}>1$, thus
Theorem~\ref{thmIIbif}(i) gives the existence of a leg of type II bimodal solutions for $K_1>L_{00}$,
resulting in a fold bifurcation.
Figs.~\ref{figA454amp} and~\ref{figA454per} show that the $\epsilon=0.02$ solution branch
also has a fold bifurcation for $K_1\approx1.9892$ slightly less than $L_{00}=2.04$,
with the solution profile changing from unimodal to bimodal (see insets in Fig.~\ref{figA454amp} for
$K_1=1.9351$ and $K_1=2.3991$).
Similarly, for $A=4.54$, Theorem~\ref{thmIbif} indicates a second fold bifurcation of singular solutions at
$K_1=M_{01}^+=2.2034$ where the solution profile also transitions between a unimodal and a
type I bimodal solution, and Figs.~\ref{figA454amp} and~\ref{figA454per} show that the
numerically computed branch for $\epsilon=0.02$ has a similar bifurcation at $K_1\approx2.079$. Insets in Fig.~\ref{figA454amp} for
$K_1=1.9923$ and $K_1=2.0306$ show the resulting unimodal and bimodal solutions each side of this fold.
Trimodal solutions of both types are also observed on the $\epsilon=0.02$ branch for both $A=4.48$ and $A=4.54$ for
parameters in the gap between the $m=1$ bimodal type I and type II solutions, and these are illustrated in insets
in Figs.~\ref{figA448per} and~\ref{figA454per}.

One important aspect of this cusp-like bifurcation is that it has the potential to create stable bimodal and multimodal
periodic solutions. In Section~\ref{secnumerics} 
all of the solutions
with more than one local maxima per period were unstable occurring on the leg of the
bifurcation branch between the fold bifurcations. The bimodal and trimodal solutions occurring between
the fold bifurcations for $A=4.54$ illustrated in Figs.~\ref{figA454amp} and~\ref{figA454per}
are also unstable. However before the cusp bifurcation with
$A\leappr4.5$
stable periodic solutions with more than one local maxima per period
occur close to $K_1=L_{00}$ on the principal branch. Fig.~\ref{figA448amp} illustrates stable
bimodal solutions for $\epsilon=0.02$ and $K_1=1.7508$ and $K_1=1.9140$ which correspond to $\epsilon=0$
type I and type II bimodal singular solutions.
Interestingly the type I and II trimodal solutions for
$A=4.48$ and $K_1=1.8361$ and $K_1=1.8622$ illustrated in Fig.~\ref{figA448per} are unstable,
even though they are not between fold bifurcations. Both these periodic orbits have a pair of
complex conjugate unstable Floquet multipliers, indicating a possible torus bifurcation.

The agreement between the singular solution legs and the $\epsilon=0.02$ branch is not as good for the smaller values of $K_1$
shown in Figs.~\ref{figA448amp} and~\ref{figA448per} when $A=4.48$. In particular for the singular solution there is a
leg of type I bimodal solutions with $m=1$ for $K_1\in(M_{01}^-,M_{01}^+)\approx(1.3711,1.9391)$, but
for $\epsilon=0.02$ the corresponding bimodal solution only exists in the interval $K_1\in(1.6571,1.8259)$.
To explain
this note that the interval $(M_{01}^-,M_{01}^+)$ is derived from the roots of a quadratic with parameters such
that it is close to its double root and is thus very sensitive to the value of $A$; decreasing $A$ to $4.41$ causes
this interval and the associated type I bimodal singular solutions to vanish.
Computations with other values of $A$ (not shown)
suggest that for $\epsilon=0.02$ the fold bifurcation associated with the point $K_1=L_{00}$ disappears
at about $A=4.52$, whereas this occurs at $A=4.5$ for the singular solution, hence the $\epsilon=0.02$ solution
branch for $A=4.48$ is actually twice as far from its critical value as the singular solutions shown in the same figures, so
it is not surprising that $\epsilon=0.02$ bimodal solution exists on a smaller interval.


\begin{figure}[t!]
\vspace{-5ex}
\begin{center}
\scalebox{0.55}{\mbox{}\hspace{-3cm}\includegraphics{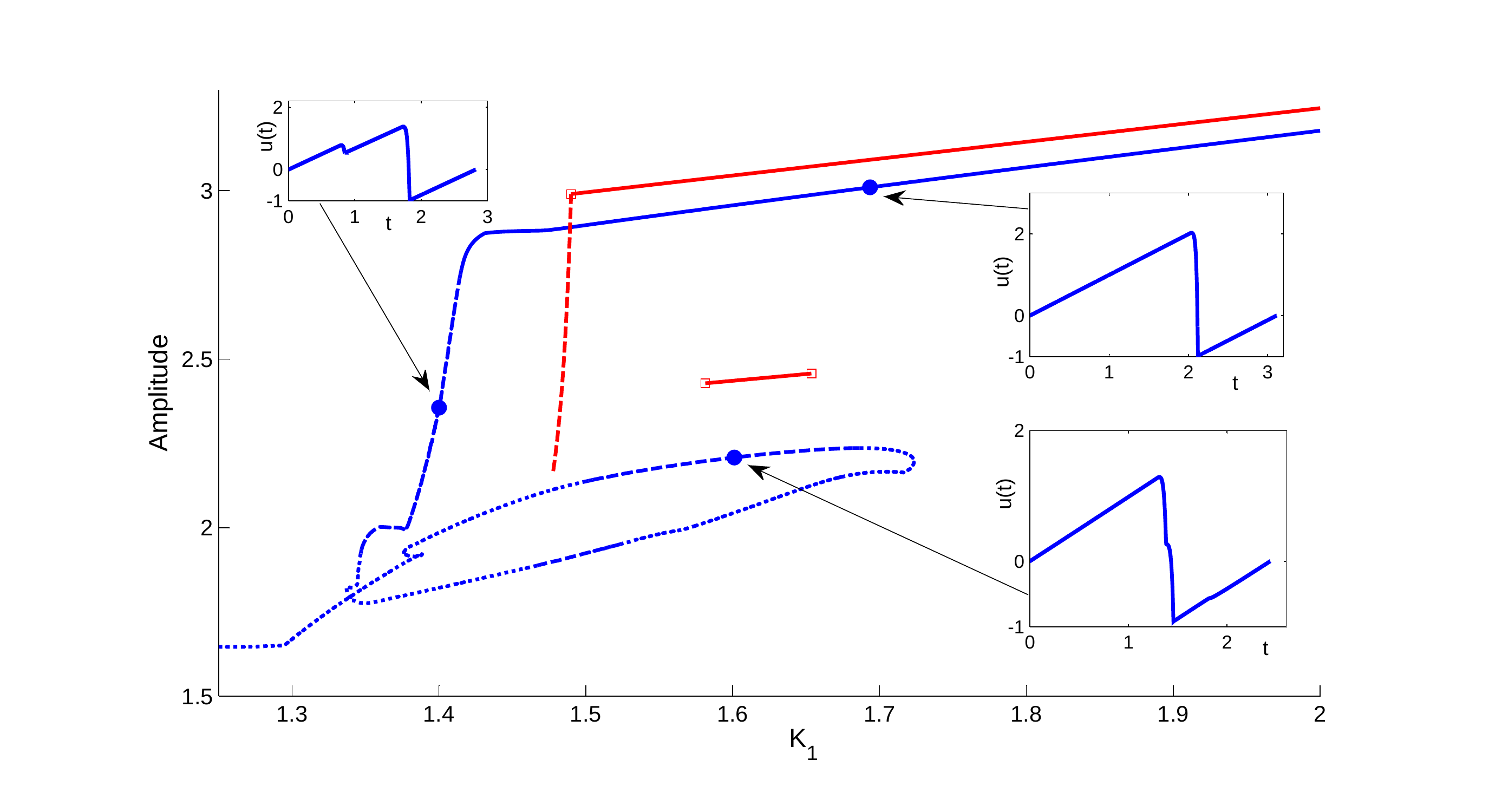}
\put(-240,188){\Large $K_1=1.6010$}
\put(-637,368){\Large $K_1=1.4001$}
\put(-240,315){\Large $K_1=1.6936$}

\put(-420,345){\rotatebox{7}{\Large $\epsilon\!=\!0$, unimodal, $m\!=\!1$}}
\put(-532,253){\Large $\epsilon=0$,\; type II}
\put(-547,238){\Large bimodal,\; $m=2$}
\put(-430,235){\rotatebox{5}{\Large $\epsilon\!=\!0$, unimodal, $m\!=\!2$}}
}
\vspace{-8ex}
\end{center}
\caption{Amplitudes of Type I and II bimodal singular solutions with $m=2$ and unimodal solutions
with $m=1$ for $A=6.98$, $K_2=0.5$ and $n=0$. Also shown is the numerically computed
principal branch of periodic orbits for the same parameters except $\epsilon=0.02$. Unimodal,
bimodal and multimodal periodic orbits for $\epsilon=0.02$ occur on the solid, dashed or dotted parts of
the curve, respectively. Insets with $K_1=1.6936$ and $K_1=1.4001$ show profiles of $\epsilon=0.02$
stable unimodal and type II bimodal periodic orbits. See text for discussion of the $K_1=1.601$
unstable bimodal solution.}
\label{figamp698}
\end{figure}

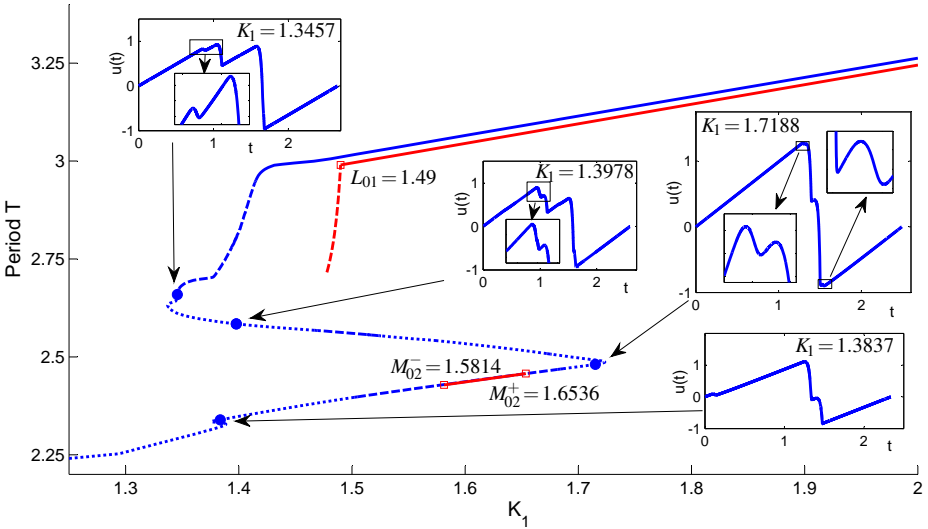
\begin{figure}
\vspace{-5ex}
\begin{center}
\scalebox{0.55}{\mbox{}\hspace{-3cm}\includegraphics{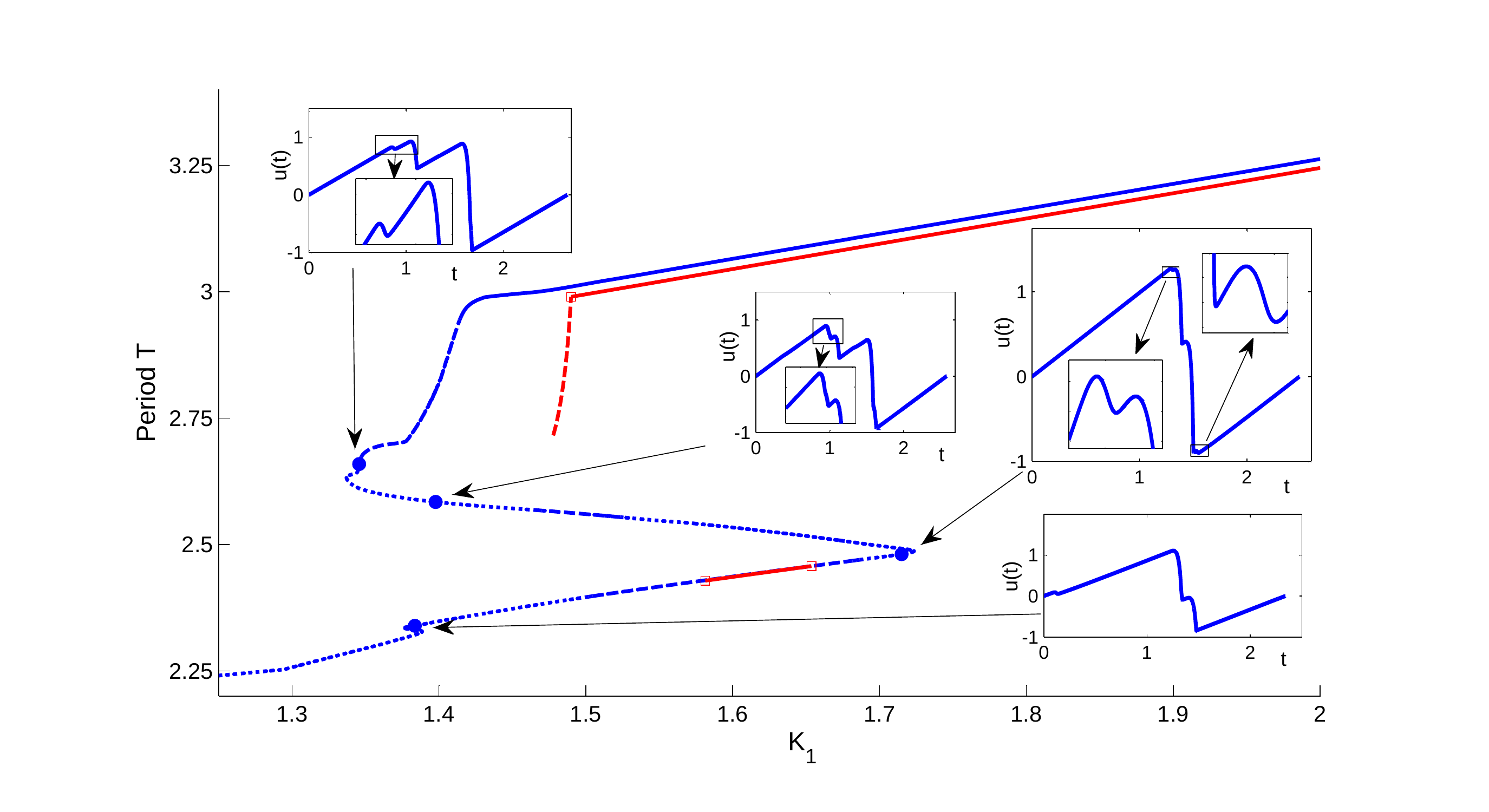}
\put(-175,146){\Large $K_{\!1}\!=\!1.3837$}
\put(-562,363){\Large $K_{\!1}\!=\!1.3457$}
\put(-355,266){\Large $K_{\!1}\!=\!1.3978$}
\put(-240,298){\Large $K_{\!1}\!=\!1.7188$}

\put(-483,262){\Large $L_{01}\!=\!1.49$}
\put(-387,113){\Large $M_{\!02}^+\!=\!1.6536$}
\put(-455,132){\Large $M_{\!02}^-\!=\!1.5814$}
}
\vspace{-8ex}
\end{center}
\caption{Periods of $\epsilon=0$ singular solutions and the numerically computed
$\epsilon=0.02$ solution branch for the same parameters as Fig.~\protect\ref{figamp698}. Insets show profiles of
$\epsilon=0.02$ trimodal and quadrimodal solutions. The solution with $K_1=1.3457$ is
stable, the others are unstable.}
\label{figper698}
\end{figure}

Although, as expected, Figs.~\ref{figA448amp}-\ref{figA448per} do not display a fold bifurcation between the
unimodal and bimodal solutions near $K_1=L_{00}$ with $A=4.48$, we note that two fold bifurcations
are visible earlier on the branch in this case at $K_1\approx1.6571$ and $K_1\approx1.7282$.
These folds are not associated with unimodal solutions but with bimodal and multimodal solutions. An inset for $K_1=1.7248$
in Fig.~\ref{figA448per} shows a periodic solution with 4 local maxima per period
close to one of these folds.

\begin{figure}[t]
\vspace{-5ex}
\begin{center}
\scalebox{0.55}{\mbox{}\hspace{-3cm}\includegraphics{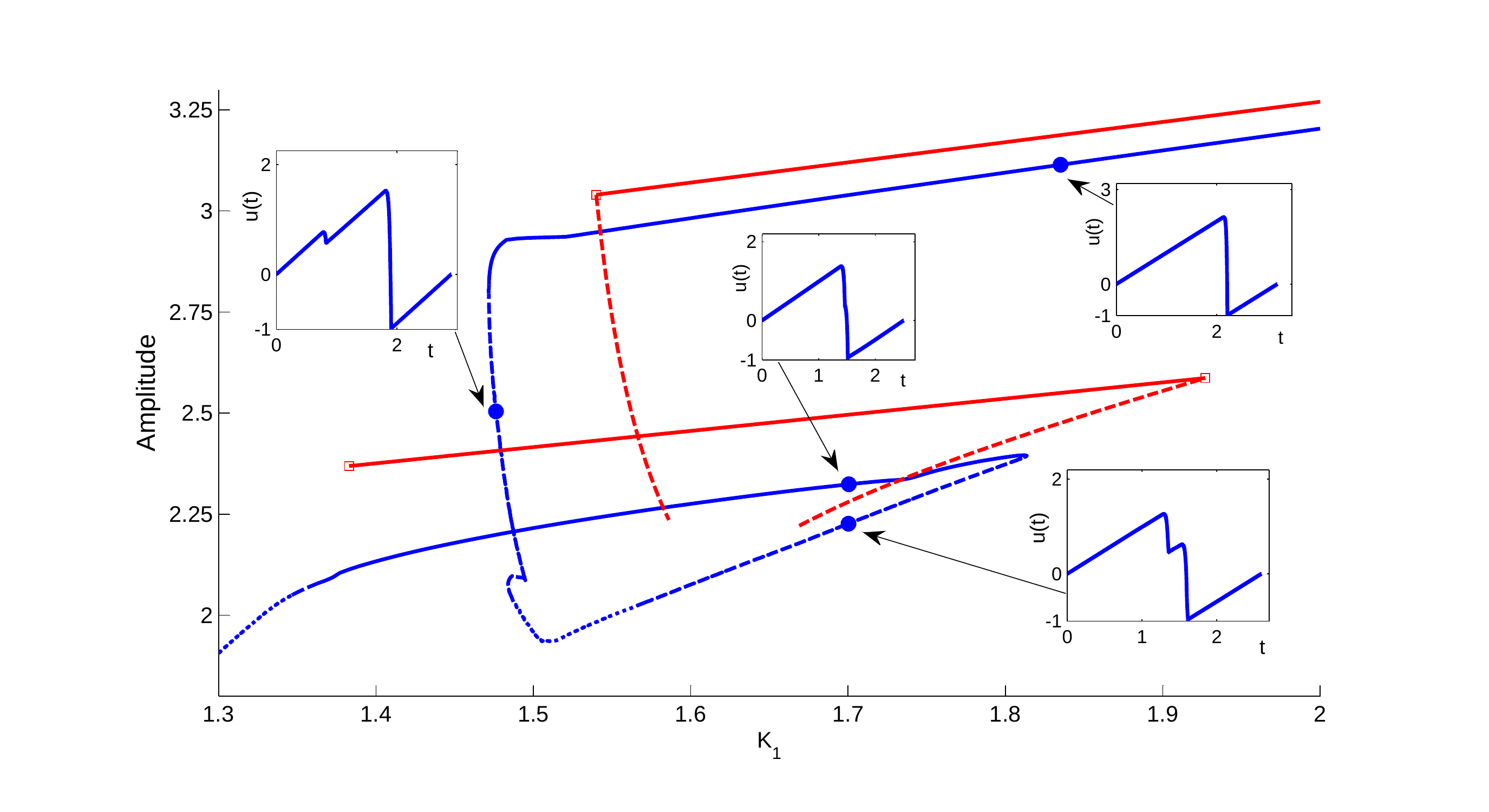}
\put(-623,340){\Large $K_1\!=\!1.4763$}
\put(-176,323){\Large $K_1\!=\!1.8351$}
\put(-191,168){\Large $K_1\!=\!1.7003$}
\put(-378,296){\Large $K_1\!=\!1.7005$}
\put(-420,345){\rotatebox{7}{\Large $\epsilon\!=\!0$, unimodal, $m\!=\!1$}}
\put(-270,205){\Large $\epsilon\!=\!0$,}
\put(-225,200){\Large type I bimodal, $m\!=\!2$}
\put(-502,253){\Large $\epsilon\!=\!0$,\; type II}
\put(-517,238){\Large bimodal,\; $m\!=\!2$}
\put(-285,225){\rotatebox{5}{\Large $\epsilon\!=\!0$, unimodal, $m\!=\!2$}}
}
\vspace{-8ex}
\end{center}
\caption{Amplitudes of Type I and II bimodal singular solutions with $m=2$ and unimodal solutions
with $m=1$ and $m=2$ for $A=7.08$, $K_2=0.5$ and $n=0$. The numerically computed
principal branch of periodic orbits with $\epsilon=0.02$ is also shown.
Insets show profiles of $\epsilon=0.02$ periodic orbits which
correspond to stable unimodal, and unstable type I and type II bimodal solutions.}
\label{figamp708}
\end{figure}

\begin{figure}
\vspace{-5ex}
\begin{center}
\scalebox{0.55}{\mbox{}\hspace{-3cm}\includegraphics{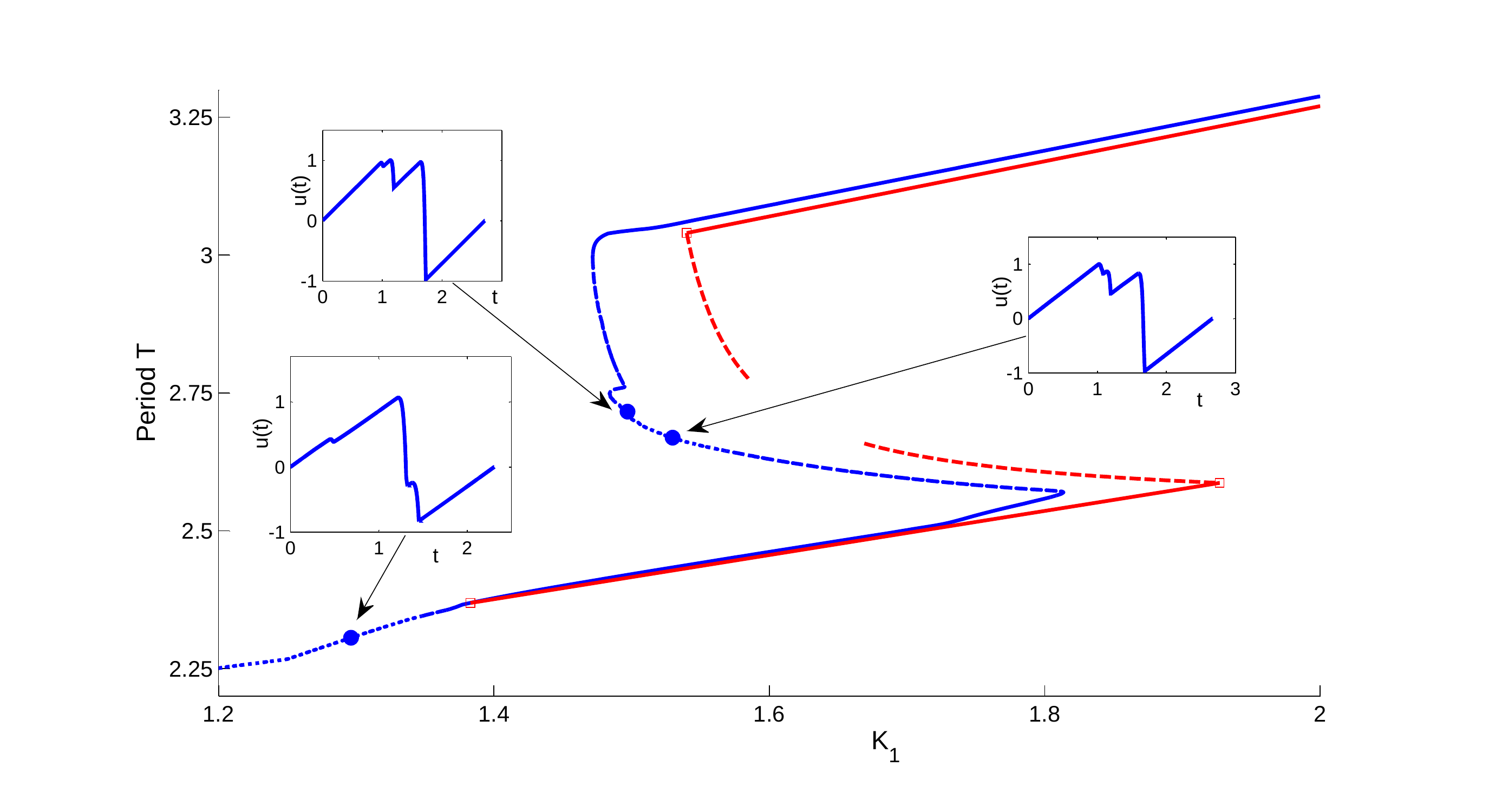}
\put(-600,230){\Large $K_1=1.2962$}
\put(-205,295){\Large $K_{\!1}\!=\!1.5298$}
\put(-600,353){\Large $K_{\!1}\!=\!1.4970$}

\put(-410,298){\Large $L_{01}=1.54$}
\put(-165,155){\Large $M_{02}^+=1.9271$}
\put(-550,95){\Large $M_{02}^-=1.3829$}
}
\vspace{-8ex}
\end{center}
\caption{Periods of $\epsilon=0$ singular solutions and the numerically computed
$\epsilon=0.02$ solution branch for the same parameters as Fig.~\protect\ref{figamp708}. Insets show profiles of
unstable trimodal solutions on the $\epsilon=0.02$ branch.}
\label{figper708}
\end{figure}
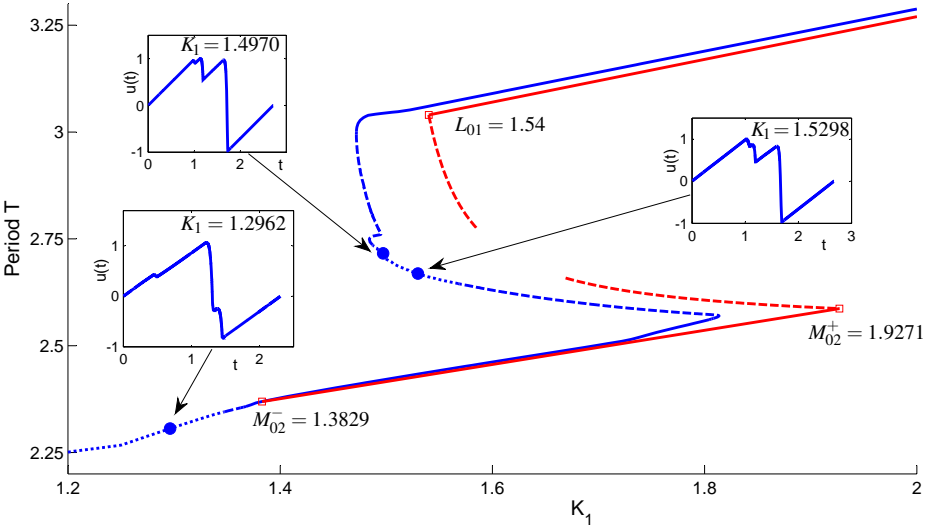

In Figs.~\ref{figamp698}-\ref{figper708} we illustrate the change in the dynamics near to $K_1=L_{01}$ as
$A$ passes through $7$; the second cusp-like bifurcation indicated by \eq{cusp}. Figs.~\ref{figamp698}-\ref{figper698}
demonstrate that for $\epsilon=0.02$ and $A=6.98$ there is a transition from a bimodal to a unimodal solution
close to $K_1=L_{01}$ without a fold bifurcation, while for $A=7.08$ the same transition is associated with a fold
bifurcation. As was the case with the first cusp-like bifurcation for $\epsilon=0.02$ this transition occurs for a value of
$K_1$ slightly less than $L_{01}$ both when $A=6.98$ and $A=7.08$. Also, whereas the fold appears when $A=7$ for
the singular solution with $\epsilon=0$, additional computations
with other values of $A$ (not shown) suggest that for $\epsilon=0.02$ the fold bifurcation associated with the point $K_1=L_{01}$
disappears at about $A=7.04$.

Figs.~\ref{figamp708}-\ref{figper708} also indicate good agreement between the singular and $\epsilon=0.02$
solutions near to the fold point $K_1=M_{02}^+$, with the $\epsilon=0.02$ solution having a fold bifurcation
associated with the solution profile transitioning from unimodal to bimodal with
$K_1\approx1.8137$ slightly less than $M_{02}^+=1.9271$.
The insets with $K_1=1.4970$ and $1.5298$ in Fig.~\ref{figper708} for $A=7.08$ show that trimodal
solutions again occur for $\epsilon=0.02$ in the gap between the two intervals of bimodal solutions,
just as was previously seen for $A=4.54$ in Fig.~\ref{figA454per}.

Figs.~\ref{figamp698}-\ref{figper698} show a significant difference between the
dynamics near to the second cusp-like bifurcation compared to the first one. Although for $\epsilon=0.02$
and $A\leappr7$ there is
no longer a fold bifurcation near to $K_1=L_{01}$ and there are no folds associated with
transitions from unimodal to bimodal solutions, there are still fold bifurcations on the branch. Insets in
Fig.~\ref{figper698} show trimodal solution profiles close to each of these folds.
The $\epsilon=0$ singular solutions also show differences between the first and second cusp-like bifurcation,
since when
$A=4.48$ there are type I bimodal solutions for $K_1<L_{00}$ whereas for
$A=6.98$ there are no type I bimodal solutions, but there is a small interval of
$m=2$ unimodal solutions which coexist with the $m=1$ unimodal solutions. Fig.~\ref{figper698} shows
that the $\epsilon=0.02$ branch has solutions whose $K_1$ values and periods almost exactly agree with those
of the unimodal $m=2$ singular solutions while the inset for $K_1=1.601$ in Fig.~\ref{figamp698} shows that
while the $\epsilon=0.02$ solution has smaller amplitude and is bimodal, its profile is close to unimodal.

The inset for $K_1=1.2962$ shows a trimodal solution for $\epsilon=0.02$ occurring before the pair of fold bifurcations.
Such a solution is also seen in Fig.~\ref{figA454per}, and they seem to be ubiquitous, also arising even
when the folds disappear (see inset for $K_1=1.5997$ in Fig.~\ref{figA448per}).

\section{Other Solutions and Bifurcations}
\label{sec:othersols}

\begin{figure}[ht]
\vspace{-5ex}
\mbox{}\hspace{-1.5cm}\scalebox{0.55}{\includegraphics{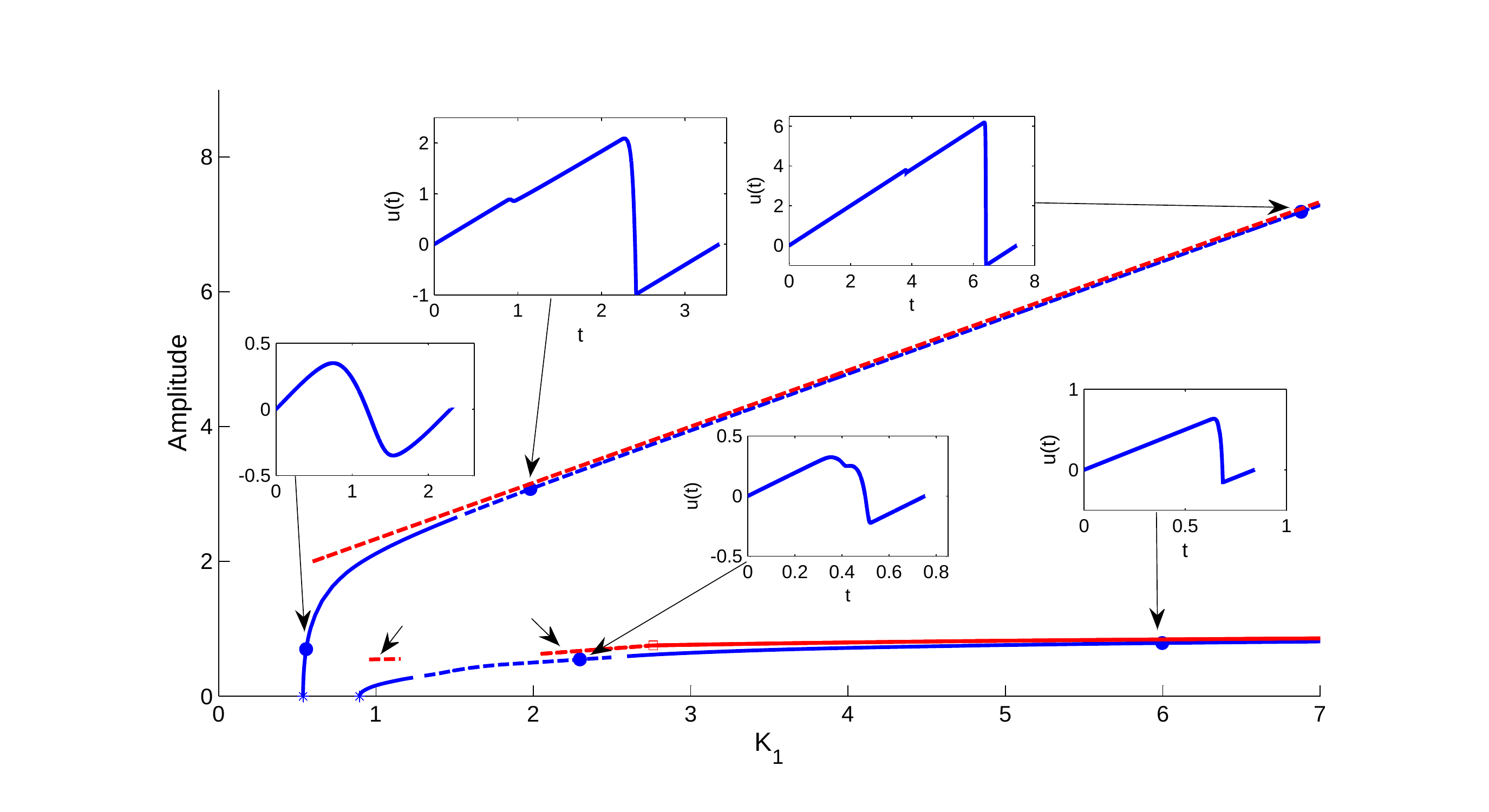}
\put(-610,237){\Large $K_{\!1}\!=\!0.5570$}
\put(-555,355){\Large $K_1=1.9817$}
\put(-365,357){\Large $K_1\!=\!6.8797$}
\put(-392,142){\Large $K_{\!1}\!=\!2.2963$}
\put(-214,168){\Large $K_{\!1}\!=\!5.9966$}

\put(-605,120){\Large $\epsilon\!=\!0$ bimodal, $n\!=\!1$, $m\!=\!0$}
\put(-540,106){\Large type I}
\put(-600,105){\Large type II}

\put(-480,195){\rotatebox{20}{\Large $\epsilon\!=\!0$, type II bimodal, $n\!=\!m\!=\!0$}}
\put(-380,95){\Large $\epsilon\!=\!0$, unimodal, $n\!=\!1$, $m\!=\!0$}
\put(-450,70){\Large $M_{\!10}^+\!=\!2.7625$}
}
\vspace{-5ex}
\caption{The first two bifurcation branches with $A=7/6$, $\epsilon=0.05$ and $\gamma=a_1=c=1$, $K_2=0.5$.
Also shown are the corresponding $\epsilon=0$ singular solution branches, and insets show
profiles of periodic solutions for $\epsilon=0.05$ at different points on the branches. The $\epsilon=0.05$ solutions
on the first branch are all stable, and those on the second branch all unstable; DDEBiftool does not detect any
secondary bifurcations.}
\label{figamp76}
\end{figure}
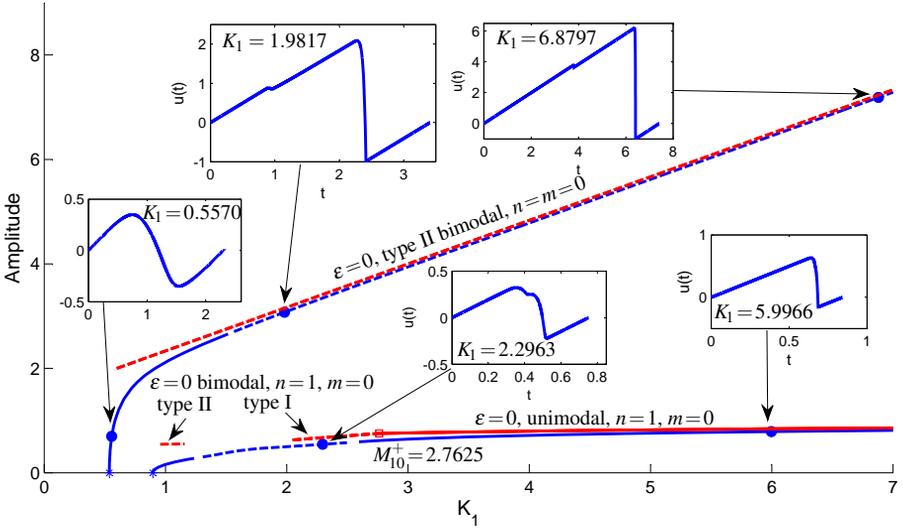

In this section we study some of the other solutions and bifurcations that can arise with \eq{eps2del}.
In Sections~\ref{secnumerics} and~\ref{sec:cusp} we were mainly concerned with
the fold and cusp-like bifurcations predicted by Theorems~\ref{thmIbif}
and~\ref{thmIIbif}. The folds occurred between legs of unimodal solutions, and as noted in the discussion
after Theorem~\ref{thmunilegs}, we require $A>3$ to have more than one leg of unimodal solutions on the
principal $n=0$ branch of periodic solutions. So we begin this section
by considering the dynamics when $A\in(1,3)$.
Noting that Theorem~\ref{thmunilegs} guarantees the existence of unimodal solutions for all $K_1$
sufficiently large \emph{unless} $m=m^0(n)=n(A-1)$ and $A\leq1+K_2/(1+K_2)$ we first consider this exceptional case.
On the principal branch this occurs when $m=n=0$ and taking $K_2=0.5$ with $A\in[1,4/3]$. Consequently
we consider the dynamics with $A=7/6$, as shown in Fig.~\ref{figamp76}.

Verifying the conditions of Theorem~\ref{thm2delII} we find that there is a type II bimodal
solution with $n=m=0$ when $A=7/6$ for all $K_1>0.59816$ which is shown in Fig.~\ref{figamp76}.
With $\epsilon=0.05$,
DDEBiftool finds a Hopf bifurcation at $K_1\approx0.5373$ leading to a branch of stable periodic solutions which exist
for all larger values of $K_1$. Close to the Hopf bifurcation these solutions are unimodal and
sinusoidal, but for all $K_1>1.5167$ these solutions have two local maxima per period, and closely
resemble type II bimodal solutions (see $K_1=1.9817$ and $K_1=6.8797$ insets in Fig.~\ref{figamp76}).
There is also very good agreement for $K_1>1$ between the amplitude of the $\epsilon=0$ singular solutions
given by Theorem~\ref{thm2delII} and the numerically found $\epsilon=0.05$ solutions.
Type II bimodal singular solutions and their $\epsilon>0$ counterparts are also
found for all $K_1$ sufficiently large for other values of $A\in(1,1+K_2/(1+K_2))$ when $m=m^0(n)=n(A-1)$.


Fig.~\ref{figamp76} also shows the $n=1$ branch for $A=7/6$.
By Theorem~\ref{thmunilegs}(i)
there is a unimodal singular solution with $n=1$, $m=0$ for all $K_1>M_{10}^+=2.7625$, while
verifying the conditions of Theorem~\ref{thm2delI} reveals that there is a type I bimodal solution
for $K_1\in(2.0481,M_{10}^+)$. At $K_1=M_{10}^+$ the two solutions coincide, with $T_2\to0$ as $K_1\to M_{10}^+$
for the type I bimodal solution. Theorem~\ref{thmIbif} deals with unimodal and type I bimodal solutions
coinciding at $K_1=M_{nm}^+$ in a fold-like bifurcation. That theorem does not apply here because we have
$m=0<n(A-1)=1/6$ outside its range of validity, nevertheless we still have a transition between the
two types of solutions, but here it occurs without the fold-like bifurcation. With $\epsilon=0.05$,
DDEBiftool finds all the solutions on the corresponding branch to be unstable with a transition between
bimodal and unimodal solutions close to  $K_1=M_{10}^+$.

\begin{figure}
\vspace{-5ex}
\mbox{}\hspace{-1.5cm}
\scalebox{0.55}{\includegraphics{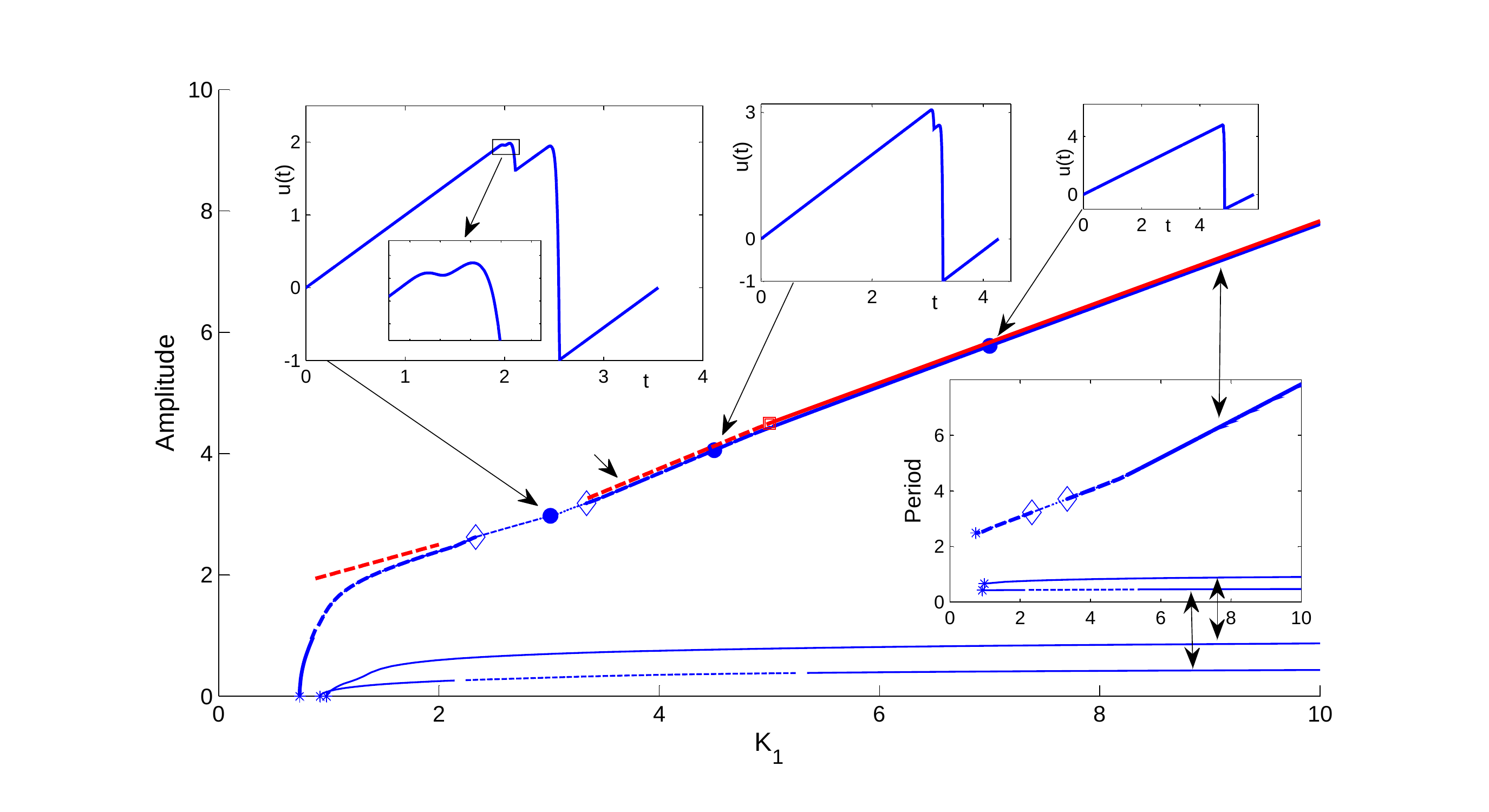}
\put(-615,360){\Large $K_1=3.0137$}
\put(-375,290){\Large $K_1\!=\!4.4998$}
\put(-213,366){\Large $K_{\!1}\!=\!6.9996$}

\put(-538,130){\Large $D_{1}=2.3347$}
\put(-475,150){\Large $D_{2}=3.3414$}
\put(-380,190){\Large $M_{00}=5$}

\put(-249,263){\rotatebox{20}{\Large $\epsilon=0$, unimodal, $n=m=0$}}

\put(-667,138){\rotatebox{13}{\Large $\epsilon=0$, type II bimodal}}
\put(-630,131){\rotatebox{13}{\Large $n=m=0$}}

\put(-542,207){\Large $\epsilon=0$, type I bimodal}
\put(-515,193){\Large $n=m=0$}
}
\vspace{-5ex}
\caption{Amplitude and period plots for the first three branches of periodic solutions with insets
showing unimodal, bimodal and trimodal solutions on the prinicpal branch for
$A=a_2=1.5$, $\epsilon=0.05$ and $\gamma=a_1=c=1$, $K_2=0.5$. DDEBiftool detects period doubling
bifurcations at $K_1=D_1$ and $D_2$, and the solutions on the principal branch are stable
except for $K_1$ between these values. Also shown are the unimodal and type I bimodal singular solutions
corresponding to the principal branch.}
\label{figampA15eps005}
\end{figure}

\begin{figure}
\vspace{-5ex}
\mbox{}\hspace{-1.2cm}\scalebox{0.55}{\includegraphics{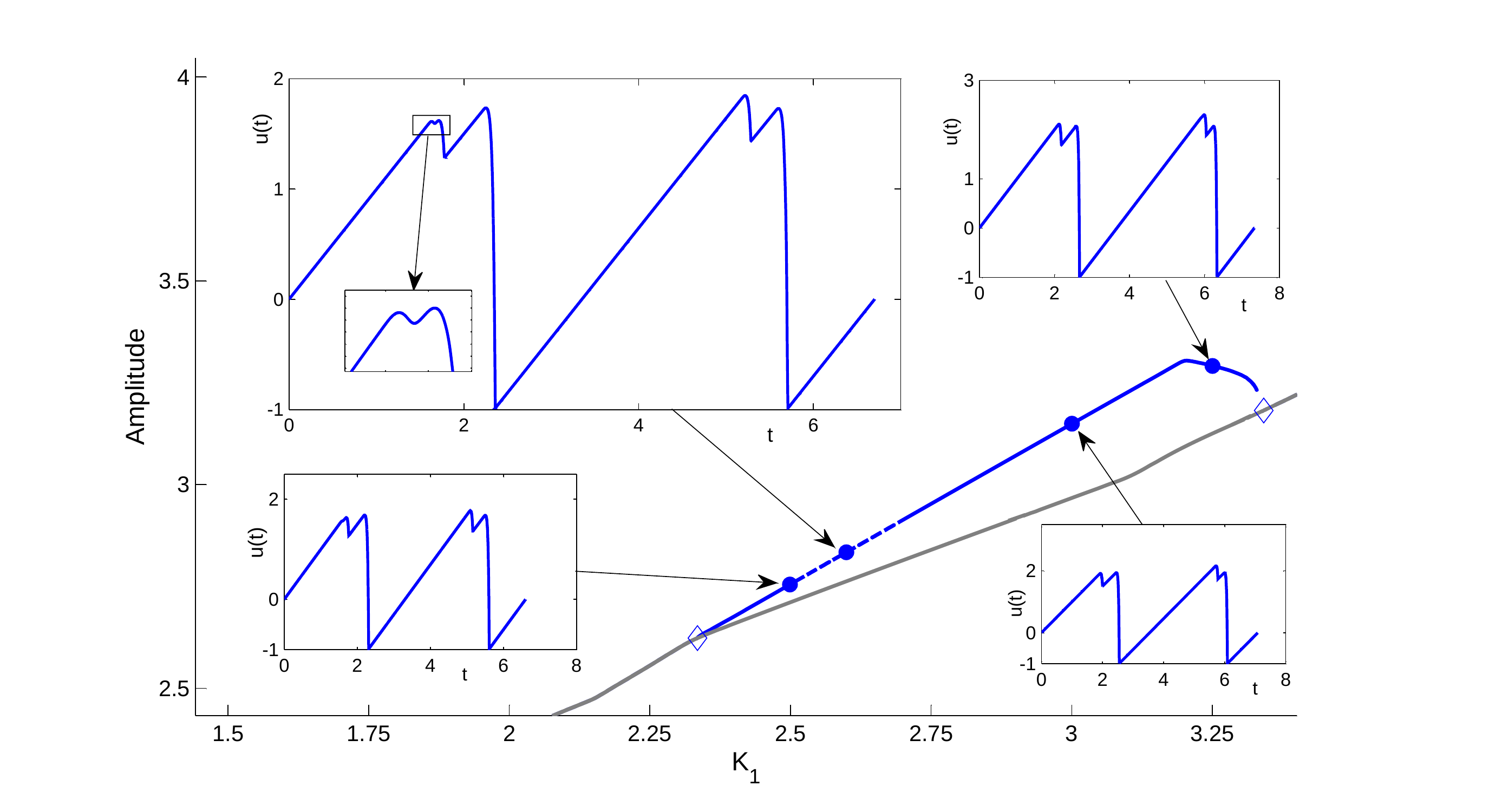}
\put(-510,370){\Large $K_1=2.5994$}
\put(-225,376){\Large $K_1=3.2504$}
\put(-605,167){\Large $K_1\!=\!2.4990$}
\put(-210,138){\Large $K_1=3.0005$}

\put(-415,85){\Large $D_{1}=2.3347$}
\put(-135,196){\Large $D_{2}=$}
\put(-145,180){\Large $3.3414$}
}
\vspace{-5ex}
\caption{Branch of stable period doubled solutions for $K_1\in(D_1,D_2)$ with insets showing solution profiles.
The parameters are the same as in Fig.~\ref{figampA15eps005}, and the principal branch of periodic
solutions from Fig.~\ref{figampA15eps005} is redrawn in gray here.}
\label{figamppdA15eps005}
\end{figure}
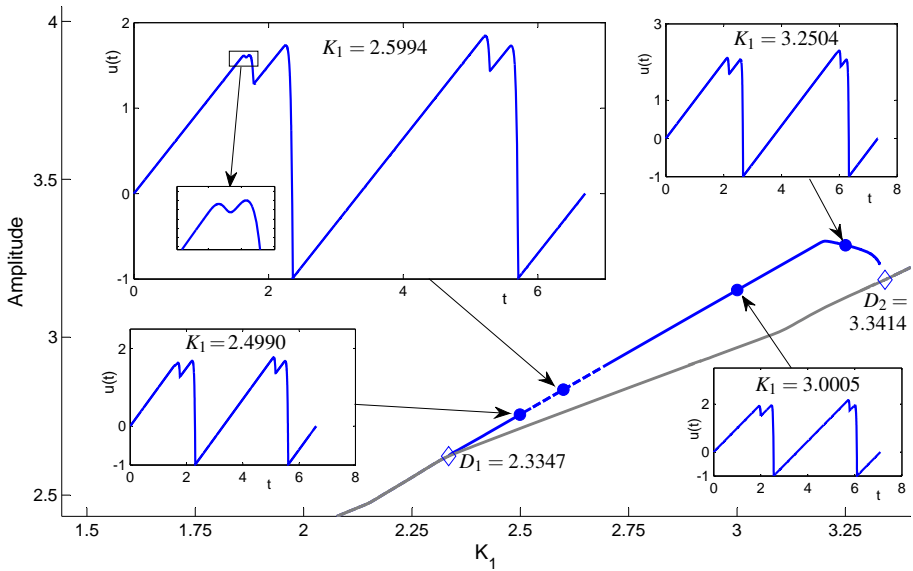

Next we consider $A\in(4/3,3)$ for which Theorem~\ref{thmunilegs} guarantees the existence of a unimodal
singular solution with $n=m=0$ on the principal branch for all $K_1$ sufficiently large, but for which with
$n=0$ there is no value of $m$ which satisfies the conditions of Theorem~\ref{thmIbif} or~\ref{thmIIbif}, and so
we do not expect fold bifurcations of periodic orbits.
Taking $A=1.5$, Theorem~\ref{thmunilegs}(i) gives the existence of unimodal singular solutions
with $n=m=0$ for $K_1>M_{00}=5$. Similarly to the $n=1$ branch of the previous example, verifying the
conditions of Theorem~\ref{thm2delI} we find a type I bimodal singular solution with $n=m=0$ for
$K_1\in(3.3508,M_{00})$.
For $\epsilon=0.05$, using DDEBiftool
we numerically compute the principal solution branch from the
Hopf bifurcation (at $K_1\approx0.7363$), finding stable bimodal periodic solutions for
$K_1\in(3.3414,5.0543)$, and stable unimodal solutions for $K_1>5.0543$ as shown in Fig.~\ref{figampA15eps005}.
The parameter ranges and amplitudes of the solutions with $\epsilon=0.05$ are seen to be very close to the
$\epsilon=0$ singular solutions. Solutions are also stable on the initial part of the branch and
unimodal for $K_1\in(0.7363,0.8586)$ and bimodal for $K_1\in(0.8586,2.3347)$.
However, the solutions on the
principal branch are unstable in the range $K_1\in(D_1,D_2)=(2.3347,3.3414)$ with DDEBiftool detecting period
doubling bifurcations (characterized by a real Floquet multiplier passing through $-1$) at both ends of
this interval. On the principal branch trimodal solutions are found for
$K_1\in(2.9016,3.3112)$ while the solutions are bimodal in the rest of the interval $(D_1,D_2)$.

In Fig.~\ref{figamppdA15eps005} we show the resulting branch of $\epsilon=0.05$ stable period-doubled solutions for
$K_1\in(D_1,D_2)$. The branch is computed from $K_1=D_1$ and appears to terminate at $K_1=D_2$, though
numerical computation of the branch is very difficult near to $K_1=D_2$. Insets show profiles of the
resulting stable periodic solutions which all have period close to $7$ and mainly have four local maxima
per period, except for $K_1\in(2.5247,2.6909)$ where the first peak splits into two
(reminiscent of the transitions from bimodal to trimodal solutions seen in Section~\ref{secnumerics})
resulting in periodic solutions with five local maxima per period.

We do not have a characterization from the singular solutions of when to expect period doubling bifurcations.
To determine the parameter ranges where period doubled orbits can occur
we could parameterise the period doubled singular periodic orbits. This would be
similar to a perturbation of two copies of the parameterisations illustrated in Figs.~\ref{figtypeI}
and~\ref{figtypeII} and
would involve twenty parametrisation intervals and more algebraic manipulation than we care to
contemplate. We note that at the end of the interval of validity of
the type I bimodal solution at $K_1=3.3508$ the left inequality in \eq{tIthbds} is tight, and fails for smaller
values of $K_1$. This is the same inequality that failed at the transition between
type I bimodal and trimodal solutions between the fold bifurcations at $K_1=L_{nm-1}$ and $K_1=M_{nm}^+$
which we studied in Section~\ref{secnumerics}.
Given the proximity of the end of the interval of type I bimodal singular solution at $K_1=3.3508$
to the period doubling bifurcation with $\epsilon=0.05$ at $K_1=3.3414$ and the transition from
bimodal to trimodal solutions at $K_1=3.3112$ we suspect that in the limit as $\epsilon\to0$ the singular solutions
undergo a period doubling bifurcation at the same parameter value where the periodic solution transitions
from type I bimodal to trimodal. The behaviour seen in this example contrasts with the examples
in Sections~\ref{secnumerics} and~\ref{sec:cusp} where
no period doubling bifurcations were detected between the fold bifurcations.


\begin{figure}
\vspace{-5ex}
\mbox{}\hspace{-1.5cm}\scalebox{0.55}{\includegraphics{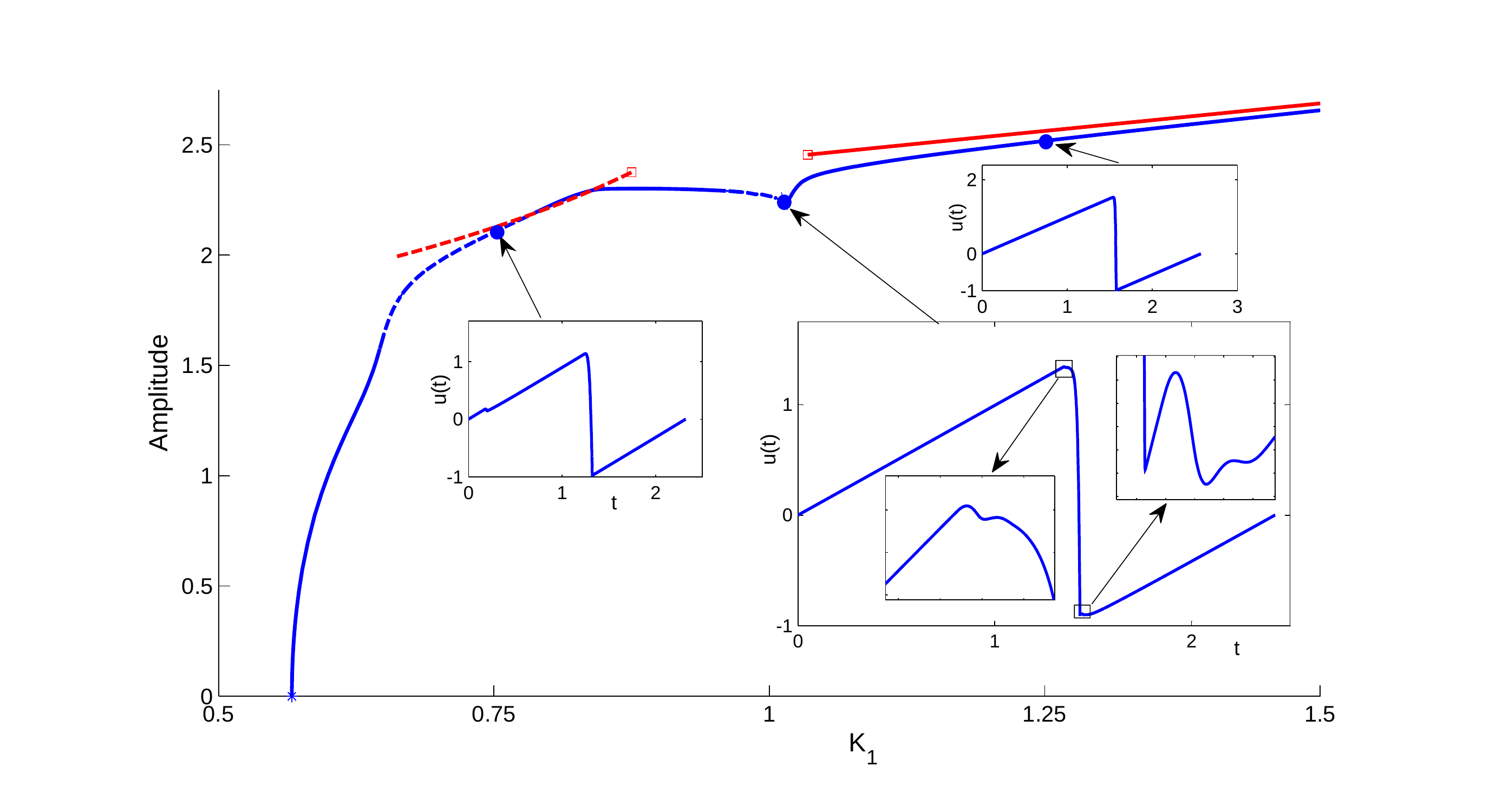}

\put(-392,360){\Large $M_{01}^-\!=\!1.0348$}
\put(-450,340){\Large $L_{01}\!=\!0.875$}

\put(-225,370){\rotatebox{6}{\Large $\epsilon\!=\!0$, unimodal, $m\!=\!1$}}
\put(-610,302){\rotatebox{18}{\Large $\epsilon\!=\!0$, type II bimodal, $m\!=\!2$}}

\put(-541,250){\Large $K_{\!1}\!=\!0.7522$}
\put(-267,333){\Large $K_{\!1}\!=\!1.2504$}
\put(-360,245){\Large $K_{\!1}\!=\!1.0115$}
}
\vspace{-6ex}
\caption{Amplitude plot of type II bimodal solution with $n=0$ and $m=2$ which exists for $K_1<L_{01}<1$
and a unimodal singular solution with $n=0$ and $m=1$ which exists for $K_1>M_{01}^->1$. Here $A=a_2=5.75$ and
as before $K_2=0.5$, $a_1=1$ and $c=1$. Also shown is the corresponding branch of periodic orbits
with $\epsilon=0.005$ computed using DDEBiftool. This branch bifurcates from the steady state solution at $K_1\approx0.5665$.
Insets show examples of $\epsilon=0.005$ stable solutions for $K_1<L_{01}$ and $K_1>M_{01}^-$ and an unstable solution
for $K_1=1.0115\in(L_{01},M_{01}^-)$.}
\label{figA575septypeII}
\end{figure}
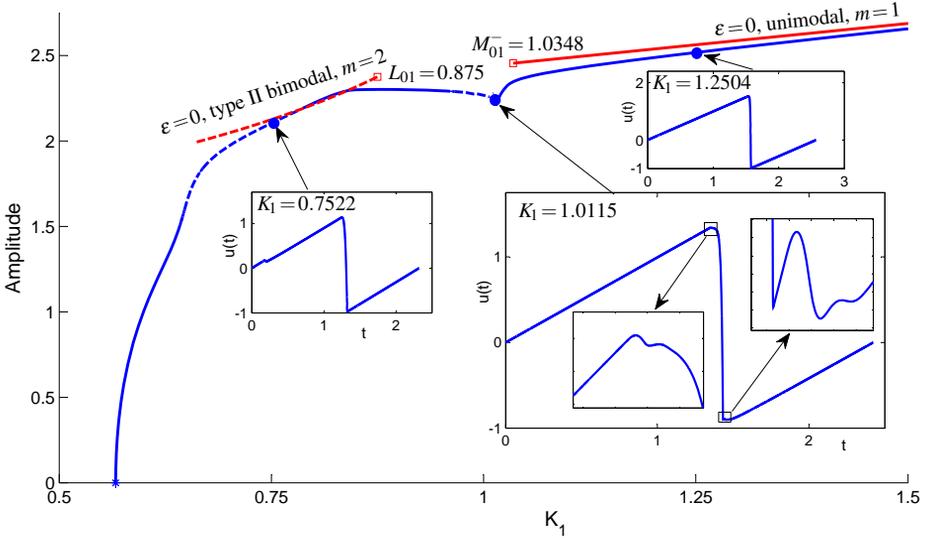

We have already seen examples corresponding to Theorem~\ref{thmIIbif}(i) and (ii) with unimodal and type II bimodal
solutions which coincide at $K_1=L_{nm-1}$ with or without a fold bifurcation. Fig.~\ref{figA575septypeII}
illustrates Theorem~\ref{thmIIbif}(iii), showing a type II bimodal singular solution
which exists for $K_1<L_{nm-1}$ and a unimodal solution for $K_1>M_{nm-1}^-$, where
$M_{nm-1}^->1>L_{nm-1}$. When $n=0$ have $L_{nm^*}=1$ at $A=1+m(2+K_2)$ so $A=3+K_2,5+2K_2,7+3K_2,\ldots$
and so separated unimodal and type II bimodal solutions
will occur for $A$ slightly smaller than these values.
In Fig.~\ref{figA575septypeII} we consider $K_2=0.5$ and $A=5.75<6=5+2K_2$.

Theorems~\ref{thm2delunimod} and~\ref{thm2delI}
require $K_1>1$ for unimodal and type I bimodal solutions to exist, but Theorem~\ref{thm2delII} only requires
$K_1>1-K_2$ for the construction type II bimodal solutions, and Fig.~\ref{figA575septypeII} shows an
example of type II bimodal solutions which exist for $K_1<1$. Fig.~\ref{figA575septypeII} also shows a numerically computed
branch of periodic orbits with $\epsilon=0.005$ which passes very close to the legs of bimodal type I and
unimodal singular solutions. While the type II bimodal solutions exist for $K_1\in(0.6621,0.875)$
the $\epsilon=0.005$ branch has bimodal solutions
for $K_1\in(0.6444,0.7822)$ with the inset solution profile for $K_1=0.7522$ showing that these resemble type II singular solutions.
The $\epsilon=0.005$ branch also has unimodal solutions for $K_1>1.0126$ which approximate the unimodal singular solutions
existing for $K_1>1.0348$. It was found that the period of the solutions on the $\epsilon=0.005$
branch increases monotonically
from $T\approx2.2478$ at the Hopf bifurcation and crosses $T=2.375$ at
$K_1\approx0.84265$. At this value of $K_1$ the period $T$ satisfies $2T=4.75=a_2-a_1$,
that is the difference between the delays is exactly two periods. For
the singular solutions $2T=a_2-a_1$ when $K_1=L_{01}=0.875$ at the end of the interval of bimodal type II solutions.
Fig.~\ref{figA575septypeII} suggests that there is not a bifurcation near to $K_1=L_{nm-1}<1$ when the conditions
of Theorem~\ref{thmIIbif}(iii) are satisfied, but neither do the solutions transition directly from type II bimodal to
unimodal solutions, as occurs in Theorem~\ref{thmIIbif}(i) and (ii). We do not have an explanation for the
dynamics for $K_1\in(L_{01},M_{01}^-)$ in Fig.~\ref{figA575septypeII}, but note with $\epsilon=0.005$
bimodal solutions are seen for $K_1\in(0.9606,1.0066)$ and multimodal solutions for $K_1\in(1.0066,1.0126)$,
and the solution transitions directly from multimodal to unimodal at the kink in the bifurcation branch with
$K_1=1.0126$.

Another example where bimodal type II singular solutions could exist for $K_1<1$ was already seen in
Fig.~\ref{figAmpeps005}. Fig.~\ref{figAmpeps005} illustrated the boundary between Theorem~\ref{thmIIbif}(ii) and (iii)
with $m=m^{*\!}(n)$ and $L_{nm}=M_{nm}^-=1$.

\begin{figure}
\vspace{-5.25ex}
\mbox{}\hspace{-1.5cm}\scalebox{0.55}{\includegraphics{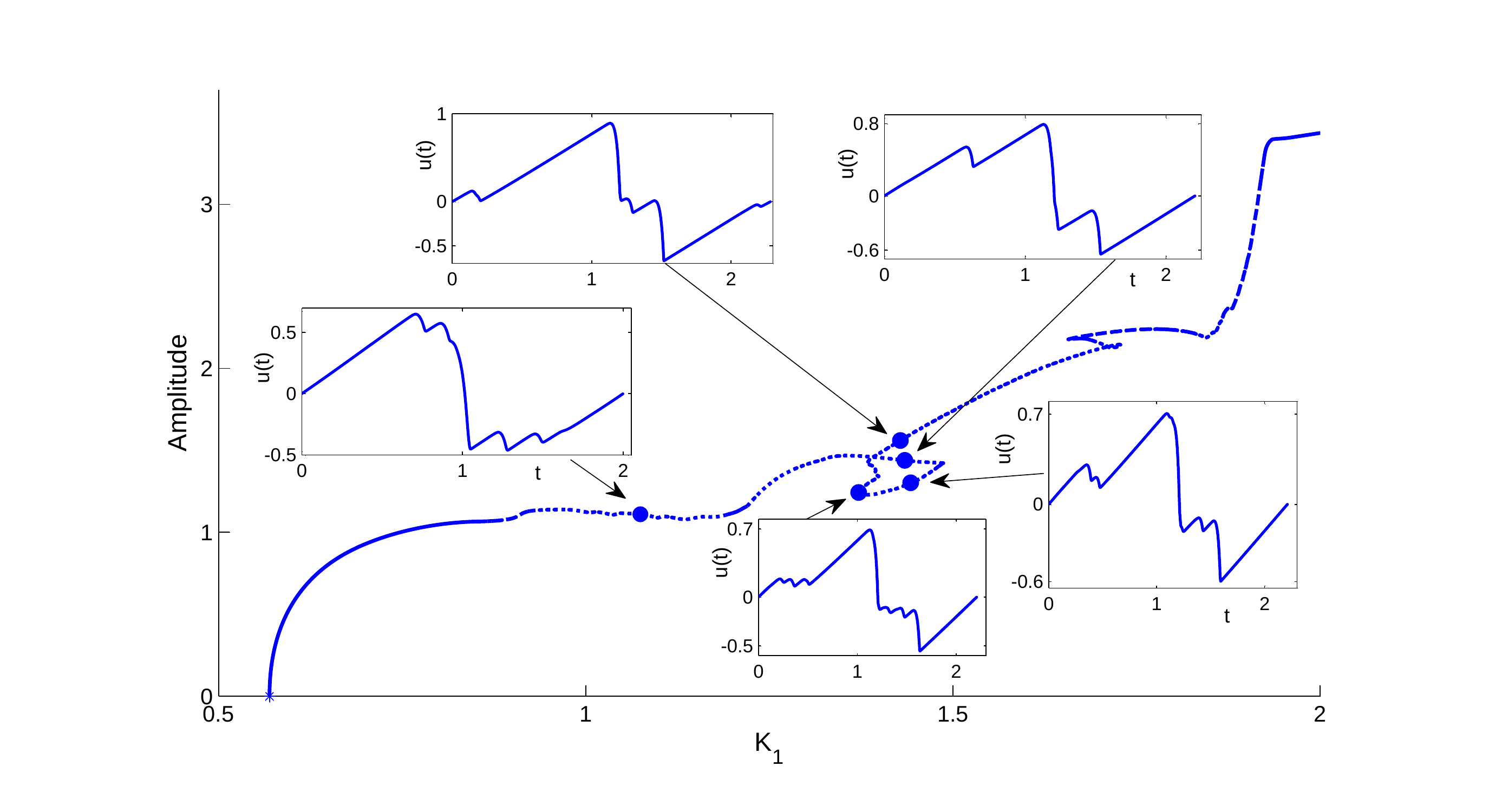}
\put(-525,254){\Large $K_1\!=\!1.074$}
\put(-220,356){\Large $K_1\!=\!1.435$}
\put(-232,124){\Large $K_1\!=\!1.443$}
\put(-385,88){\Large $K_1\!=\!1.372$}
\put(-450,358){\Large $K_1\!=\!1.429$}
}
\vspace{-6ex}
\caption{Amplitude plot of principal branch of periodic solutions with $\epsilon=0.02$, $A=a_2=4.48$,
$a_1=c_1=1$ with insets showing profiles of unstable multimodal solutions that occur on the branch. Recall that dotted
lines indicate multimodal solutions, dashed indicate bimodal and solid lines show unimodal solutions.}
\label{figampA448bizarre}
\end{figure}
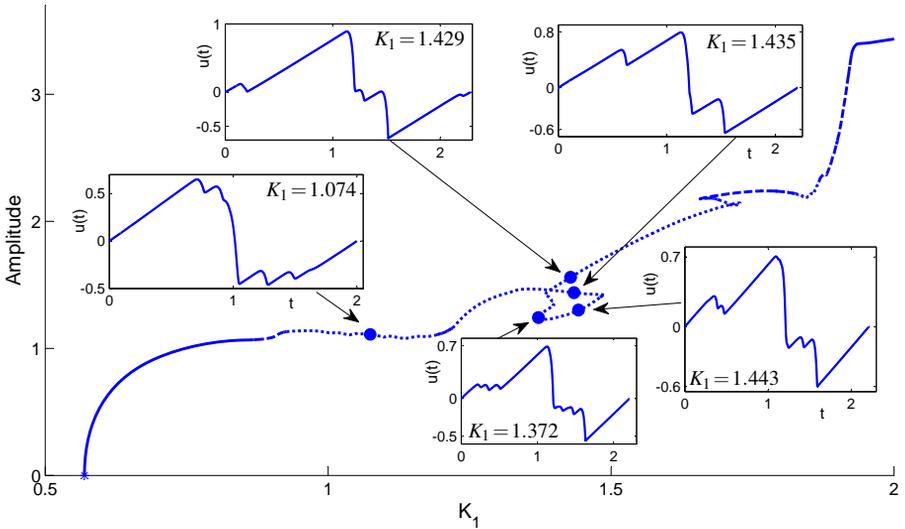

Thus far, we have concentrated our attention on unimodal and bimodal solutions, but noted that trimodal and
quadrimodal solutions arise between legs of type I and II bimodal solutions. Fig.~\ref{figampA448bizarre}
shows examples of multimodal solutions with up to seven local minima
per period (see the $K_1=1.372$ inset) for $\epsilon=0.02$. The parameters in  Fig.~\ref{figampA448bizarre}
are the same as those considered in Figs.~\ref{figA448amp}-\ref{figA448per}, where we studied the cusp-like
bifurcation at $K_1=L_{00}=2$ with $A=4.5$. Fig.~\ref{figampA448bizarre} shows that
even for $A<4.5$ when there are no fold bifurcations near to $K_1=L_{00}$,
there are still six fold bifurcations earlier on the principal branch, and there are solutions with multimodal profiles
near to each of these folds. The multimodal solution profiles shown in the figure for $K_1\in(1.3,1.44)$
all have well-defined `sawteeth' with the periodic solution profile having gradient close to $1/c=1$ before
each local maxima and large negative gradient afterwards. It seems plausible that the
fold bifurcations associated with the transitions between unimodal and bimodal solution profiles that
we studied earlier are just the simplest example of a sequence of such bifurcations that occur at points
where the number of local maxima in the periodic solution profile changes. In principle, Definition~\ref{defadmiss}
and our techniques could be used to locate such bifurcations in the $\epsilon\to0$ limit.

\begin{figure}
\mbox{}\hspace{-0.5cm}\scalebox{0.35}{\includegraphics{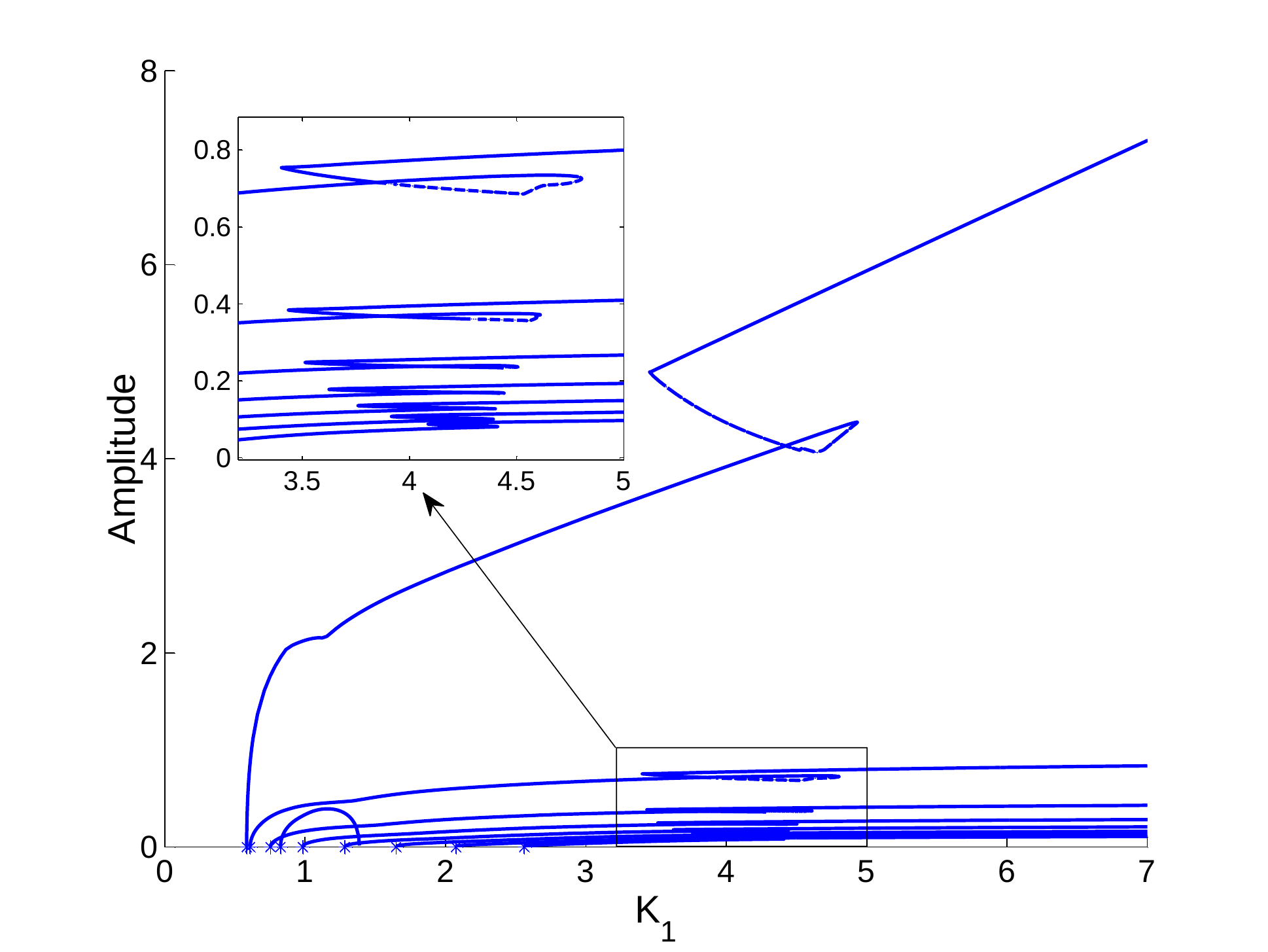}\hspace{-1cm}\includegraphics[trim=1cm 0cm 0cm 0cm,clip]{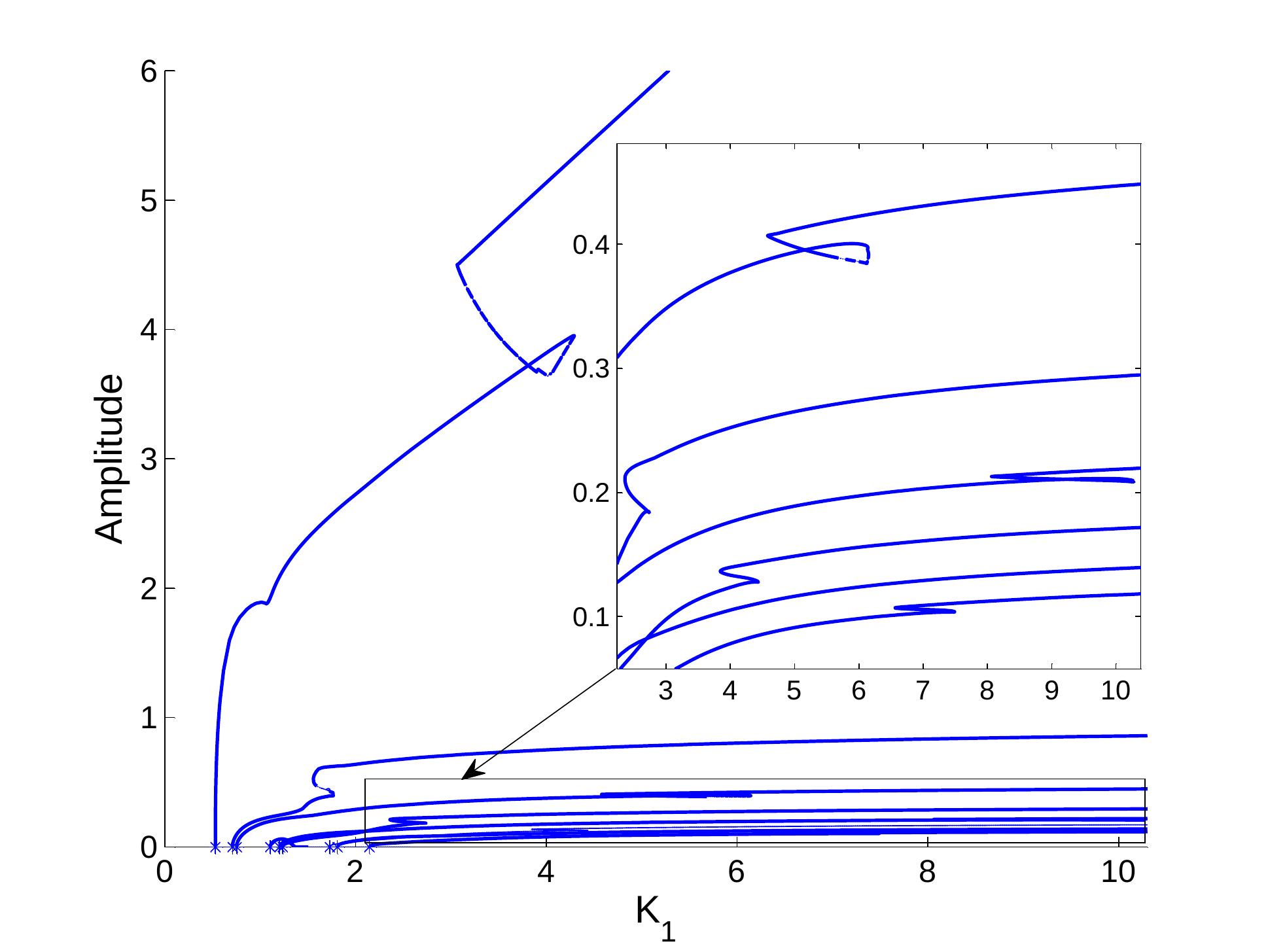}}
\put(-150,130){(ii)}
\put(-260,130){(i)}
\caption{Amplitudes of several branches of periodic solutions with $\epsilon=0.05$, $K_2=0.5$, $a_1=c=1$
and (i) $A=6$  and (ii) $A=(9+\sqrt{5})/2$.}
\label{figA6align}
\end{figure}

\begin{figure}
\vspace{-5ex}
\mbox{}\hspace{-1.5cm}\scalebox{0.55}{\includegraphics{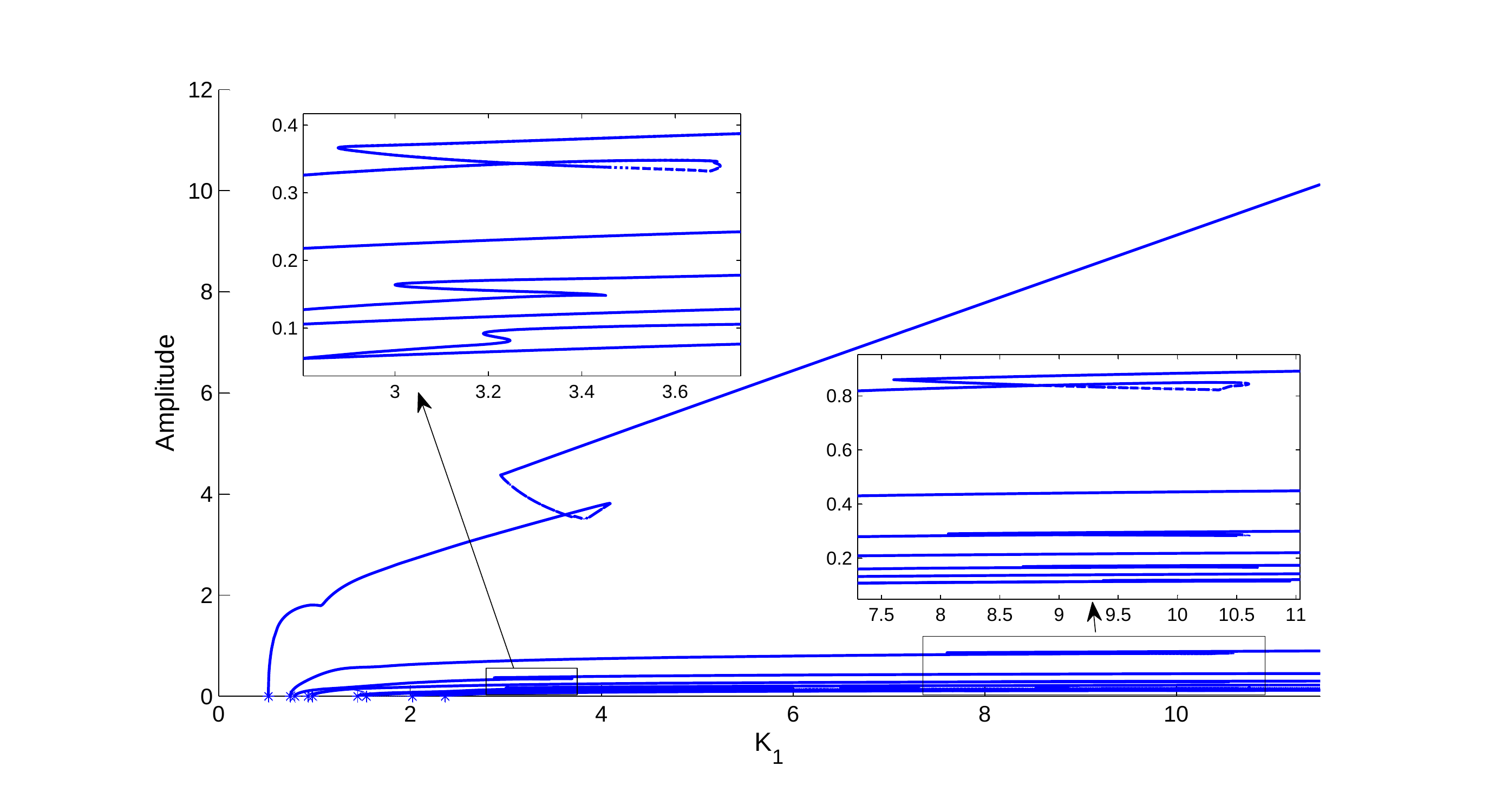}}
\vspace{-5ex}
\caption{Amplitudes of several branches of periodic solutions with $\epsilon=0.05$, $K_2=0.5$, $a_1=c=1$
and $A=5.5$.}
\label{figA55align}
\end{figure}

We have mainly considered the principal branch of periodic solutions corresponding to singular solutions
with $n=0$, since this is the branch on which stable solutions can be observed, but
it was demonstrated
in \cite{DCDSA11} that there are infinitely many Hopf bifurcations for $\epsilon>0$
and we finish this work by considering the
alignment of the bifurcations on the different branches.
This is illustrated in Figs.~\ref{figA6align}
and~\ref{figA55align} for $\epsilon=0.05$.
With $A=6$ we saw earlier that in the limit as $\epsilon\to0$ the fold
bifurcations occur on the principal $n=0$ branch at $K_1=L_{00}=3.5$ and $K_1=M_{01}^+=5$.
Fig.~\ref{figA6align}(i) suggests that the folds on the other (unstable) branches of periodic
solutions all occur between the same $K_1$ values. Contrast this with Fig.~\ref{figA55align} where with $A=5.5$
there seems to be an alignment between the bifurcations on every second branch,
and Fig.~\ref{figA6align}(ii) where with $A$ equal to $4$ plus the golden ratio there does not appear to be
any alignment between the bifurcations on different branches.

To explain this alignment notice that in the singular limit $\epsilon\to0$ by Theorems~\ref{thmIbif}
and~\ref{thmIIbif}(i) for suitable integer value(s) of $m$ there are fold bifurcations on the $n$-th branch
at $K_1=L_{nm}$ and $K_1=M_{nm-1}^+$. But $L_{nm}$ defined by \eq{Lnm} and $M_{nm}^+$ defined as the larger
zero of $G_{nm}(K_1)$ (see \eq{Gnm}) both depend on $n$ and $m$ only through the common term $m-n(A-1)$.
Hence if $A=p/q$ is rational then defining
$$n_k=n_0+kq, \qquad m_k=m_0+k(p-q),\qquad k\in\N$$
we see that
$$m_k-n_k(A-1)=m_0-n_0(A-1),$$
and hence $L_{n_km_k}=L_{n_0m_0}$ and $M_{n_km_k}^+=M_{n_0m_0}^+$ for each integer $k$ and for each
$n_0=0,1,2,\ldots,q-1$. Hence these singular fold bifurcations align on every $q$-th branch
when $A=p/q$ is rational. Thus, when $A$ is integer these bifurcations align on
all the branches (eg $A=6$, see Fig.~\ref{figA6align}(i)), when $A=p/2$ the bifurcations align on
every second branch (eg $A=5.5$, see Fig.~\ref{figA55align}), and when $A$ is irrational there
is no alignment between the bifurcations (see Fig.~\ref{figA6align}(ii)).

Moving to the $\epsilon>0$ case, we see from the figures with $\epsilon=0.05$ that the fold bifurcations
which should align exactly in the limit as $\epsilon=0$, actually appear to occur within shrinking subintervals
of $[L_{n_0m_0},M_{n_0m_0-1}^+]$ as $n_k$ is increased, and for sufficiently large $n_k$ the fold bifurcations
disappear entirely. Although for each fixed $n_k$ the folds occur for all $\epsilon$ sufficiently
small and converge to $K_1=L_{n_0m_0}$ and $K_2=M_{n_0m_0-1}^+$ as $\epsilon\to0$ the convergence is clearly
not uniform, with smaller values of $\epsilon$ required to create the fold bifurcations for larger
values of $n_k$. This is not surprising, since the larger the value of $n_k$ the smaller the period
and amplitude of the $\epsilon=0$ singular solutions defined by Theorems~\ref{thm2delunimod}, \ref{thm2delI}
and~\ref{thm2delII}. But, when solving with $\epsilon>0$ the smaller amplitude solutions appear smoother
and more sinusoidal than the larger amplitude solutions and the fold bifurcations do not occur unless
$\epsilon$ is reduced sufficiently to resolve the sawteeth in the solution.

\section{Conclusions} \label{sec:conc}

Through Definitions~\ref{defadmiss} and~\ref{defsingsol} we have
introduced a new definition of singular solution via a double
parametrisation which allows us to define a continuous parametrisation even when the limiting profile is not continuous.
This reduces the problem of constructing singular solutions to a purely
algebraic problem.
For the DDE \eq{eps2del} with two state-dependent delays we constructed three different
solution profiles in Section~\ref{secsingsols} and in
Theorems~\ref{thm2delunimod}, \ref{thm2delI} and~\ref{thm2delII}
identified parameter constraints for these unimodal and type I and type II
bimodal singular solutions to exist. In Section~\ref{secbifs} we investigated the parameter constraints for the
singular solutions constructed in Section~\ref{secsingsols}, and treating $K_1$ as a
bifurcation parameter in Theorem~\ref{thmunilegs} identified intervals of $K_1$ for which unimodal singular solutions exist.
Theorem~\ref{thmIbif} identifies intervals on which type I bimodal solutions exist, and also a singular fold bifurcation where
the solution profile also transitions between unimodal and type I bimodal.
Theorem~\ref{thmIIbif} identifies intervals on which type II bimodal solutions exist,
and a point where the solution profile transitions between unimodal and type II bimodal, with or without a
singular fold bifurcation, and we hence identify a singular codimension-two bifurcation.

The results in Sections~\ref{secsingsols} and~\ref{secbifs} all follow from our definition of singular solution
following purely algebraic arguments. Although we do not prove analytically that the singularly perturbed DDE
\eq{eps2del} has corresponding periodic solutions for $0<\epsilon\ll1$, in Section~\ref{secnumerics} we
demonstrate numerically using DDEBiftool that the singular periodic solutions that we found do persist
for $\epsilon>0$. Moreover, we find that there is very good agreement between the parts of the bifurcation
diagram determined by the unimodal and bimodal singular solutions, and the numerically computed small $\epsilon$
branches and profiles. The $\epsilon>0$ computations also reveal intervals of bistability of unimodal periodic solutions
and unstable solutions with two, three and more local maxima per period.

In Section~\ref{secnumerics} we saw that for $0<\epsilon\ll1$ fold bifurcations occur close to
$K_1=L_{nm}$ and $K_1=M_{nm}^+$. In Section~\ref{sec:cusp}, we considered the codimension-two bifurcations
predicted by Theorem~\ref{thmIIbif}
where the fold at $K_1=L_{nm}$ vanishes and for $0<\epsilon\ll1$ found the predicted cusp bifurcations,
and associated stable bimodal periodic orbits (see Figs.~\ref{figA448amp} and \ref{figamp698})

In Section~\ref{sec:othersols}, Theorem~\ref{thmunilegs}(i) led us to find stable periodic orbits with
two local maxima per period when $\epsilon>0$ (Fig.~\ref{figamp76}).
A period doubling bifurcation also gives rise to stable periodic orbits with up to 5 local maxima per period
(see Fig.\ref{figamppdA15eps005}). We were also able to use our singular solution theory to predict the alignment of the fold
bifurcations on different solution branches.

In addition to the fold bifurcations associated with the transition between unimodal and bimodal solutions
predicted by Theorems~\ref{thmIbif} and~\ref{thmIIbif} we also found many examples of solutions of \eq{eps2del}
with three or more local maxima per period sometimes with fold bifurcations associated to the transitions between
such solutions. In contrast the one delay DDE \eq{eps1del} has only been seen to have
periodic orbits with one local maxima per period, and no secondary bifurcations on the branches of periodic orbits
\cite{DCDSA11}.


In conclusion, the state-dependent DDE \eq{eps2del} has very rich and interesting dynamics in the $\epsilon\to0$ singular
limit, and the concept of singular solution that we introduce in
Definitions~\ref{defadmiss} and~\ref{defsingsol} is a useful tool in the study of those dynamics. While we have not proved
rigourously that the singular solutions that we construct persist for $\epsilon>0$, we have shown numerically that they do,
and identified where bifurcations occur. A useful first step in proving convergence as $\epsilon\to0$ is to identify
what the singularly perturbed solutions should converge to. With this work that first step is resolved.

\section{Acknowledgments}

Tony Humphries thanks John Mallet-Paret and Roger Nussbaum for introducing him
to this problem and patiently explaining their results in the one delay case. He is also grateful to NSERC (Canada) for funding
through the Discovery Grants program.
Renato Calleja thanks NSERC and the Centre de recherches math\'ematiques, Montr\'eal for funding and to
FQRNT, Qu\'ebec for a PBEEE award.
Daniel Bernucci and Michael Snarski are grateful to NSERC for Undergraduate Student Research Awards.
Namdar Homayounfar thanks the Institut des Sciences Math\'ematiques, Montr\'eal for an
Undergraduate Summer Scholarship.

\end{document}